\documentclass{article}

\usepackage{arxiv}

\usepackage{fancyhdr,lipsum}
\makeatletter

\usepackage[cp1252]{inputenc}
\usepackage[T1]{fontenc}    
\usepackage{hyperref}       
\usepackage{url}            
\usepackage{booktabs}       
\usepackage{amsfonts}       
\usepackage{nicefrac}       
\usepackage{microtype}      
\usepackage{lipsum}
\usepackage{graphicx}

\usepackage[cp1252]{inputenc}
\usepackage{textcomp}
\usepackage{amsmath}
\usepackage{amssymb}
\usepackage{graphicx}
\usepackage{adjustbox}
\usepackage[abs]{overpic}
\usepackage{hyphenat}
\usepackage{inputenc}
\usepackage{amsthm}
\usepackage{booktabs}
\usepackage{bbold}
\usepackage{enumitem}
\setenumerate[1]{label = \arabic*.}
\setenumerate[2]{label*= \arabic*}
\setenumerate[3]{label*=.\arabic*}
\setenumerate[4]{label*=.\arabic*}
\setenumerate[4]{label*=.\arabic*}
\usepackage{xcolor}
\usepackage{multirow}
\usepackage[titletoc]{appendix}
\usepackage{anyfontsize}
\usepackage{framed}
\usepackage{footnote}
\usepackage{mathbbol}
\usepackage{subcaption}
\usepackage{placeins}

\usepackage{adjustbox}

\usepackage{algorithm}
\usepackage{algpseudocode}

\usepackage{url}
\usepackage{comment}
\usepackage{xcolor}

\usepackage{marginnote}

\usepackage{tabulary}
\usepackage{bbm} 
\newtheoremstyle{dotless}{}{}{\itshape}{}{\bfseries}{}{ }{}
\theoremstyle{dotless}

\makeatletter

\renewcommand{\mathbf}{\boldsymbol}

\graphicspath{{figures/}}

\usepackage{color}
\usepackage{soul}
\usepackage{amsthm}

\usepackage{units}

\newcommand{\ba}{\mbox{\boldmath{$a$}}}
\newcommand{\bA}{\mbox{\boldmath{$A$}}}
\newcommand{\bb}{\mbox{\boldmath{$b$}}}

\newcommand{\bd}{\mbox{\boldmath{$d$}}}

\newcommand{\bg}{\mbox{\boldmath{$g$}}}
\newcommand{\bI}{\mbox{\boldmath{$I$}}}
\newcommand{\bK}{\mbox{\boldmath{$K$}}}

\newcommand{\bn}{\mbox{\boldmath{$n$}}}
\newcommand{\bN}{\mbox{\boldmath{$N$}}}
\newcommand{\bp}{\mbox{\boldmath{$p$}}}

\newcommand{\bs}{\mbox{\boldmath{$s$}}}
\newcommand{\bt}{\mbox{\boldmath{$t$}}}

\newcommand{\bu}{\mbox{\boldmath{$u$}}}
\newcommand{\bU}{\mbox{\boldmath{$U$}}}

\newcommand{\bv}{\mbox{\boldmath{$v$}}}

\newcommand{\bx}{\mbox{\boldmath{$x$}}}

\newcommand{\bepsilon}{\mbox{\boldmath{$\varepsilon$}}}

\newcommand{\bsigma}{\mbox{\boldmath{$\sigma$}}}

\newcommand{\bzero}{\mbox{$\bf 0$}}

\newfont{\twelvemsb}{msbm10 at 11.6pt}

\renewcommand{\div}{\mathop{\rm div}\nolimits}

\newcommand{\tr}{\mathop{\rm tr}}

\title{On mesh coarsening procedures for the virtual element method}

\author{
	Daniel van Huyssteen \\
	Institute of Applied Mechanics, 
	Friedrich-Alexander-Universität Erlangen-Nürnberg
	Erlangen, 91058, Germany \\
	\texttt{daniel.van.huyssteen@fau.de} \\
	\And
	Felipe Lopez Rivarola \\
	Facultad de Ingeniería,
	Universidad de Buenos Aires,
	Buenos Aires, C1127AAR, Argentina,
	CONICET\\
	\And
	Guillermo Etse \\
	Facultad de Ingeniería,
	Universidad de Buenos Aires,
	Buenos Aires, C1127AAR, Argentina,
	CONICET\\
	\AND
	Paul Steinmann \\
	Institute of Applied Mechanics, 
	Friedrich-Alexander-Universität Erlangen-Nürnberg
	Erlangen, 91058, Germany \\
}

\begin{document}
	
	\maketitle
	
	
	\begin{abstract}
	In the context of adaptive remeshing, the virtual element method provides significant advantages over the finite element method. The attractive features of the virtual element method, such as the permission of arbitrary element geometries, and the seamless permission of `hanging´ nodes, have inspired many works concerning error estimation and adaptivity. However, these works have primarily focused on adaptive refinement techniques with little attention paid to adaptive coarsening (i.e. de-refinement) techniques that are required for the development of fully adaptive remeshing procedures.
	In this work novel indicators are proposed for the identification of patches/clusters of elements to be coarsened, along with a novel procedure to perform the coarsening. The indicators are computed over prospective patches of elements rather than on individual elements to identify the most suitable combinations of elements to coarsen. The coarsening procedure is robust and suitable for meshes of structured and unstructured/Voronoi elements. 
	Numerical results demonstrate the high degree of efficacy of the proposed coarsening procedures and sensible mesh evolution during the coarsening process. It is demonstrated that critical mesh geometries, such as non-convex corners and holes, are preserved during coarsening, and that meshes remain fine in regions of interest to engineers, such as near singularities.
	\end{abstract}
	
	
	\newpage
	\noindent
	
	\section{Introduction}
	\label{S:Introduction}
	Adaptive remeshing techniques are a critical tool in engineering analysis that facilitate automatic and localized manipulation of a mesh to improve its approximation properties for a given problem. In contrast to conventional uniform remeshing, adaptive remeshing allows degrees of freedom to concentrate in, and be removed from, the most and least critical regions of a problem domain respectively. Thus, improving the accuracy of numerical simulation methods while reducing their computational load.
	A typical adaptive procedure has the well-known ${\text{`Solve}\rightarrow\text{Estimate}\rightarrow\text{Mark}\rightarrow\text{Remesh´}}$ structure. These four main steps involve; generating an approximate numerical solution to a problem using some pre-existing mesh, using the approximate solution to estimate the error, marking/flagging elements to be refined or coarsened/de-refined based on the estimated error, and finally creating an updated mesh. These steps are then performed iteratively until some user-defined termination criteria is met.
	
	In the context of the finite element method (FEM) adaptive remeshing techniques are already well-established. There exists a wide range of approaches to \textit{a-posteriori} error estimation and a variety of tools/packages for the creation of updated meshes. 
	Some of the most widely used approaches to \textit{a-posteriori} error estimation include residual-based \cite{babuvvska1978error,Babuska1979} and recovery-based error estimators \cite{Zienkiewicz1987,Zienkiewicz1992,Zienkiewicz1992a,Zienkiewicz1995}. 
	In general, there is greater focus on adaptive refinement techniques over adaptive coarsening/de-refinement techniques. However, there exists a range of open-source software libraries capable of performing both adaptive refinement and coarsening of finite element meshes \cite{Arndt2022,Kirk2006}. 
	Performing localized adaptation of finite element meshes is non-trivial as significant manipulation of not only the elements being adapted but also of the surrounding elements is required to preserve the method's conformity. In general, coarsening of finite element meshes is more complex than refinement. As such, most coarsening processes performed using finite elements only reverse previously performed refinement to return a mesh to its initially coarser state, see for example \cite{Park2012}. The complexities involved in coarsening finite element meshes mean that true coarsening procedures are rarely used in practical applications. 
	
	The introduction of the virtual element method (VEM) gave rise to many new opportunities in the context of adaptive remeshing.
	The VEM is an extension of the FEM that permits arbitrary polygonal and polyhedral element geometries in two- and three-dimensions respectively \cite{VEIGA2012,Veiga2014}. A feature of the VEM of particular interest in the context of adaptive remeshing is the permission of arbitrarily many nodes along an element's edge. That is, nodes that would be considered `hanging´ in a finite element context are trivially incorporated into the virtual element formulation \cite{VEMContactWriggers2016,Wriggers2019}. 
	The geometric robustness of the VEM has been demonstrated with the method exhibiting optimal convergence behaviour in cases of challenging, including strongly non-convex, element geometries \cite{Sorgente2021,Sorgente2022,Huyssteen2020,Huyssteen2021}. Additionally, in cases of distorted, and possibly stretched, element geometries that could arise during adaptive remeshing (particularly during anisotropic remeshing) the VEM stabilization term can be easily tuned to improve the accuracy of the method \cite{ReddyvanHuyssteen2021,DvH_BDR_MeshQuality}.
	Furthermore, the robustness of the VEM under challenging numerical conditions, such as near-incompressibility and near inextensibility, is increasingly well reported \cite{WriggersIsotropic2017,Wriggers2023,Tang2020,Reddy2019,Huyssteen2020,Huyssteen2021}. 
	Due to its geometric flexibility, as well as its geometric and numerical robustness, the VEM is an obvious candidate for problems involving adaptive remeshing. 
	However, in contrast to the FEM, the VEM basis/shape functions are not explicitly defined over an element domain. Thus, the typical error estimators used in a FEM context usually cannot be trivially applied in a VEM setting. Additionally, the freedom of element geometry permitted by the VEM requires the development of more versatile mesh refinement and coarsening techniques than those used in FEM applications.
	
	Adaptive refinement is currently a very popular topic in the VEM context with a rapidly growing literature focusing on residual-based \cite{Cangiani2017,Veiga2015a,Berrone2017,Mora2017} and recovery-based \cite{Chi2019,Guo2019,Wei2023,NguyenThanh2018} \textit{a-posteriori} error estimation. 
	Furthermore, several approaches for localized refinement of the unstructured polygonal element geometries permitted by the VEM have been presented \cite{NguyenThanh2018,Huyssteen2022,Berrone2021}.
	A comparatively unexplored topic is the development of techniques for adaptive coarsening (i.e. de-refinement) of virtual element meshes. Since a fully adaptive remeshing procedure comprises both refinement and coarsening capabilities, the development of coarsening techniques for the VEM is of great significance. 
	Furthermore, since the VEM permits arbitrarily many nodes along an element's edge, preserving the conformity of the method during coarsening is trivial. Thus, the VEM is particularly well-suited to problems involving truly adaptive coarsening processes. 
	
	To the best of the authors' knowledge, to date only one contribution exists focusing on adaptive coarsening in a virtual element context \cite{Choi2020}. The work \cite{Choi2020} uses the element-level displacement gradient error to identify elements to group. Specifically, neighbouring elements whose error estimators fall below a certain threshold are collected into groups. The mesh is then coarsened by combing the groups of elements into larger individual elements using a simple edge straightening procedure. The results presented are promising and of great interest. Specifically, the method yields coarser meshes in less critical portions of a problem domain while retaining finer discretizations near singularities. Additionally, the coarsening approach improves the efficiency of the VEM solution. That is, the amount of error per degree of freedom exhibited by the coarsened meshes is lower than that of uniform meshes.
	
	The importance of developing coarsening procedures for the VEM to be used in fully adaptive remeshing procedures, and the promising results of the coarsening approach presented in \cite{Choi2020}, strongly motivate further development and investigation into adaptive coarsening procedures for the VEM.
	In addition to its importance for fully adaptive remeshing procedures there are numerous other benefits of, even standalone, adaptive coarsening procedures. For example, a simulation using an adaptively coarsened mesh should have a similar level of accuracy to, and should execute significantly faster than, a uniform mesh simulation. This improved efficiency can be exploited to speed up problems that require multiple solutions or solution steps such as; sensitivity analysis, Monte Carlo simulations, solution procedures in non-linear analysis, dynamic problems, and fracture analysis. The improved efficiency offered by adaptively coarsened meshes can, thus, reduce both the computational time and energy consumption of the aforementioned problem types.
	
	In this work a novel patch-based approach to computing coarsening indicators for meshes of virtual elements is presented. The indicators are computed over prospective patches of elements to coarsen in contrast to the conventional element-level indicators used in adaptive refinement procedures. The patch-based approach is motivated by seeking to identify the most suitable combinations of elements to coarsen, and can even be used to approximately predict the element-level error after coarsening. This prediction could be used to determine if the coarsening of a particular patch might result in unsatisfactorily high local or global error after coarsening and should not be performed. Two approaches to the computation of the coarsening indicator over a patch are proposed and are based on the displacement field and approximate energy error norm computed using a recovery procedure. 
	Furthermore, a novel approach to determining the geometry of the coarsened element groupings is presented. Specifically, the updated/coarsened element geometry is created by constructing a convex hull around an element patch. To this end, a novel edge straightening procedure is presented and used in creating the geometry of the convex hull. The robustness of the proposed edge straightening procedure is demonstrated through its application to complex groups of edges.
	The performance of the proposed indicators and coarsening procedures is investigated through a range of benchmark problems of varying complexity.
	For each problem the mesh evolution during coarsening is analysed along with the behaviour in the $\mathcal{H}^{1}$ error norm for structured and unstructured/Voronoi meshes with a range of initial uniform mesh densities.
	Finally, to facilitate a thorough investigation of the proposed coarsening procedures, and critically for brevity, it is chosen in this work to analyse only adaptive coarsening procedures. The analysis of fully adaptive remeshing procedures thus represents a future contribution.
	
	The structure of the rest of this work is as follows. The governing equations of linear elasticity are set out in Section~\ref{sec:GovEq}. This is followed in Section~\ref{sec:VEM} by a description of the first-order virtual element method. 
	The procedures used to generate and coarsen meshes are presented in Section~\ref{sec:MeshGenerationAndCoarsening}. This is followed, in Section~\ref{sec:MeshCoarseningIndicators}, by a description of the various mesh coarsening indicators along with the procedure used to identify patches of elements qualifying for coarsening.
	Section~\ref{sec:Results} comprises a set of numerical results through which the performance of the various coarsening procedures is evaluated.
	Finally, the work concludes in Section~\ref{sec:Conclusion} with a discussion of the results.
	
	\section{Governing equations of linear elasticity} 
	\label{sec:GovEq}
\begin{figure}[ht!]
		\centering
		\includegraphics[width=0.25\textwidth]{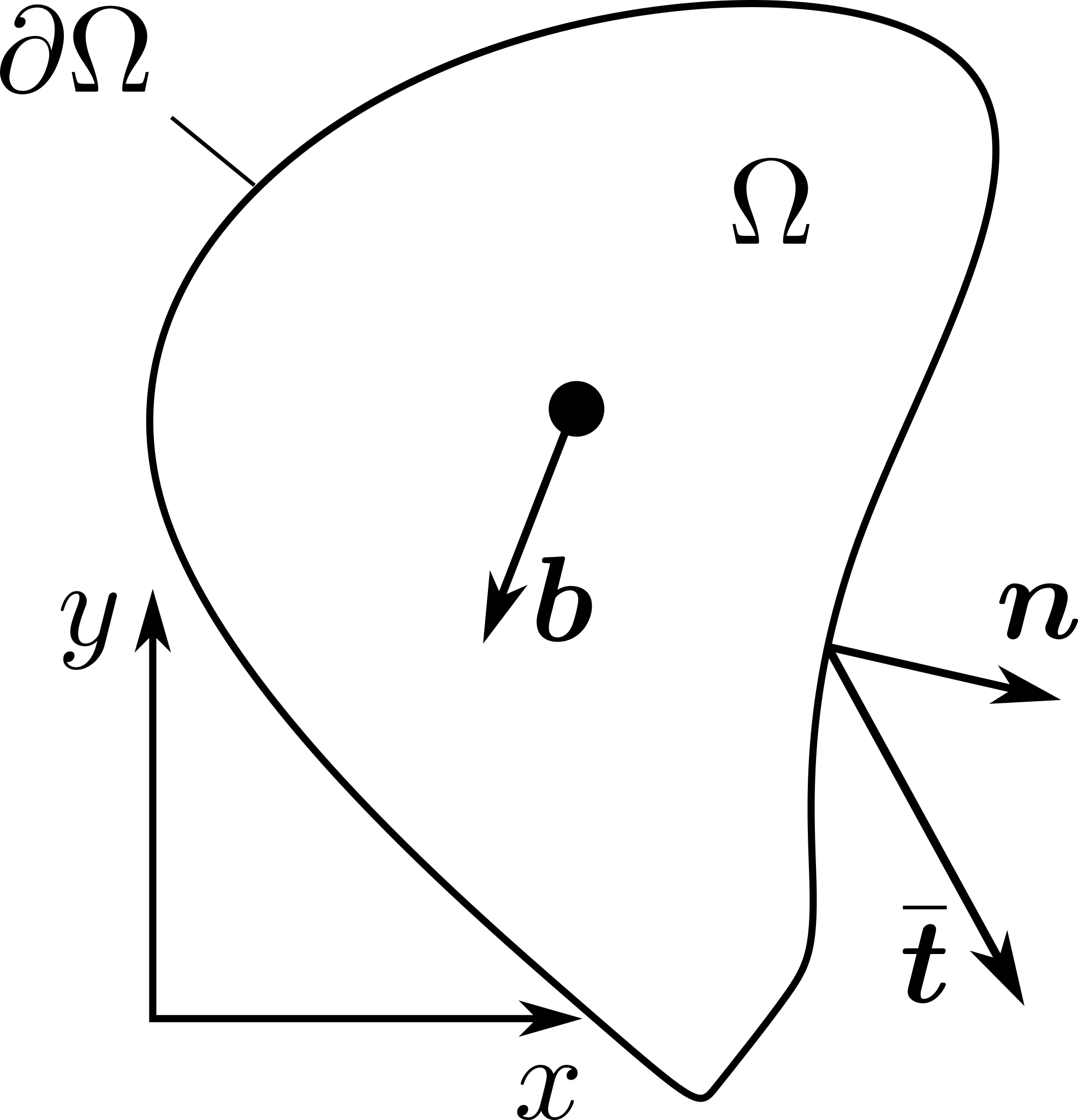}
		\caption{Arbitrary elastic body subject to body force and traction.
			\label{fig:ElasticBody}}
	\end{figure}

	Consider an arbitrary elastic body occupying a plane, bounded, domain ${\Omega \subset \mathbb{R}^{2} }$ subject to a traction ${\bar{\bt}}$ and body force ${\bb}$ (see Figure~\ref{fig:ElasticBody}).
	The boundary ${\partial \Omega}$ has an outward facing normal denoted by $\bn$ and comprises a non-trivial Dirichlet part $\Gamma_{D}$ and a Neumann part $\Gamma_{N}$ such that ${\Gamma_{D} \cap \Gamma_{N} = \emptyset}$ and ${\overline{\Gamma_{D} \cup \Gamma_{N}}=\partial \Omega}$.
	
	\FloatBarrier
	\FloatBarrier
	
	In this work small displacements are assumed and the strain-displacement is relation given by 
	\begin{equation}
		\bepsilon \left( \bu \right) = \frac{1}{2} \left[\nabla\,\bu + \left[ \nabla\,\bu \right]^{T} \right] \, . \label{eqn:DisplacementStrain}
	\end{equation}
	Here the displacement is denoted by $\bu$, ${\bepsilon}$ is the symmetric infinitesimal strain tensor and ${\nabla \left( \bullet \right) = \frac{\partial \left( \bullet \right)_{i} }{\partial \, x_{j}} \, \boldsymbol{e}_{i} \otimes \boldsymbol{e}_{j} }$ is the gradient of a vector quantity. Additionally, linear elasticity is assumed and the stress-strain relation is given by
	\begin{equation}
		\bsigma = \mathbb{C} : \bepsilon \,. \label{eqn:StressStrain1}
	\end{equation}
	Here, ${\bsigma}$ is the Cauchy stress tensor and ${\mathbb{C}}$ is a fourth-order constitutive tensor. For a linear elastic and isotropic material (\ref{eqn:StressStrain1}) is given by
	\begin{eqnarray}
		\bsigma = \lambda \tr \left( \bepsilon \right)  \bI + 2\mu \, \bepsilon \, , \label{eqn:StressStrain2}
	\end{eqnarray}
	where ${\tr \left( \bullet \right)}$ denotes the trace, $\bI$ is the second-order identity tensor, and $\lambda$ and $\mu$ are the well-known Lam\'{e} parameters.
	
	For equilibrium it is required that 
	\begin{equation}
		\div \, \bsigma + \bb = \bzero \, , \label{eqn:Equilibrium}
	\end{equation}
	where ${\div \left( \bullet \right) = \frac{\partial \left( \bullet \right)_{ij} }{\partial \, x_{j}} \boldsymbol{e}_{i} }$ is the divergence of a tensor quantity.
	The Dirichlet and Neumann boundary conditions are given by 
	\begin{align}
		&\bu = \bg \quad \text{on } \Gamma_{D} \, , \text{ and} \label{eqn:DirichletBC} \\
		&\bsigma \cdot \bn = \bar{\bt} \quad \text{on } \Gamma_{N} \, , \label{eqn:NeumannBC}
	\end{align}
	respectively, with $\bg$ and $\bar{\bt}$ denoting prescribed displacements and tractions respectively.
	Equations (\ref{eqn:StressStrain2})-(\ref{eqn:NeumannBC}), together with the displacement-strain relationship (\ref{eqn:DisplacementStrain}), constitute the boundary-value problem for a linear elastic isotropic body.
	
	\subsection{Weak form}
	\label{subsec:WeakForm}
	The space of square-integrable functions on $\Omega$ is hereinafter denoted by ${\mathcal{L}^{2}\left(\Omega\right)}$. The Sobolev space of functions that, together with their first derivatives, are square-integrable on $\Omega$ is hereinafter denoted by ${\mathcal{H}^{1}\left(\Omega\right)}$. Additionally, the function space $\mathcal{V}$ is introduced and defined such that
	\begin{align}
		\mathcal{V} = \left[ \mathcal{H}^{1}_{D} \left(\Omega\right) \right]^{d} 
		=
		\left\{ \bv \, | \, v_{i} \in \mathcal{H}^{1}\left(\Omega\right), \, \bv = \boldsymbol{0} \,\, \text{on} \,\, \Gamma_{D}  \right\} \, 
	\end{align}
	where ${d=2}$ is the dimension.
	Furthermore, the function ${\bu_{g}\in \left[ \mathcal{H}^{1} \left(\Omega\right) \right]^{d} }$ is introduced satisfying (\ref{eqn:DirichletBC}) such that ${\bu_{g}|_{\Gamma_{D}}=\bg}$.
	
	The bilinear form ${a\left(\cdot,\cdot\right)}$, where ${a : \left[ \mathcal{H}^{1} \left(\Omega\right) \right]^{d}  \times \left[ \mathcal{H}^{1} \left(\Omega\right) \right]^{d} \rightarrow \mathbb{R}}$, and the linear functional ${\ell \left(\cdot\right)}$, where ${\ell : \left[ \mathcal{H}^{1} \left(\Omega\right) \right]^{d} \rightarrow \mathbb{R}}$, are defined respectively by
	\begin{equation}
		a\left(\bu, \, \bv \right) = \int_{\Omega} \bsigma \left( \bu \right) : \bepsilon \left( \bv \right) \, dx \, , \label{eqn:BilinearForm}
	\end{equation}
	and 
	\begin{equation}
		\ell \left( \bv \right) = \int_{\Omega} \bb \cdot \bv \, dx + \int_{\Gamma_{N}} \bar{\bt} \cdot \bv \, ds - a\left(\bu_{g},\, \bv \right) \, . \label{eqn:LinearFuntional} 
	\end{equation}
	
	The weak form of the problem is then: given ${\bb \in \left[ \mathcal{L}^{2}\left(\Omega\right) \right]^{d} }$ and ${ \bar{\bt} \in \left[ \mathcal{L}^{2}\left(\Gamma_{N}\right) \right]^{d} }$, find ${\bU \in \left[ \mathcal{H}^{1}\left(\Omega\right) \right]^{d} }$ such that 
	\begin{equation}
		\bU = \bu + \bu_{g} \, , \quad \bu \in \mathcal{V} \, ,
	\end{equation}
	and
	\begin{equation}
		a\left(\bu , \, \bv \right) = \ell \left( \bv \right) \, , \quad \forall \bv \in \mathcal{V} \, . \label{eqn:BilinearFormFinal}
	\end{equation}
	
	\section{The virtual element method}
	\label{sec:VEM}
	
\begin{figure}[ht!]
		\centering
		\includegraphics[width=0.33\textwidth]{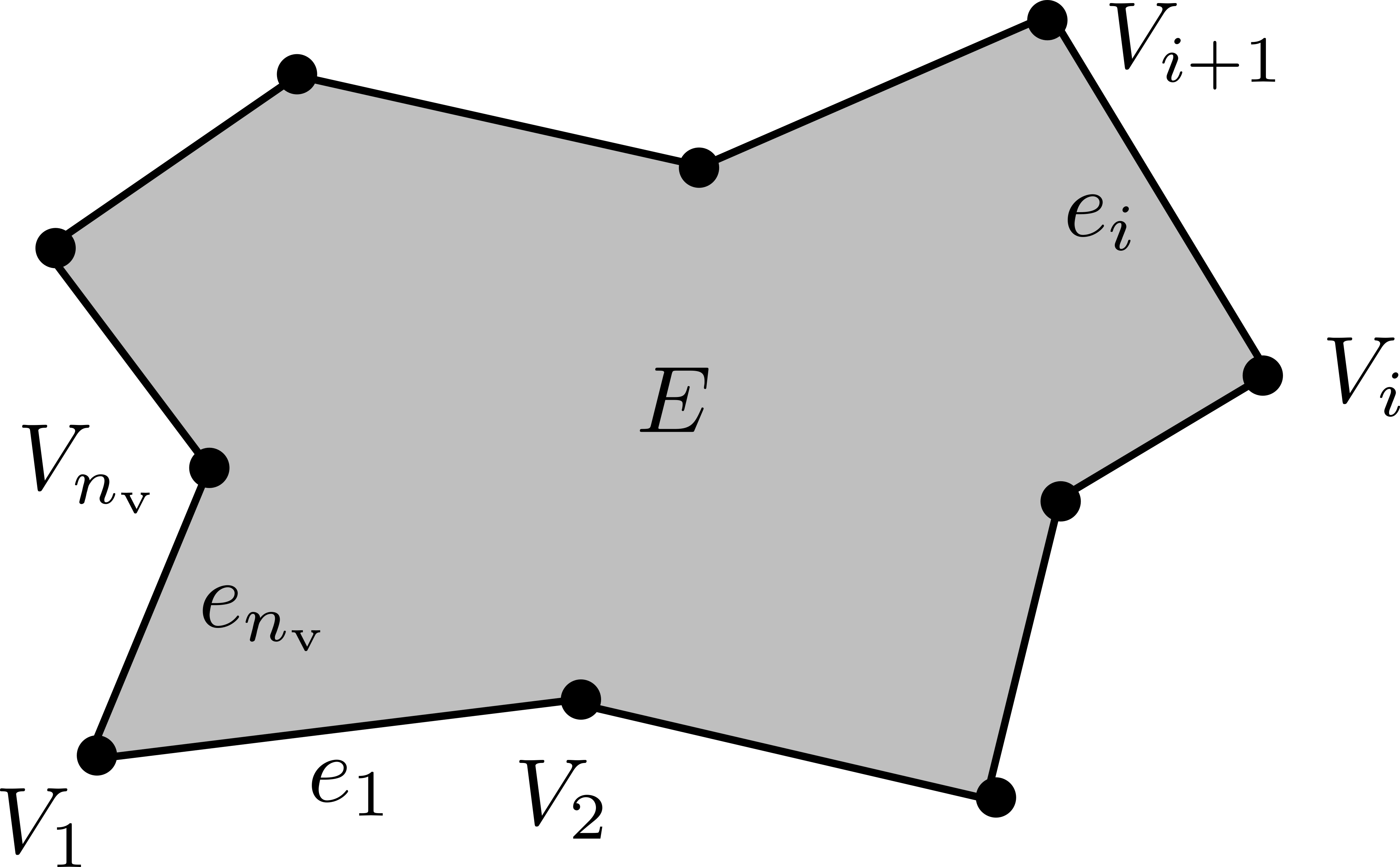}
		\caption{Sample virtual element.
			\label{fig:SampleElement}}
	\end{figure}

	The domain $\Omega$ is partitioned into a mesh of non-overlapping arbitrary polygonal elements\footnote{If $\Omega$ is not polygonal the mesh will be an approximation of $\Omega$} $E$ with $\overline{\cup E}=\overline{\Omega}$. Here $E$ denotes the element domain and $\partial E$ its boundary, with ${\overline{(\, \bullet \,)}}$ denoting the closure of a set.
	An example of a typical first-order element is depicted in Figure~\ref{fig:SampleElement} with edge $e_{i}$ connecting vertices $V_{i}$ and $V_{i+1}$. Here ${i=1,\dots,n_{\rm v}}$ with $n_{\rm v}$ denoting the total number of element vertices.
	
	\FloatBarrier
	\FloatBarrier
	A conforming approximation of order $k$ is constructed in a space ${\mathcal{V}^{h} \subset \mathcal{V}}$ where $\mathcal{V}^{h}$ is built-up element-wise and comprises vector valued functions $\bv_{h}$.
	The functions $\bv_{h}$ are those that are $\mathcal{C}^{0}$ continuous on the domain $\Omega$, are polynomials of degree ${\leq \, k}$ on element edges, and whose strain gradient divergence is a polynomial of degree ${\leq \, k-2}$ on an element (see \cite{Artioli2017}). 
	For the most general case of an approximation of arbitrary order $k$ the space $\mathcal{V}^{h}|_{E}$ is defined as 
	\begin{equation}
		\mathcal{V}^{h}|_{E} = \left\{ \bv_{h} \in \mathcal{V} \, | \, \bv_{h} \in \left[ \mathcal{C}^{0}(E) \right]^{2} \, , \, \nabla^{2} \, \bv_{h}  \in \mathcal{P}_{k-2} \text{ on } E \, , \, \bv_{h}|_{e} \in \mathcal{P}_{k}(e)  \right\} \,. \label{eqn:ArbitraryVEMSpace}
	\end{equation}
	Here ${\mathcal{P}_{k}(X)}$ is the space of polynomials of degree ${\leq \, k}$ on the set ${X \, \subset \, \mathbb{R}^{d} }$ with ${d=1,\,2}$ and ${\nabla^{2}=\nabla\cdot\nabla}$ is the Laplacian operator. 
	In this work a first-order, i.e. ${k=1}$, approximation is considered, thus (\ref{eqn:ArbitraryVEMSpace}) simplifies to 
	\begin{equation}
		\mathcal{V}^{h}|_{E} = \left\{ \bv_{h} \in \mathcal{V} \, | \, \bv_{h} \in \left[ \mathcal{C}^{0}(E) \right]^{2} \, , \, \nabla^{2} \, \bv_{h}  = \boldsymbol{0} \text{ on } E \, , \, \bv_{h}|_{e} \in \mathcal{P}_{1}(e)  \right\} \,. 
	\end{equation}
	
	All computations will be performed on element edges and it is convenient to write, for element $E$,
	\begin{equation}
		\bv_{h}|_{\partial E} = \bN \cdot \bd^{E} \,. \label{eqn:DisplacementTrace}
	\end{equation}
	Here, $\bN$ is a matrix of standard linear Lagrangian basis functions and $\bd^{E}$ is a ${{2}n_{\rm v} \times 1}$ vector of the degrees of freedom associated with $E$. The virtual basis functions are not known, nor required to be known on $E$; their traces, however, are known and are simple Lagrangian functions.
	
	The virtual element projection for a first-order formulation ${\Pi \, : \, \mathcal{V}^{h}|_{E} \rightarrow \mathcal{P}_{0}(E)  }$ is required to satisfy
	\begin{equation}
		\int_{E} \Pi \, \bv_{h} \cdot \bepsilon\left( \bp \right) \, dx = \int_{E} \bepsilon\left(\bv_{h}\right) \cdot \bepsilon\left( \bp \right) \, dx \quad \forall \bp \in \mathcal{P}_{1} \,, \label{eqn:Projection}
	\end{equation}
	where ${\Pi \, \bv_{h}}$ represents the ${\mathcal{L}^{2}}$ projection of the symmetric gradient of ${\bv_{h}}$ onto constants \cite{Artioli2017}. Since the projection is constant at element-level, after applying integration by parts to (\ref{eqn:Projection}), and considering (\ref{eqn:DisplacementTrace}), the components of the projection can be computed as
	\begin{align}
		\left(\Pi\,\bv_{h}\right)_{ij} &= \frac{1}{2}\frac{1}{|E|}  \sum_{e\in\partial E} \int_{e} \left[ N_{iA} \, d_{A}^{E} \, n_{j} + N_{jA} \, d_{A}^{E} \, n_{i}\right] ds \,, \label{eqn:ComputeProjection}
	\end{align}
	where summation is implied over repeated indices.
	
	The virtual element approximation of the bilinear form (\ref{eqn:BilinearForm}) is constructed by writing 
	\begin{align}
		a^{E}\left(\bu,\,\bv\right) :&= a\left(\bu,\,\bv\right)|_{E} 
		= \int_{E} \bepsilon\left(\bv_{h}\right) : \left[ \mathbb{C} : \bepsilon\left(\bu_{h}\right) \right] dx \, , \label{eqn:ElementBilinearForm}
	\end{align}
	where ${a^{E}\left(\cdot,\cdot\right)}$ is the contribution of element $E$ to the bilinear form ${a\left(\cdot,\cdot\right)}$. Consideration of (\ref{eqn:ComputeProjection}) allows (\ref{eqn:ElementBilinearForm}) to be written as (see \cite{Reddy2019})
	\begin{align}
		a^{E}\left(\bu_{h},\,\bv_{h}\right) 
		&= \underbrace{\int_{E} \Pi\,\bv_{h} : \left[ \mathbb{C} : \Pi\,\bu_{h} \right] dx }_{\text{Consistency term}} 
		+ \underbrace{\int_{E} \left[ \bepsilon\left( \bv_{h} \right) : \left[ \mathbb{C} : \bepsilon \left( \bu_{h} \right) \right] - \Pi\,\bv_{h} : \left[ \mathbb{C} : \Pi \, \bu_{h} \right] \right] dx }_{\text{Stabilization term}} \,, \label{eqn:ExpandedBilinearForm}
	\end{align}
	where the remainder term is discretized by means of a stabilization.
	
	\subsection{The consistency term} 
	\label{subsec:ConsistencyTerm}
	
	The projection (\ref{eqn:ComputeProjection}), and thus the consistency term, can be computed exactly yielding
	\begin{equation}
		a_{\rm c}^{E}\left(\bu_{h},\,\bv_{h}\right) \, = \, \int_{E} \Pi\,\bv_{h} : \left[ \mathbb{C} : \Pi\,\bu_{h} \right] dx \, = \, \widehat{\bd}^{E} \cdot \left[ \bK_{\rm c}^{E} \cdot \bd^{E} \right] \,.
	\end{equation}
	
	Here $\bK_{\rm c}^{E}$ is the consistency part of the stiffness matrix of element $E$ with ${\widehat{\bd}^{E}}$ and ${\bd^{E}}$ the degrees of freedom of $\bv_{h}$ and $\bu_{h}$ respectively that are associated with element $E$.
	
	\subsection{The stabilization term} 
	\label{subsec:Stab}
	The remainder term cannot be computed exactly and is approximated by means of a discrete stabilization term \cite{Gain2014,Veiga2015}.
	The approximation employed in this work is motivated by seeking to approximate the difference between the element degrees of freedom $\bd^{E}$ and the nodal values of a linear function that is closest to $\bd^{E}$ in some way (see \cite{Artioli2017,Reddy2019}). 
	The nodal values of the linear function are given by
	\begin{equation}
		\widetilde{\bd} = \boldsymbol{\mathcal{D}} \cdot \bs \,. \label{eqn:LinearApprox}
	\end{equation}
	Here $\bs$ is a vector of the degrees of freedom of the linear function and $\boldsymbol{\mathcal{D}}$ is a matrix relating $\widetilde{\bd}$ to $\bs$ with respect to a scaled monomial basis. For the full expression of $\boldsymbol{\mathcal{D}}$ see \cite{Artioli2017,Reddy2019}.
	After some manipulation (see, again, \cite{Reddy2019}) the stabilization term of the bilinear form can be approximated as
	\begin{equation}
		a_{\text{stab}}^{E}\left(\bu_{h},\,\bv_{h}\right) \, = \,
		\int_{E} \left[ \bepsilon\left( \bv_{h} \right) : \left[ \mathbb{C} : \bepsilon \left( \bu_{h} \right) \right] - \Pi\,\bv_{h} : \left[ \mathbb{C} : \Pi \, \bu_{h} \right]  \right] dx \, \approx \, \widehat{\bd}^{E} \cdot \bK_{\rm s}^{E} \cdot \bd^{E} \, ,
	\end{equation}
	where $\bK_{\rm s}^{E}$ is the stabilization part of the stiffness matrix of element $E$ and is defined as
	\begin{equation}
		\bK_{\rm s}^{E} = \mu \left[ \bI - \boldsymbol{\mathcal{D}} \cdot \left[ \boldsymbol{\mathcal{D}}^{T} \cdot \boldsymbol{\mathcal{D}}\right]^{-1} \cdot  \boldsymbol{\mathcal{D}}^{T} \right] \, .
	\end{equation}
	The total element stiffness matrix ${\bK^{E}}$ is then computed as the sum of the consistency and stabilization matrices.
	
	\section{Mesh generation and  mesh coarsening} 
	\label{sec:MeshGenerationAndCoarsening}
	In this section the procedures used to generate meshes and coarsen patches of elements are described.
	
	\subsection{Mesh generation} 
	\label{subsec:MeshGeneration}
\begin{figure}[ht!]
		\centering
		\includegraphics[width=0.95\textwidth]{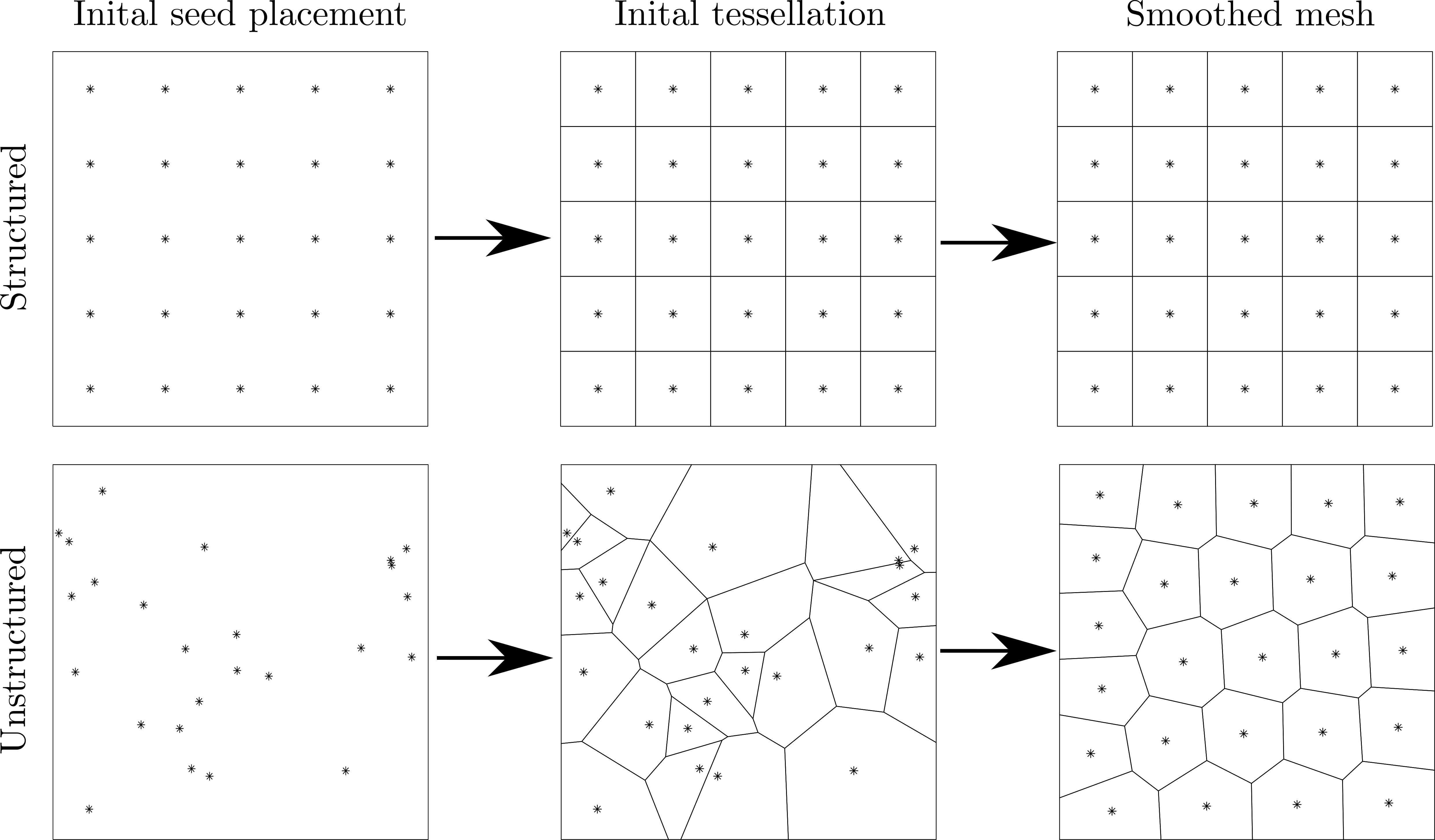}
		\caption{Mesh generation procedure for structured and unstructured/Voronoi meshes.
			\label{fig:MeshGeneration}}
	\end{figure}

	The mesh generation procedure used in this work is identical to that described in \cite{Huyssteen2022}.
	All meshes are created by Voronoi tessellation of a set of seed points. Seed points will be generated in both structured and unstructured sets to create structured and unstructured meshes respectively. 
	In the case of structured meshes seeds points are placed to form a structured grid, while in the case of unstructured/Voronoi meshes seeds are placed arbitrarily within the problem domain. Hereinafter the terms 'unstructured' and 'Voronoi' meshes will be used interchangeably to refer to meshes created from arbitrarily placed seed points.
	An initial Voronoi tessellation of the seed points is created using PolyMesher \cite{PolyMesher}. Then, a smoothing algorithm in PolyMesher is used to iteratively modify the locations of the seed points to create a mesh in which all elements have approximately equal areas.
	The mesh generation procedure is illustrated in Figure~\ref{fig:MeshGeneration} where the top and bottom rows depict the generation of structured and unstructured/Voronoi meshes respectively.  
	
	\FloatBarrier
	\FloatBarrier
	
	\subsection{Mesh coarsening} 
	\label{subsec:MeshCoarsening}
	The patches of elements qualifying for coarsening are identified using the procedure described in the next section. Once a patch of elements has been marked for coarsening, the coarsening process is performed by grouping the elements into a single larger element formed by a convex hull.
	
	\subsubsection{Overview of coarsening procedure} 
	\label{subsubsec:MeshCoarseningOverview}
\begin{figure}[ht!]
		\centering
		\includegraphics[width=0.825\textwidth]{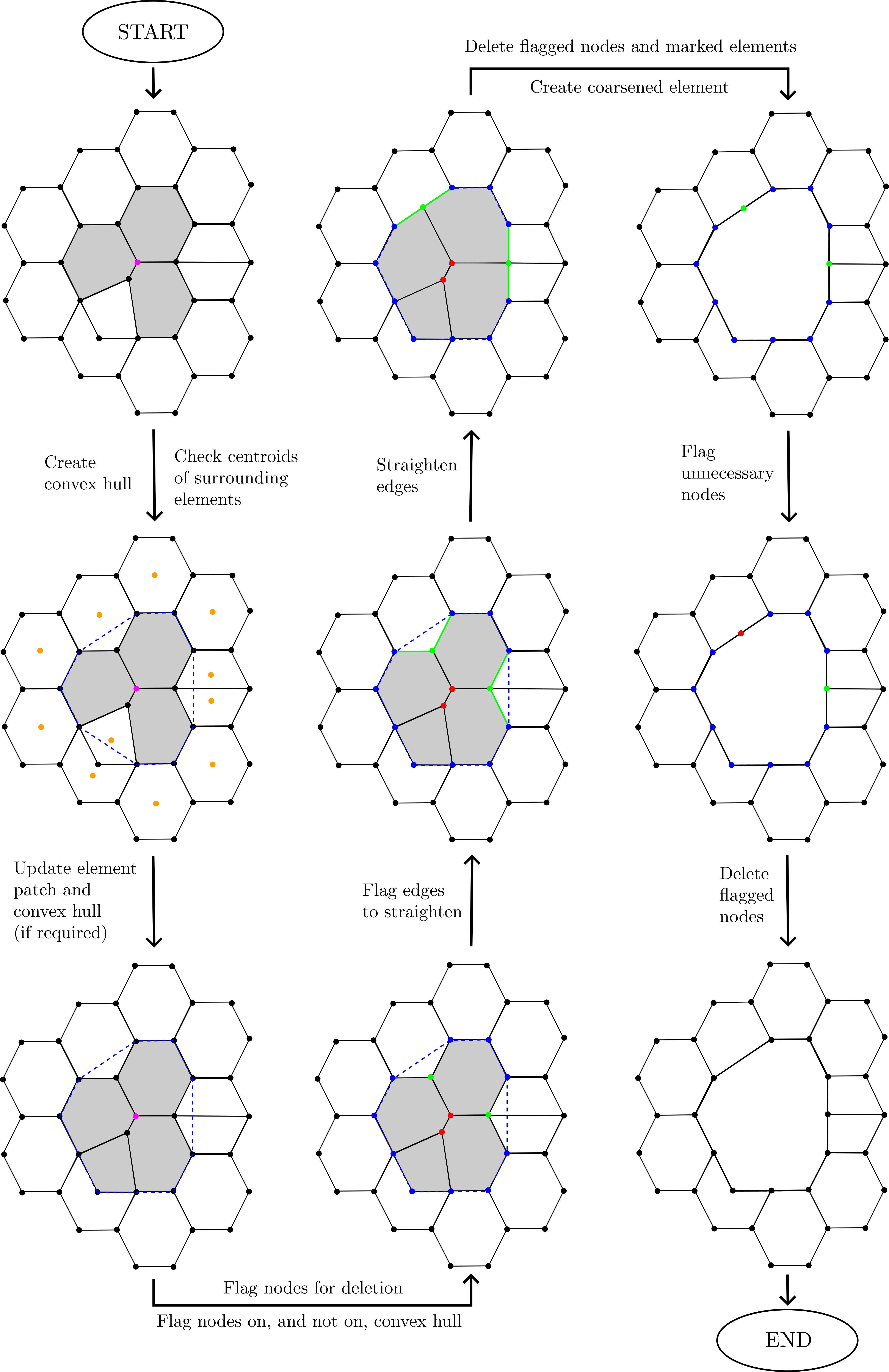}
		\caption{Coarsening procedure overview.
			\label{fig:MeshCoarseningOverview}}
	\end{figure}

	The same coarsening procedure is used for patches of both structured and unstructured elements. An overview of the mesh coarsening procedure is illustrated in Figure~\ref{fig:MeshCoarseningOverview} for the more general case of unstructured elements. Here the patch of elements marked for coarsening is indicated in grey and the surrounding elements are indicated in white. An element is considered to be surrounding the patch of marked elements if it shares any nodes with any of the marked elements.
	
	The first step in the coarsening procedure involves creating a convex hull surrounding the patch of marked elements. The convex hull is indicated by a blue dashed line. Thereafter, the surrounding elements are checked to determine if their centroids (indicated in orange) lie within the convex hull. If any of the centroids of the surrounding elements do lie within the convex hull, those elements are added to the patch of elements to coarsen and an updated convex hull around the new patch of elements is computed.
	
	The next phase of the coarsening procedure involves categorising the nodes in the element patch. Firstly, nodes that will not form part of the coarsened element geometry are identified and flagged for deletion. These nodes are those which are completely surrounded by marked elements. That is, the relevant node does not lie on the domain boundary, and every element connected to the node has been marked for coarsening (i.e. is indicated in grey). Nodes that meet these criteria are indicated in red. Secondly, the remaining nodes, i.e. non-red nodes, are divided into two groups; those that lie on the convex hull (indicated in blue), and those that do not lie on the convex hull (indicated in green). 
	
	The green nodes are used to identify which edges of the element patch require straightening to align with the convex hull. An edge is flagged for straightening if any of its nodes do not lie on the convex hull, i.e. if any of its nodes have been marked as green. The edges flagged for straightening are indicated in green. These edges are then straightened using the procedure described in the next section.
	
	After the edges have been straightened the nodes flagged for deletion (red nodes) are removed, the elements marked for coarsening (grey elements) are deleted, and a new element is created from the remaining (blue and green) nodes connected in the conventional counter-clockwise sequence.
	
	All remaining nodes (blue and green nodes) are checked to determine if they are necessary. That is, if removing them would alter the geometry or connectivity of the newly created element, or any of the surrounding elements. If all the element edges connected to a particular node are co-linear, then that node is considered as unnecessary. The unnecessary nodes are then flagged for deletion (indicated in red). Finally, the unnecessary nodes are removed.
	
	\FloatBarrier
	\FloatBarrier
	
	\subsubsection{Edge straightening procedure} 
	\label{subsubsec:EdgeStraighteningProcedure}
\begin{figure}[ht!]
		\centering
		\includegraphics[width=0.95\textwidth]{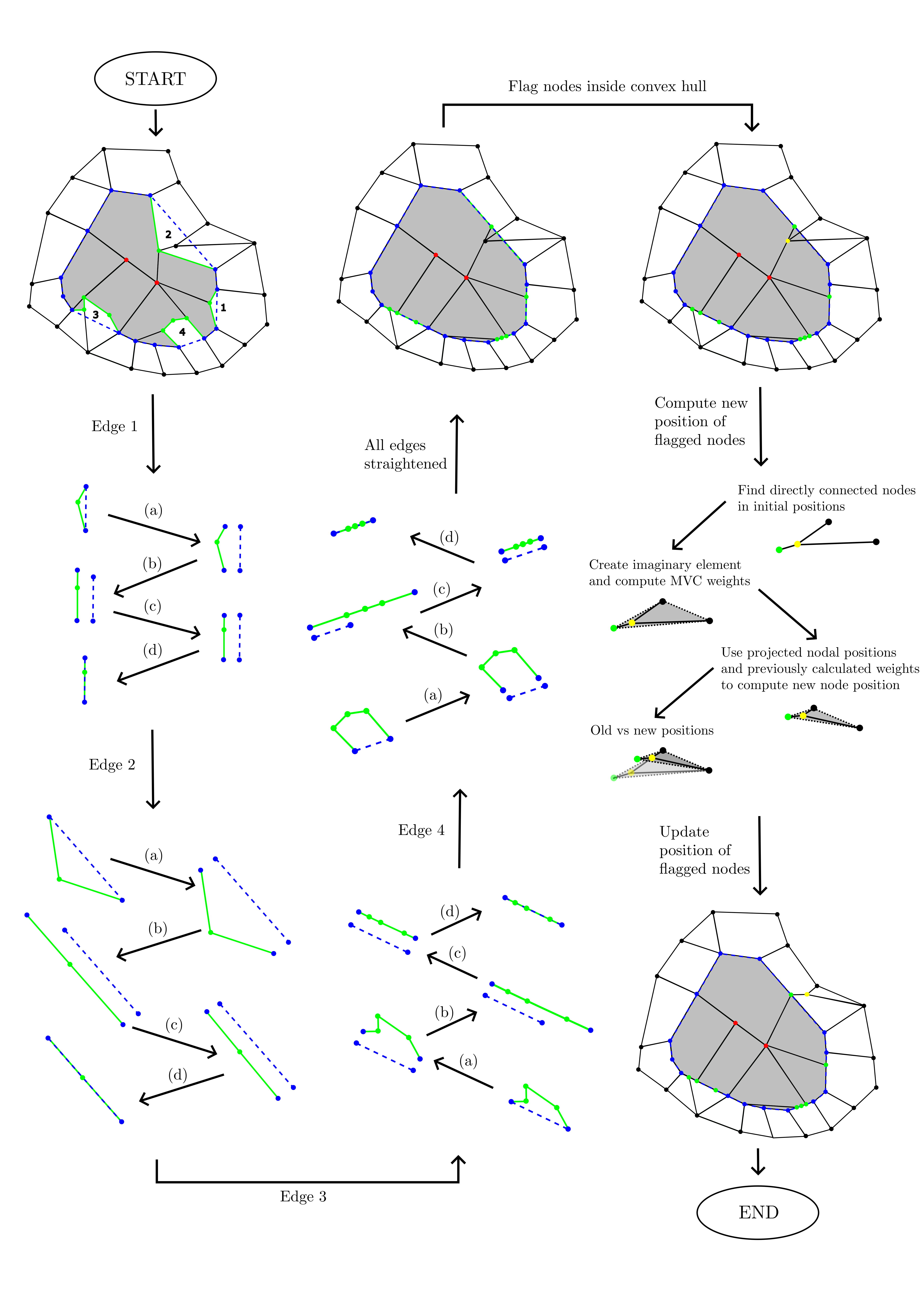}
		\vspace*{-12mm}
		\caption{Edge straightening procedure.
			\label{fig:EdgeStraighteningProcedure}}
	\end{figure}

	The edge straightening procedure is illustrated in Figure~\ref{fig:EdgeStraighteningProcedure} where the first step depicted corresponds to the step in Figure~\ref{fig:MeshCoarseningOverview} before the edge straightening is performed\footnote{It is noted that the element geometries presented in Figure~\ref{fig:EdgeStraighteningProcedure} are not from any simulation or actual results. The exaggerated geometric features are presented for illustrative purposes to demonstrate challenging scenarios that could arise during a coarsening procedure.}. 
	The edge straightening procedure is performed by grouping consecutive edges that have been marked for straightening and are not separated by any (blue) nodes lying on the convex hull. Examples of these groups of edges are depicted in Figure~\ref{fig:EdgeStraighteningProcedure} and are numbered for clarity.
	
	A group of edges to be straightened comprises some number of (green) edges that have been identified for straightening, and two (blue) nodes lying on the convex hull that form a line segment indicated by a blue dashed line. The straightening is performed by first straightening the green edges to form a line that is parallel to the line segment formed by the two blue nodes. This straightened line, together with the nodes along it, is then scaled to the size of the blue dashed line segment. This procedure creates a straightened edge along which the green nodes have the same relative spacing as they did in the unstraightened configuration. Additionally, the procedure prevents flattening of sections of the surrounding elements and has been found to be robust in the presence of challenging and highly non-convex geometries.
	
	The edge straightening procedure is illustrated in Figure~\ref{fig:EdgeStraighteningProcedure} for each of the groups of edges to be straightened. First, the edges to be straightened, the nodes on the convex hull, and the line segment connecting these nodes are depicted. In step (a) the edges to be straightened and the line segment are separated for illustrative purposes. In step (b) the (green) edges are straightened to form a line parallel to the (blue dashed) line segment. In step (c) the straightened line, together with the nodes along it, is scaled to the size of the (blue dashed) line segment. Finally, step (d) depicts the newly straightened and scaled edge overlaying the (blue dashed) line segment and, thus, depicts the result of the straightening procedure for a particular edge.
	
	After all of the edge groups have been straightened, all of the nodes of the surrounding (white) elements are checked to determine if they lie inside the convex hull, these nodes are indicated in yellow. If a node does lie inside the convex hull it is projected to an updated position using mean value coordinates (MVC) \cite{Floater2003}. This is performed by finding all other nodes that are directly connected to, i.e. share an edge with, the yellow node and using both their initial/unprojected and projected locations. Using the initial locations a fictitious element is imagined comprising the directly connected nodes. Then, relative weights are computed using MVC for each of the directly connected nodes at the location of the yellow node. This corresponds to evaluating the MVC weight functions for each directly connected node at the position of the yellow node. Finally, the new/projected location of the yellow node is computed as the weighted sum of the updated/projected locations of the directly connected nodes using the previously calculated MVC weights. This step is employed to prevent `tangling´ of elements during the coarsening procedure.
	It is noted that other approaches to relocating the (yellow) nodes trapped inside the convex hull are possible. For example, the yellow nodes could be projected onto the convex hull using a minimum distance projection. Alternatively, the geometry of the (white) surrounding elements could be trimmed so that there is no overlap with the convex hull and the yellow nodes could then be deleted. The MVC-based approach proposed here is chosen because it does not require the creation or deletion of any nodes, and does not alter the nodal connectivity of any elements.
	
	\FloatBarrier
	\FloatBarrier
	
	\section{Mesh coarsening indicators} 
	\label{sec:MeshCoarseningIndicators}
	In this section the proposed mesh coarsening indicators are presented along with the procedure used to identify the patches of elements qualifying for coarsening. 
	
	Since the coarsening of a mesh involves combining groups/patches of elements into a single larger element, it is chosen to construct and compute the coarsening indicators over patches of elements. Patches of elements are identified as all of the elements connected to a specific node. The computed coarsening indicator values are then assigned to the relevant node but reflect the behaviour over the element patch. 
	
\begin{figure}[ht!]
		\centering
		\includegraphics[width=0.75\textwidth]{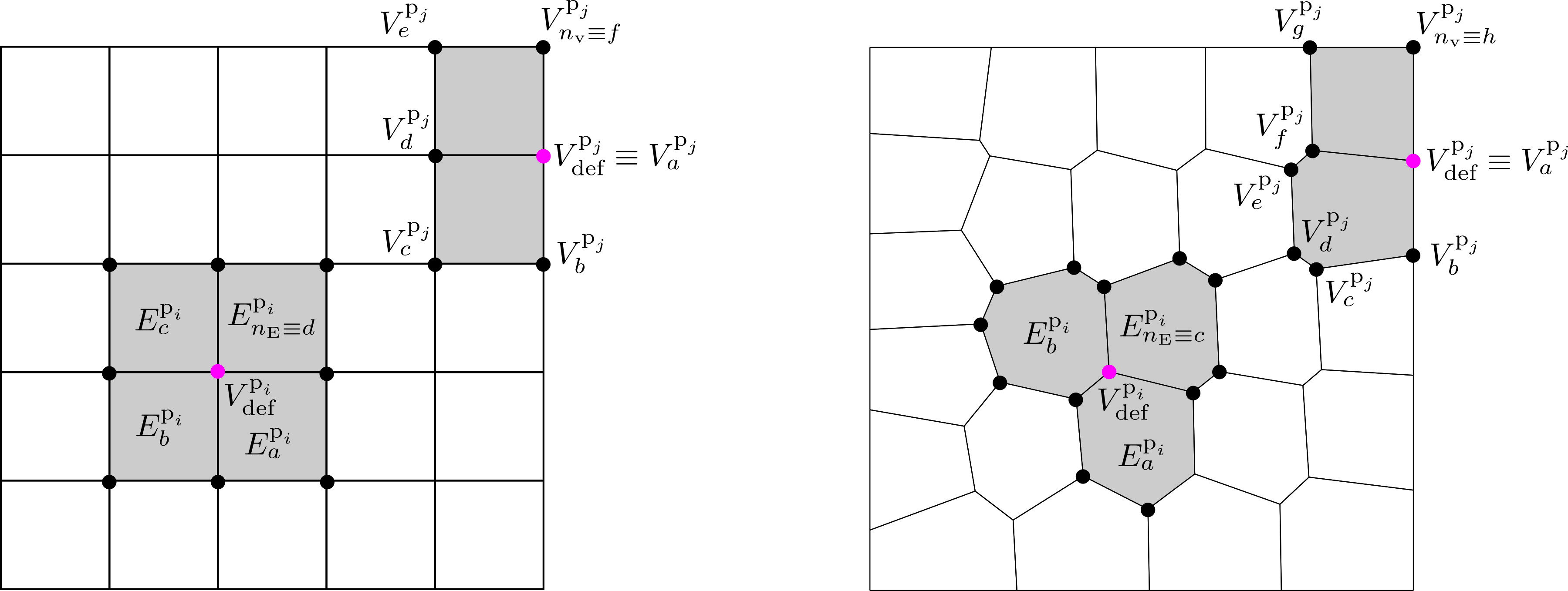}
		\caption{Examples of element patches.
			\label{fig:ExamplesOfElementPatches}}
	\end{figure}

	Examples of element patches are depicted in Figure~\ref{fig:ExamplesOfElementPatches} where the node defining patch $i$, i.e. $V_{\text{def}}^{\text{p}_{i}}$, is indicated in purple. All of the elements in the patch are indicated in grey and are labelled $E_{a}^{\text{p}_{i}},\,E_{b}^{\text{p}_{i}},\, \dots E_{n_{E}}^{\text{p}_{i}} $ where $n_{E}$ is the number of elements in the patch and the superscript $\text{p}_{i}$ indicates their association with the $i$-th patch. Furthermore, the nodes associated with the $j$-th patch are indicated on the figure and are labelled as $V_{a}^{\text{p}_{j}},\,V_{b}^{\text{p}_{j}},\, \dots V_{n_{\text{v}}}^{\text{p}_{j}} $ where $n_{\text{v}}$ is the number of nodes associated with the patch.
	
	\FloatBarrier
	\FloatBarrier
	
	\subsection{Displacement-based indicator} 
	\label{subsec:DispBasedIndicator}
	Similar to the displacement-based refinement indicator proposed in \cite{vanHuyssteen2022,Huyssteen2022}, the displacement-based coarsening indicator is motivated by seeking to quantify the deviation from coplanar of the nodal values of the displacement $\bu_{h}$ over a patch of elements. To compute the indicator for patch $i$ a least-squares best fit linear approximation of the displacement field over the patch ${\bu_{\text{p}_{i}}}$ is computed. The displacement field ${\bu_{\text{p}_{i}}}$ is computed component-wise with component $k$ described over patch $i$ as 
	\begin{equation}
		u_{k}^{\text{p}_{i}} = \bp \left(x, \, y\right) \, \ba_{k} =
		\begin{bmatrix}
			1 & x & y
		\end{bmatrix}
		\begin{bmatrix}
			a_{k}^{1} \\ a_{k}^{2} \\ a_{k}^{3} \, .
		\end{bmatrix}
	\end{equation}
	Here $\ba_{k}$ are the degrees of freedom of $u_{k}^{\text{p}_{i}}$ and are computed as
	\begin{equation}
		\ba_{k} = \bA^{-1} \bb_{k}
	\end{equation}
	where
	\begin{equation}
		\bA = \sum_{m=1}^{n_{\text{v}}} \bp\left(x_m^{\text{p}_{i}},\,y_m^{\text{p}_{i}}\right)^{T} \bp\left(x_m^{\text{p}_{i}},\,y_m^{\text{p}_{i}}\right)
		\quad \text{and} \quad 
		\bb_{k} = \sum_{m=1}^{n_{\text{v}}} \bp\left(x_m^{\text{p}_{i}},\,y_m^{\text{p}_{i}}\right)^{T} u^{h}_{k}\left(x_m^{\text{p}_{i}},\,y_m^{\text{p}_{i}}\right)
	\end{equation}
	respectively. Here $x_m^{\text{p}_{i}}$ and $y_m^{\text{p}_{i}}$ are the coordinates of the $m$-th node associated with patch $i$, and $u^{h}_{k}$ is the displacement degree of freedom at $\bx_m^{\text{p}_{i}}$.
	The displacement-based coarsening indicator on patch $i$, denoted by ${\Upsilon_{\text{DB}}^{i}}$, is then defined as the $\mathcal{L}^2$ deviation of the nodal values of the displacement $\bu_{h}$ from the least-squares displacement ${\bu_{\text{p}_{i}}}$ and is computed as
	\begin{equation}
		\Upsilon_{\text{DB}}^{i} = \left[ \sum_{j=1}^{n_{\text{v}}} \left[\bu_{\text{p}_{i}} \left(\bx^{\text{p}_{i}}_{j}\right) - \bu_{h} \left(\bx^{\text{p}_{i}}_{j}\right) \right] 
		\cdot  
		\left[\bu_{\text{p}_{i}} \left(\bx^{\text{p}_{i}}_{j}\right) - \bu_{h} \left(\bx^{\text{p}_{i}}_{j}\right) \right]
		\right]^{0.5} \, .
	\end{equation}
	
	\subsection{Energy error-based indicator} 
	\label{subsec:EnergyErrorBasedIndicator}
	The energy error-based error indicator is inspired by the well-known $Z^{2}$ error estimator originally presented in \cite{Zienkiewicz1987}, and the energy error approximation technique presented in \cite{NguyenThanh2018} for virtual elements. The energy error-based indicator is motivated by seeking to predict how much coarsening a particular patch of elements would increase the local and global approximations of the energy error. Then, those patches that are identified to increase the error the least are the most suitable for coarsening. 
	
	The error in the $\mathcal{H}^{1}$ semi-norm, i.e. the energy error norm, is defined as 
	\begin{equation}
		e_{\mathcal{H}^{1}} =  
		\left[\frac{1}{2} \int_{\Omega} 
		\left[ \bsigma^{ex} - \bsigma^{h} \right]^{T}
		\mathbb{D}^{-1}  
		\left[ \bsigma^{ex} - \bsigma^{h} \right] d\Omega
		\right]^{0.5} \, , \label{eqn:EnergyExact}
	\end{equation}
	where $\bsigma^{ex}$ is the exact/analytical stress solution and $\mathbb{D}$ is the constitutive matrix. In practical applications the exact stress is typically unknown and is replaced with an approximation $\bsigma^{\ast}$. The stress $\bsigma^{\ast}$ is usually computed as a higher-order approximation of the element stresses described in terms the displacement basis functions. However, in VEM applications the displacement basis functions are not explicitly defined over the entire domain and it is common to use node-based error approximations. Thus, it is sufficient to compute $\bsigma^{\ast}$ at the nodal positions. This is most easily done by computing a super-convergent stress at each node using a patch-based recovery technique based on super-convergent sampling points (see \cite{NguyenThanh2018}). 
	
\begin{figure}[ht!]
		\centering
		\includegraphics[width=0.8\textwidth]{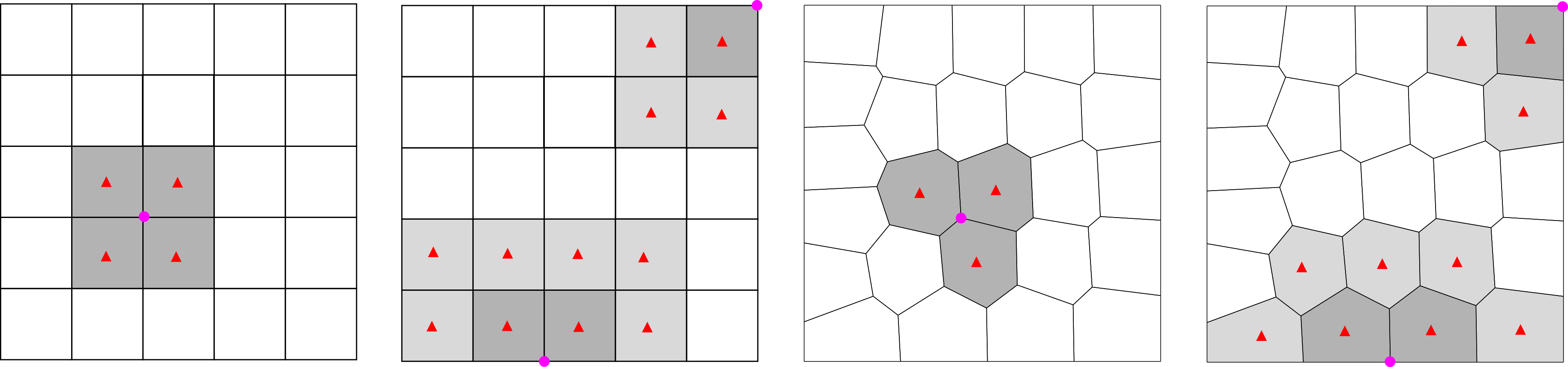}
		\caption{Super-convergent sampling points.
			\label{fig:SuperConvergentStress}}
	\end{figure}

	In this work a low-order VEM is considered where the approximation of the stress field is constant at element-level, i.e. piece-wise constant. Thus, the higher-order stress field approximation $\bsigma^{\ast}$ must be piece-wise linear and is computed at each node via a least-squares linear best fit over a patch of elements. The super-convergent stress at a node is computed by considering the patch of elements connected to the node. The location of the centroids of the elements in the patch are treated as the super-convergent sampling points and the element-level stresses are assigned as the degrees of freedom of the sampling points. Since a linear best fit is required, at least three sampling points are needed in order to determine a unique fit. Thus, in cases where a node is connected to less than three elements the patch is enlarged to increase the number of sampling points. Specifically, the patch is enlarged to include elements that are connected to any of the elements in the original patch. For clarity, a few examples of element patches and sampling points are depicted in Figure~\ref{fig:SuperConvergentStress}. Here, the node at which the super-convergent stress is to be computed is indicated in purple, the elements in the patch connected to the node are indicated in dark grey, and (if applicable) the elements included in the enlarged patch are indicated in light grey. Additionally, the locations of the sampling points are indicated as red triangles.
	
	\FloatBarrier
	\FloatBarrier 
	
	The super-convergent stress component ${\sigma^{\ast}_{i}}$ computed over a specific patch is given by
	\begin{equation}
		\sigma^{\ast}_{i} = \bp \left(x, \, y\right) \, \ba_{i} =
		\begin{bmatrix}
			1 & x & y
		\end{bmatrix}
		\begin{bmatrix}
			a_{i}^{1} \\ a_{i}^{2} \\ a_{i}^{3}
		\end{bmatrix}
	\end{equation}
	where $\ba_{i}$ are the degrees of freedom of the super-convergent stress component. The degrees of freedom are computed as 
	\begin{equation}
		\ba_{i} = \bA^{-1} \bb_{i}
	\end{equation}
	where
	\begin{equation}
		\bA = \sum_{k=1}^{n_{\text{sp}}} \bp\left(x_k,\,y_k\right)^{T} \bp\left(x_k,\,y_k\right)
		\quad \text{and} \quad 
		\bb_{i} = \sum_{k=1}^{n_{\text{sp}}} \bp\left(x_k,\,y_k\right)^{T} \sigma^{h}_{i}\left(x_k,\,y_k\right)
	\end{equation}
	respectively. Here $n_{\text{sp}}$ is the number of sampling points, $x_k$ and $y_k$ are the coordinates of the sampling points, and $\sigma^{h}_{i}$ is the stress component at the sampling point (computed via (\ref{eqn:ComputeProjection})).
	
	Using $\bsigma^{\ast}$ a node-based approximation of the energy norm can be computed as
	\begin{equation}
		e_{\mathcal{H}^{1}} \approx  
		\left[ \frac{1}{2} 
		\sum_{i=1}^{n_{\text{el}}} 
		\frac{|E_{i}|}{n_{\text{v}}^{i}} \sum_{j=1}^{n_{\text{v}}^{i}} 
		\left[
		\left[ \bsigma^{\ast}\left(\bx_{j}\right) - \bsigma^{h}\left(\bx_{j}\right) \right]^{T} \mathbb{D}^{-1}    
		\left[ \bsigma^{\ast}\left(\bx_{j}\right) - \bsigma^{h}\left(\bx_{j}\right) \right]
		\right] \right]^{0.5} \, ,
		\label{eqn:EnergyApprox}
	\end{equation}
	where $n_{\text{el}}$ is the number of elements in the domain.
	Using this node-based approach to error approximation, the `predicted´ energy error after coarsening a particular patch of elements, i.e. the energy error-based coarsening indicator, is approximated over patch $i$ as 
	\begin{equation}
		\Upsilon_{\text{EB}}^{i} =
		\left[ \frac{1}{2} 
		\frac{|E_{\text{p}_i}|}{n_{\text{v}}} \sum_{j=1}^{n_{\text{v}}} 
		\left[
		\left[ \bsigma^{\ast}\left(\bx_{j}\right) - \bar{\bsigma}^{h}_{\text{p}_i}\left(\bx_{j}\right) \right]^{T} \mathbb{D}^{-1}    
		\left[ \bsigma^{\ast}\left(\bx_{j}\right) - \bar{\bsigma}^{h}_{\text{p}_i}\left(\bx_{j}\right) \right]
		\right] \right]^{0.5} \, .
		\label{eqn:EnergyErrorIndicator}
	\end{equation}
	Here $|E_{\text{p}_i}|$ denotes the area of patch $i$. Additionally, $\bar{\bsigma}^{h}_{\text{p}_i}$ denotes the `predicted´ stress over the coarsened patch computed as the average of the element stresses on the patch.
	
	The `predicted´ energy error is not only useful for adaptive coarsening procedures but offers an interesting advantage to fully adaptive remeshing procedures too. For example, during an adaptive remeshing procedure a user selects a global and/or local error criterion or target. The `predicted´ energy error can then be used to determine whether the coarsening of a particular patch might result in unsatisfactorily high local and/or global errors after coarsening that would not meet the global and/or local error criterion and thus should not be performed.
	
	\subsection{Selecting element patches to coarsen} 
	\label{subsec:SelectElementPatchedForCoarsening}
	The procedure used to select patches of elements to coarsen comprises two steps. First, all of the element patches that are eligible for coarsening are identified. Then, from the set of eligible patches, those whose coarsening indicator value falls below a certain threshold are selected for coarsening.
	
	\subsubsection{Identifying eligible element patches} 
	\label{subsubsec:IdentifyingEligiblePatches}
	For a patch of elements to be eligible for coarsening the geometry of the element created by the coarsening must not modify the overall geometry of the problem domain. In the case of problems with convex domains this step is trivial. However, if a problem domain is non-convex, or contains a hole, care must be taken to preserve the domain geometry. 
	
\begin{figure}[ht!]
		\centering
		\includegraphics[width=0.75\textwidth]{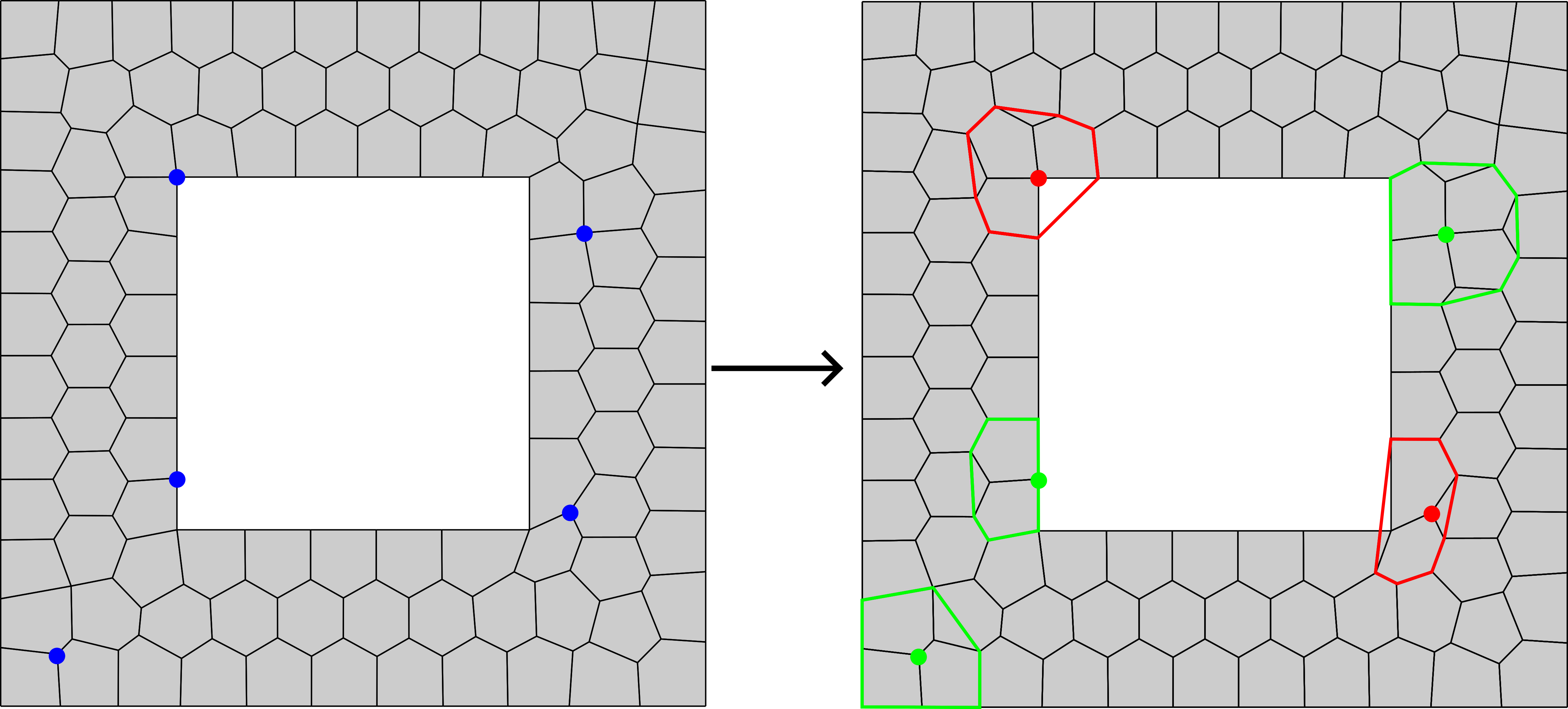}
		\caption{Examples of eligible and ineligible element patches.
			\label{fig:EligibleElementPatches}}
	\end{figure}

	To determine the eligibility of a patch all of the nodes associated with the patch are considered. The patch's nodes are checked to see if they lie on a domain boundary or corner. If no nodes lie on a boundary or corner then the patch is eligible for coarsening. If any nodes do lie on a boundary or corner then a convex hull is created from the patch's nodes. The positions of the boundary and corner nodes are then checked to see if they are coincident with the boundary of the convex hull. If all boundary and corner nodes are coincident with the boundary of the convex hull then the patch is eligible for coarsening. This process is exemplified in Figure~\ref{fig:EligibleElementPatches} where several nodes that define element patches are considered and are indicated in blue. The resulting convex hulls are indicated in green if the patch is eligible for coarsening and in red if the patch is not eligible.
	
	\FloatBarrier
	\FloatBarrier 
	
	After the set of eligible patches has been determined it must be resolved to eliminate overlapping patches. This is done by first sorting the nodes defining the patches in ascending order based on their coarsening indicator values. Then, all the nodes associated with the first patch are identified. The rest of sorted node list is then checked to see if any of the remaining nodes in the list are present in the current node patch. Any of the remaining nodes that are present in the current node patch are removed from the node list. This process is then repeated iteratively until the end of the list of eligible nodes is reached. This procedure is detailed in Algorithm~\ref{alg:RemoveOverlap}.
	
	\FloatBarrier
	\begin{algorithm}
		\caption{Remove overlapping element patches}\label{alg:RemoveOverlap}
		\begin{algorithmic}
			\State $EligibleNodes \gets \text{sort}(EligbleNodes,\, \text{ascending})$ 
			\State $nNodes \gets \text{length}(EligibleNodes)$
			\For{$i \in (1,nNodes)$}
			\State $ThisDefiningNode \gets EligibleNodes(i)$
			\State $PatchNodes \gets \text{GetPatchNodes}(ThisDefiningNode)$
			\For{$j \in (i,nNodes)$}
			\State $NodeToCheck \gets EligibleNodes(j)$
			\If{$NodeToCheck$ is member of $PatchNodes$}
			\State Mark $NodeToCheck$ for deletion
			\EndIf
			\EndFor
			\State Delete marked nodes from $EligibleNodes$
			\State $nNodes \gets \text{length}(EligibleNodes)$
			\EndFor
		\end{algorithmic}
	\end{algorithm}
	\FloatBarrier

	\subsubsection{Marking element patches} 
	\label{subsubsec:MarkingElementPatches}
	The procedure for identifying patches of elements to coarsen is similar to that presented in \cite{vanHuyssteen2022}.
	A coarsening threshold percentage ${T=X\%}$ is introduced from which an allowable threshold value ${T_{\text{val}}}$ is determined using the list of resolved eligible nodes.
	The node ${X\%}$ of the way down the resolved node list is found and the value of its coarsening indicator is set as ${T_{\text{val}}}$.
	Then any node (and associated element patch) whose coarsening indicator value is less than or equal to ${T_{\text{val}}}$ is marked for coarsening.
	
	\section{Numerical Results}
	\label{sec:Results}
	
	In this section numerical results are presented for a range of example problems of varying complexity to demonstrate the efficacy of the proposed coarsening procedures.
	The efficacy is evaluated in the $\mathcal{H}^{1}$ error norm defined by 
	\begin{align}
		||\widetilde{\bu}-\bu_{h}||_{1}=&\left[ \int_{\Omega}\left[ |\widetilde{\bu}-\bu_{h}|^{2} + |\nabla \widetilde{\bu} - \nabla \bu_{h}|^{2} \right] \, d\Omega \right]^{0.5} \, , \label{eqn:H1Error}
	\end{align}
	in which integration of $\bu_{h}$ is required over the domain. Since in the case of VEM formulations $\bu_{h}$ is only known on element boundaries a node-based approximation of (\ref{eqn:H1Error}) is used and is computed as
	\begin{align}
		\begin{split}
			||\widetilde{\bu}-\bu_{h}||_{1} \approx & \left[ \sum_{i=1}^{n_{el}}\frac{|E_{i}|}{n_{\rm v}^{i}} \sum_{j=1}^{n_{\rm v}} \Bigl[ \left[ \widetilde{\bu}(\bx_{j}) - \bu_{h}^{i}(\bx_{j}) \right] \cdot \left[ \widetilde{\bu}(\bx_{j}) -  \bu_{h}^{i}(\bx_{j}) \right] \Bigr. \right.  \\
			&\Biggl. \Bigl. + \left[ \nabla \widetilde{\bu}(\bx_{j}) - \Pi \bu_{h}^{i}(\bx_{j}) \right] : \left[ \nabla \widetilde{\bu}(\bx_{j}) - \Pi \bu_{h}^{i}(\bx_{j}) \right] \Bigr]  \Biggr]^{0.5} \, .
		\end{split} 
	\end{align}
	Here $\widetilde{\bu}$ is a reference solution generated using an overkill mesh of biquadratic finite elements, the location of the $j$-th vertex is denoted by $\bx_{j}$, and ${\Pi \bu_{h}}$ is the gradient of $\bu_{h}$ computed via the projection operator (see (\ref{eqn:ComputeProjection})).
	
	In the examples that follow the material is isotropic with a Young's modulus of ${E=1~\rm{Pa}}$, a Poisson's ratio of ${\nu=0.3}$, and where the shear modulus is computed as ${\mu = E /2 \left[1+\nu\right]}$. 
	
	\subsection{Punch}
	\label{subsec:Punch}
\begin{figure}[ht!]
		\centering
		\begin{subfigure}[t]{0.45\textwidth}
			\centering
			\includegraphics[width=0.95\textwidth]{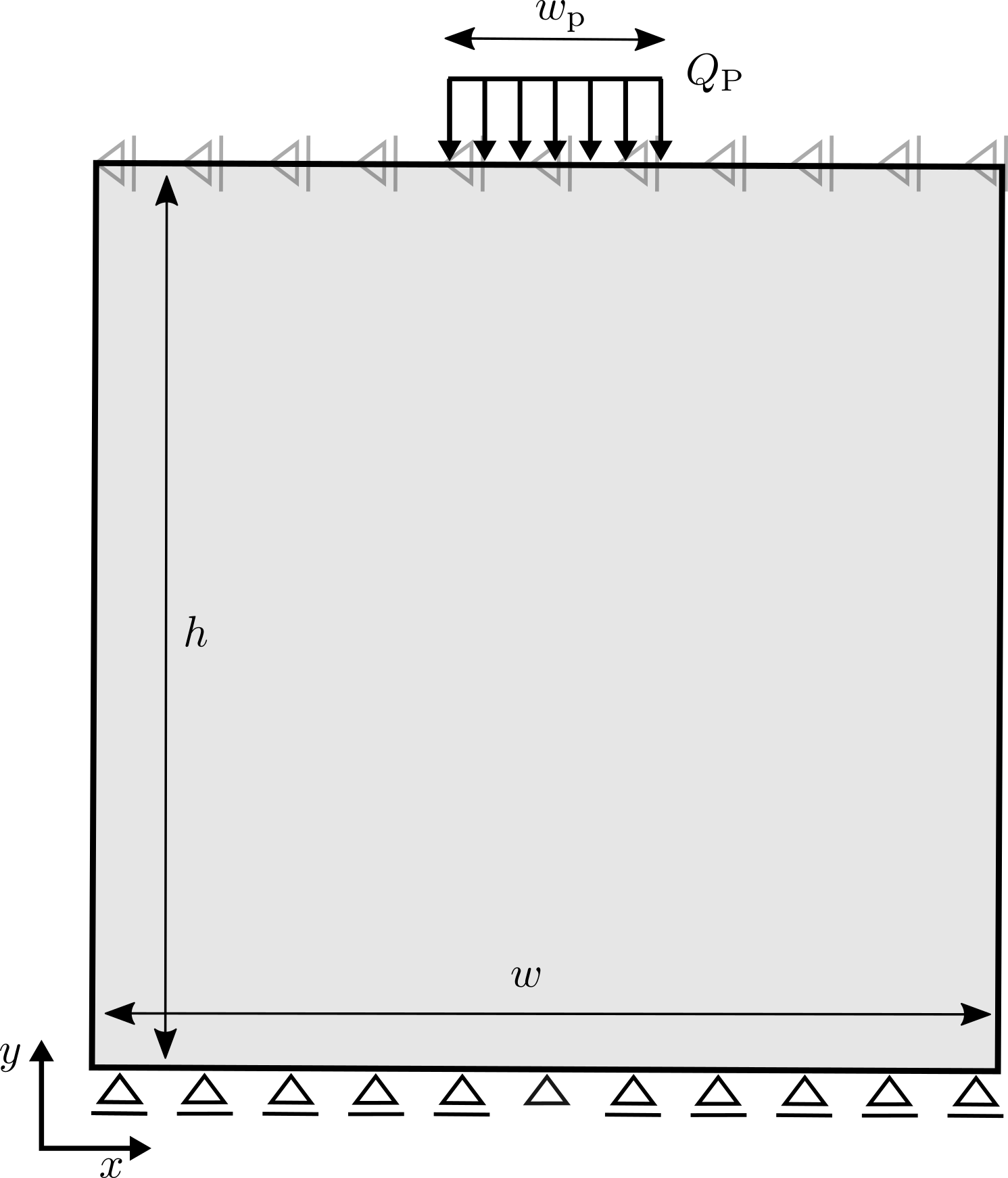}
			\caption{Problem geometry}
		\end{subfigure}%
		\begin{subfigure}[t]{0.55\textwidth}
			\centering
			\includegraphics[width=0.95\textwidth]{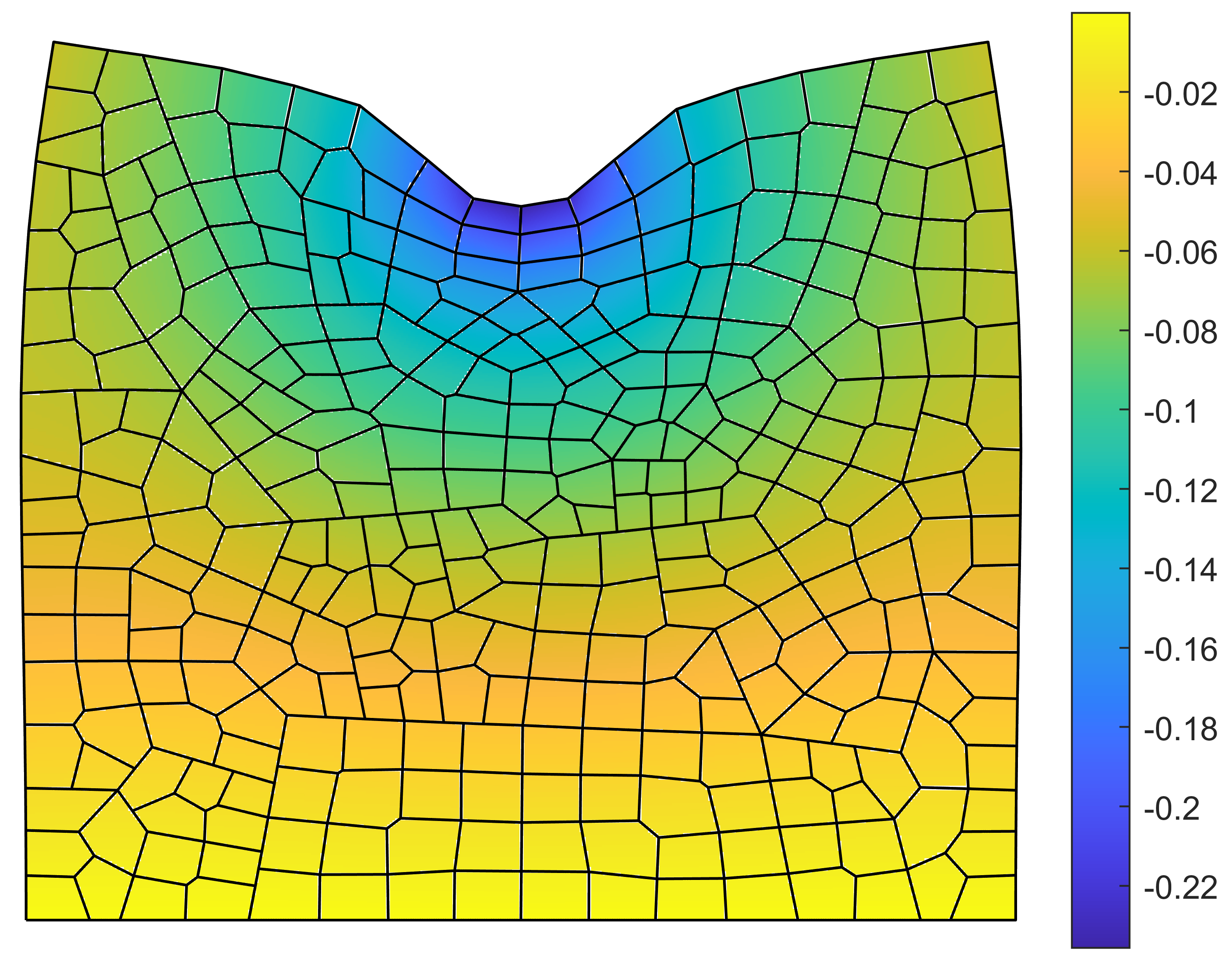}
			\caption{Deformed configuration}
		\end{subfigure}
		\caption{Punch problem (a) geometry, and (b) sample deformed configuration of a Voronoi mesh. 
			\label{fig:PunchGeometry}}
	\end{figure}

	The punch problem comprises a domain of width ${w=1~\rm{m}}$ and height ${h=1~\rm{m}}$ into which a punch of width ${w_{\rm p}=0.2~w}$ is driven into the middle of the top edge.
	The bottom edge of the domain is constrained vertically and the midpoint of the bottom edge is fully constrained.
	The top edge is constrained horizontally and the punch is modelled as a uniformly distributed load with a magnitude of ${Q_{\rm{P}}=0.675~\frac{\rm N}{\rm m}}$ (see Figure~\ref{fig:PunchGeometry}(a)). A sample deformed configuration of the body with a Voronoi mesh is depicted in Figure~\ref{fig:PunchGeometry}(b) with the vertical displacement ${u_{y}}$ plotted on the colour axis. The punch problem has a simple domain geometry which does not introduce any challenges in modelling this problem. While the punch introduces localized deformation only in its vicinity, the rest of the body experiences very little deformation (see Figure~\ref{fig:PunchGeometry}(b)). Thus, the punch problem is used to provide insight into the efficacy of the proposed coarsening procedures in cases of `less challenging´ problems.
	
	\FloatBarrier
	\FloatBarrier
	
\begin{figure}[ht!]
		\centering
		\begin{subfigure}[t]{0.33\textwidth}
			\centering
			\includegraphics[width=0.95\textwidth]{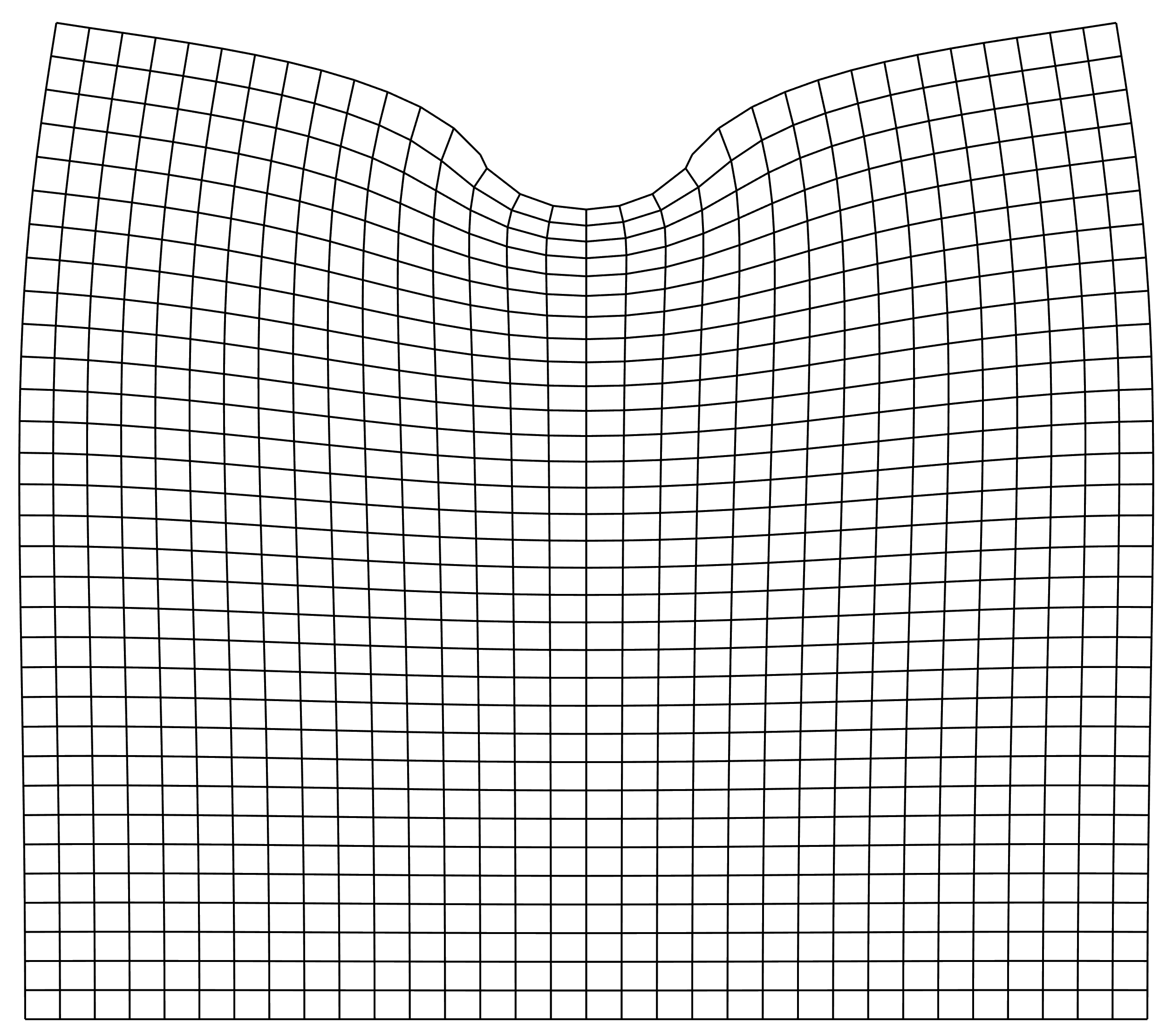}
			\caption{Structured - Step 1}
		\end{subfigure}%
		\begin{subfigure}[t]{0.33\textwidth}
			\centering
			\includegraphics[width=0.95\textwidth]{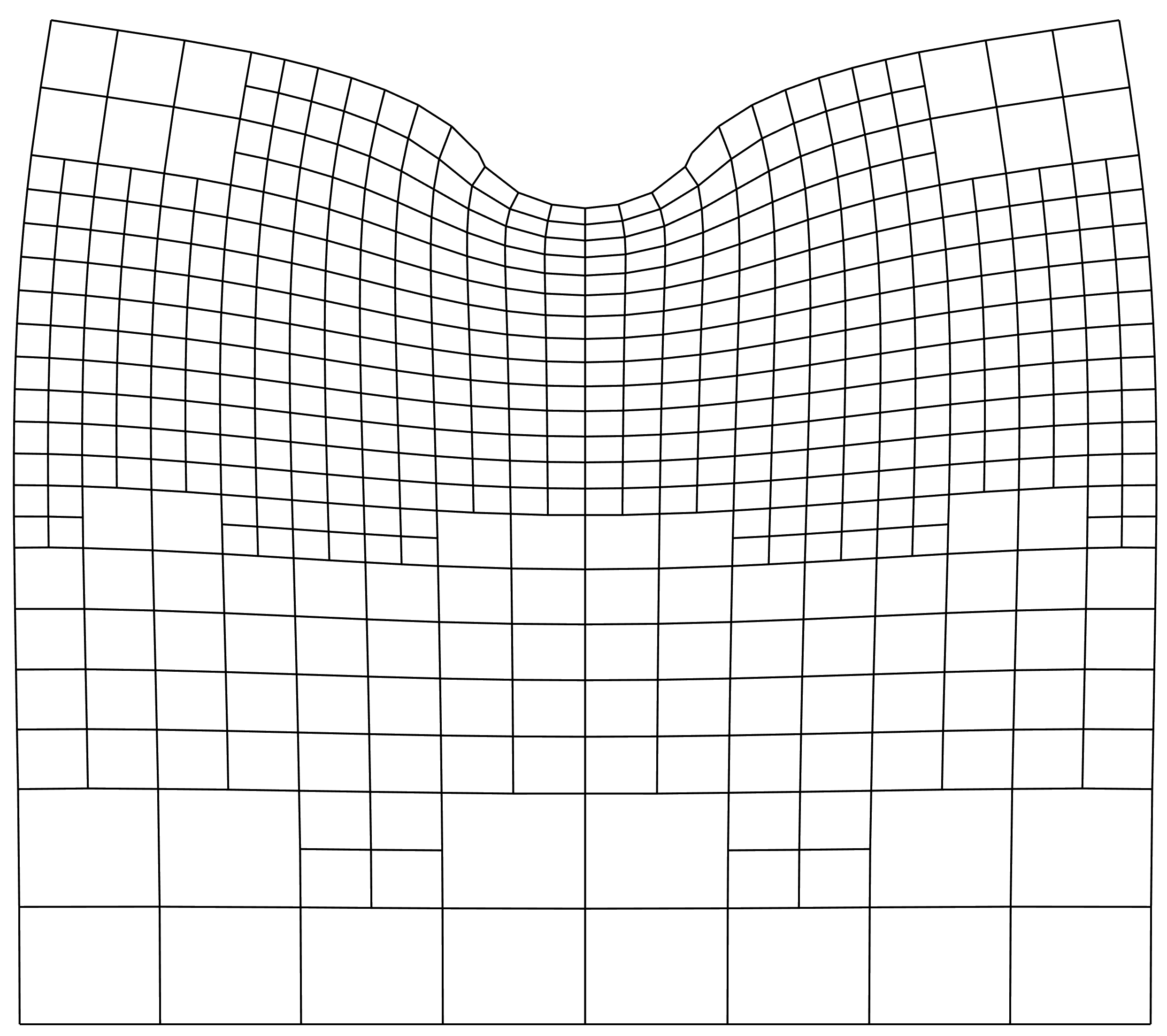}
			\caption{Structured - Step 5}
		\end{subfigure}%
		\begin{subfigure}[t]{0.33\textwidth}
			\centering
			\includegraphics[width=0.95\textwidth]{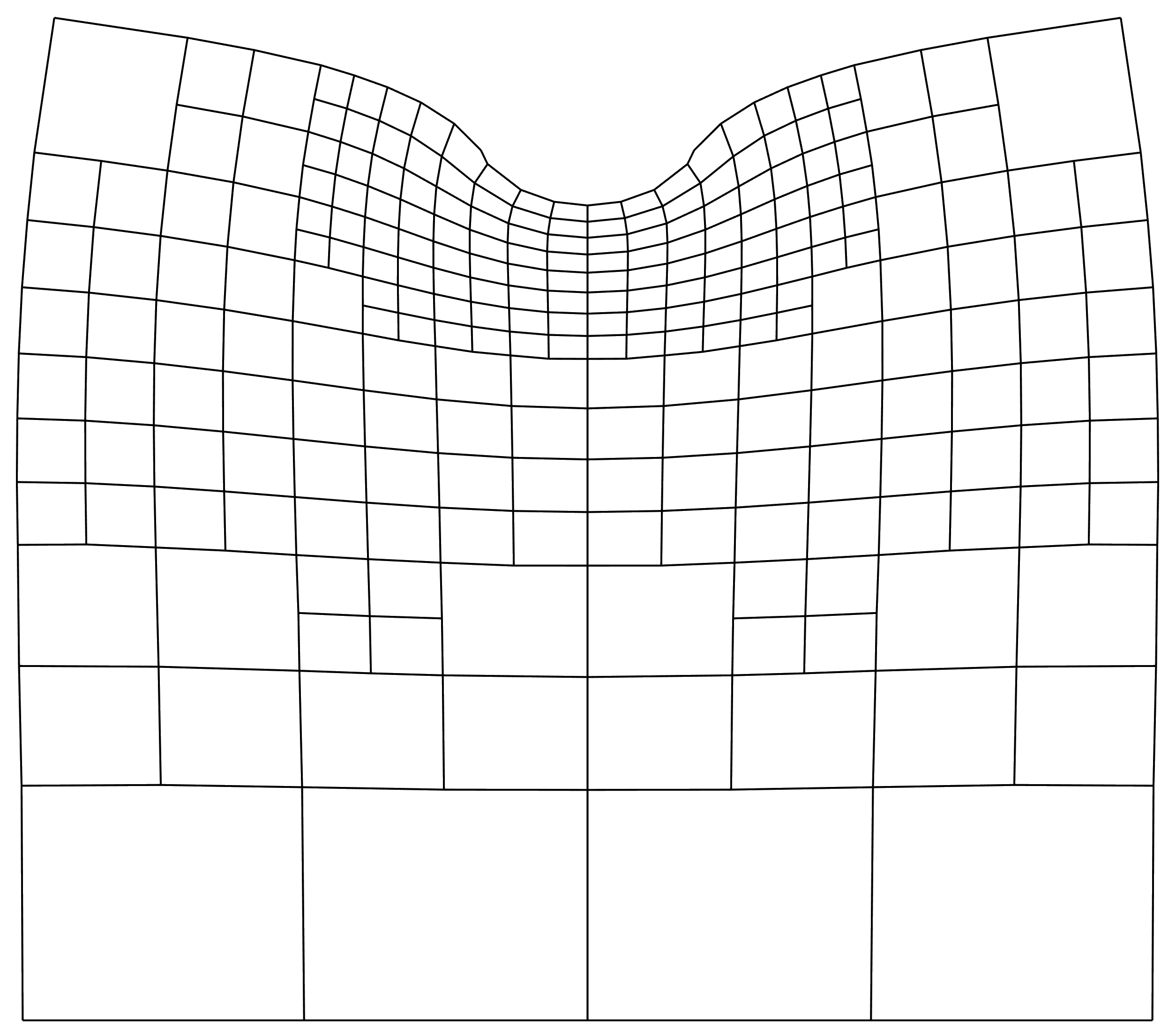}
			\caption{Structured - Step 10}
		\end{subfigure}
		\vskip \baselineskip 
		\begin{subfigure}[t]{0.33\textwidth}
			\centering
			\includegraphics[width=0.95\textwidth]{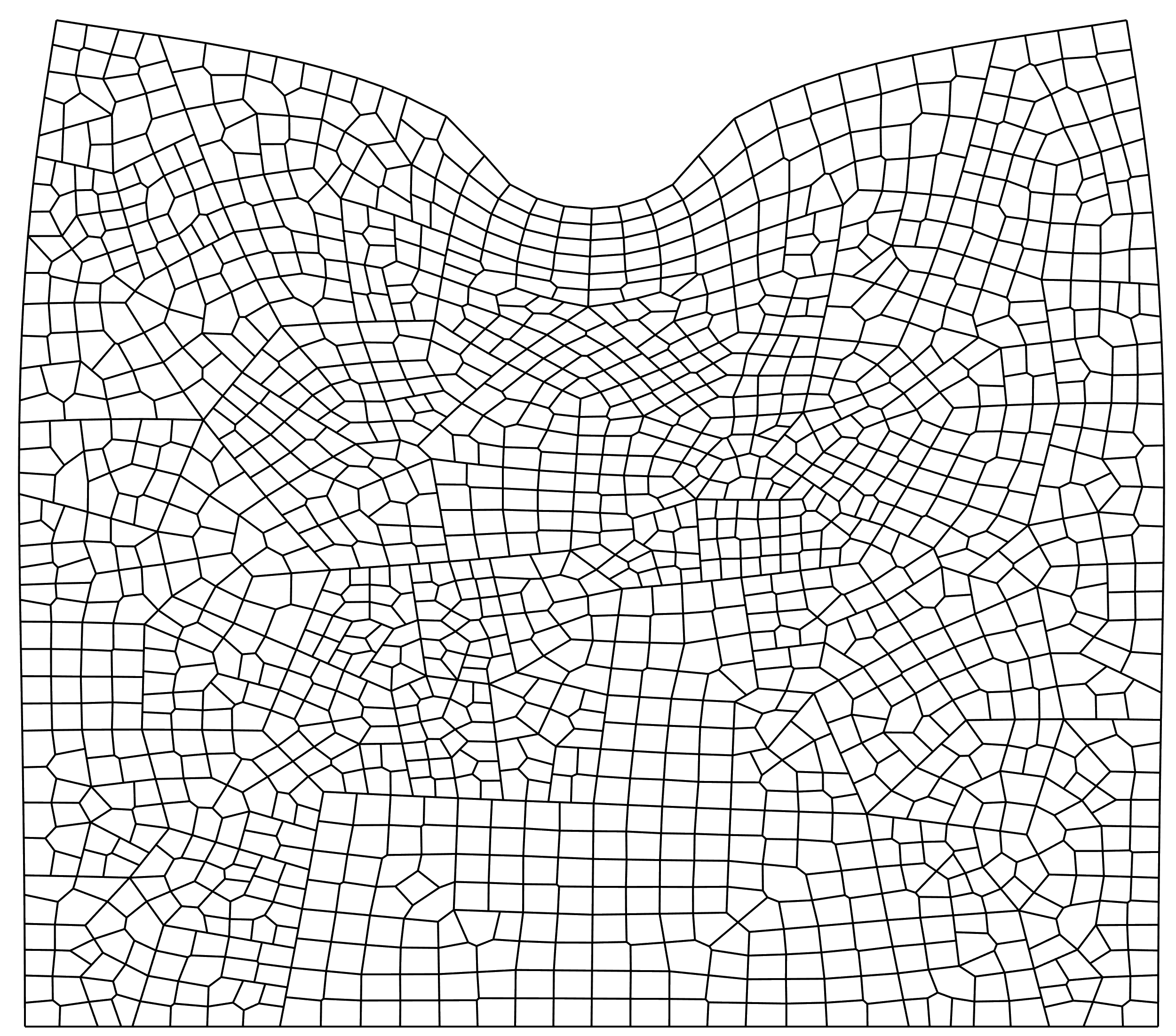}
			\caption{Voronoi - Step 1}
		\end{subfigure}%
		\begin{subfigure}[t]{0.33\textwidth}
			\centering
			\includegraphics[width=0.95\textwidth]{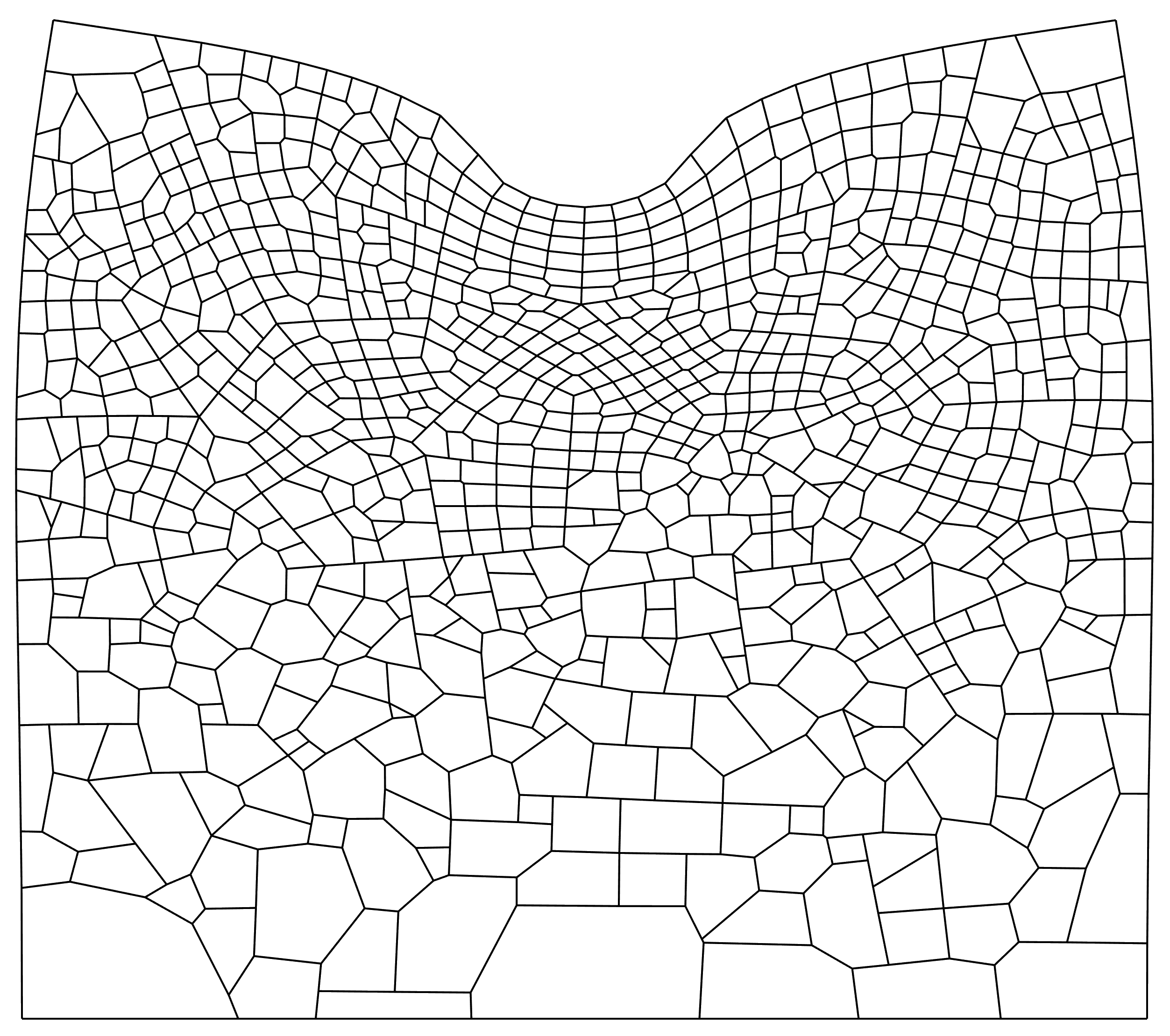}
			\caption{Voronoi - Step 6}
		\end{subfigure}%
		\begin{subfigure}[t]{0.33\textwidth}
			\centering
			\includegraphics[width=0.95\textwidth]{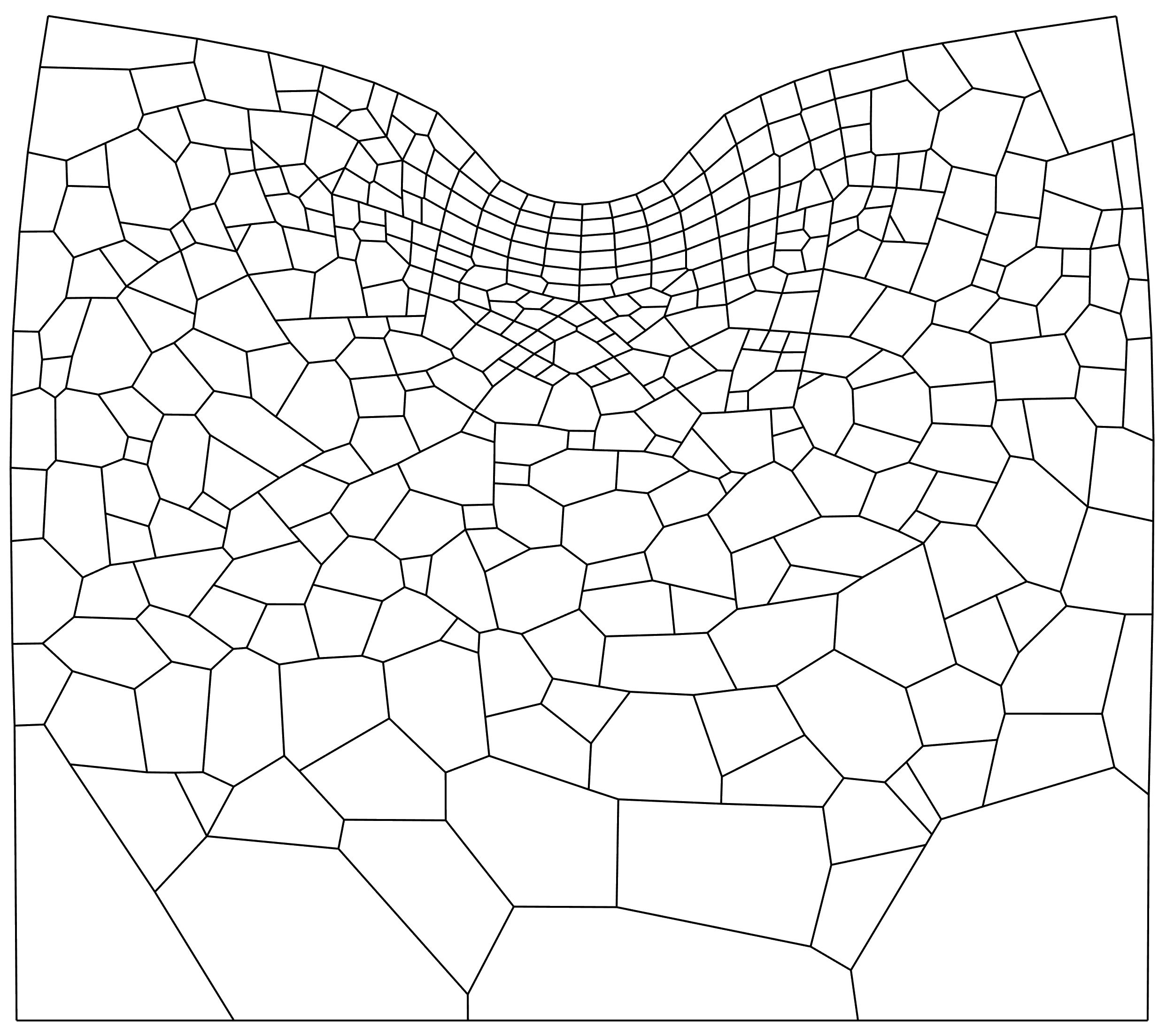}
			\caption{Voronoi - Step 12}
		\end{subfigure}
		\caption{Mesh coarsening process for the punch problem on structured and Voronoi meshes using the displacement-based coarsening procedure with ${T=20\%}$.
			\label{fig:PunchMeshes}}
	\end{figure}

	The mesh evolution during the coarsening process for the punch problem is depicted in Figure~\ref{fig:PunchMeshes} for the case of the displacement-based coarsening procedure with ${T=20\%}$ on structured and unstructured/Voronoi meshes.
	Meshes are shown at various coarsening steps with step~1 corresponding to the initial mesh. 
	Similar mesh evolution is exhibited on both structured and Voronoi meshes. The mesh becomes increasingly coarse at the bottom of the domain while remaining fine around the region of the punch. Furthermore, the mesh density exhibits a smooth and graded transition from the finest region around the punch to the coarser regions further from the punch.
	In this example problem most of the deformation occurs around the location of the punch. The rest of the body experiences comparatively little deformation with the magnitude of the deformation decreasing with increasing distance from the punch (see Figure~\ref{fig:PunchGeometry}(b)). 
	This behaviour is reflected in the mesh coarsening process, thus, the mesh evolution is sensible for this problem.
	
	\FloatBarrier
	\FloatBarrier
	
\begin{figure}[ht!]
		\centering
		\begin{subfigure}[t]{0.33\textwidth}
			\centering
			\includegraphics[width=0.95\textwidth]{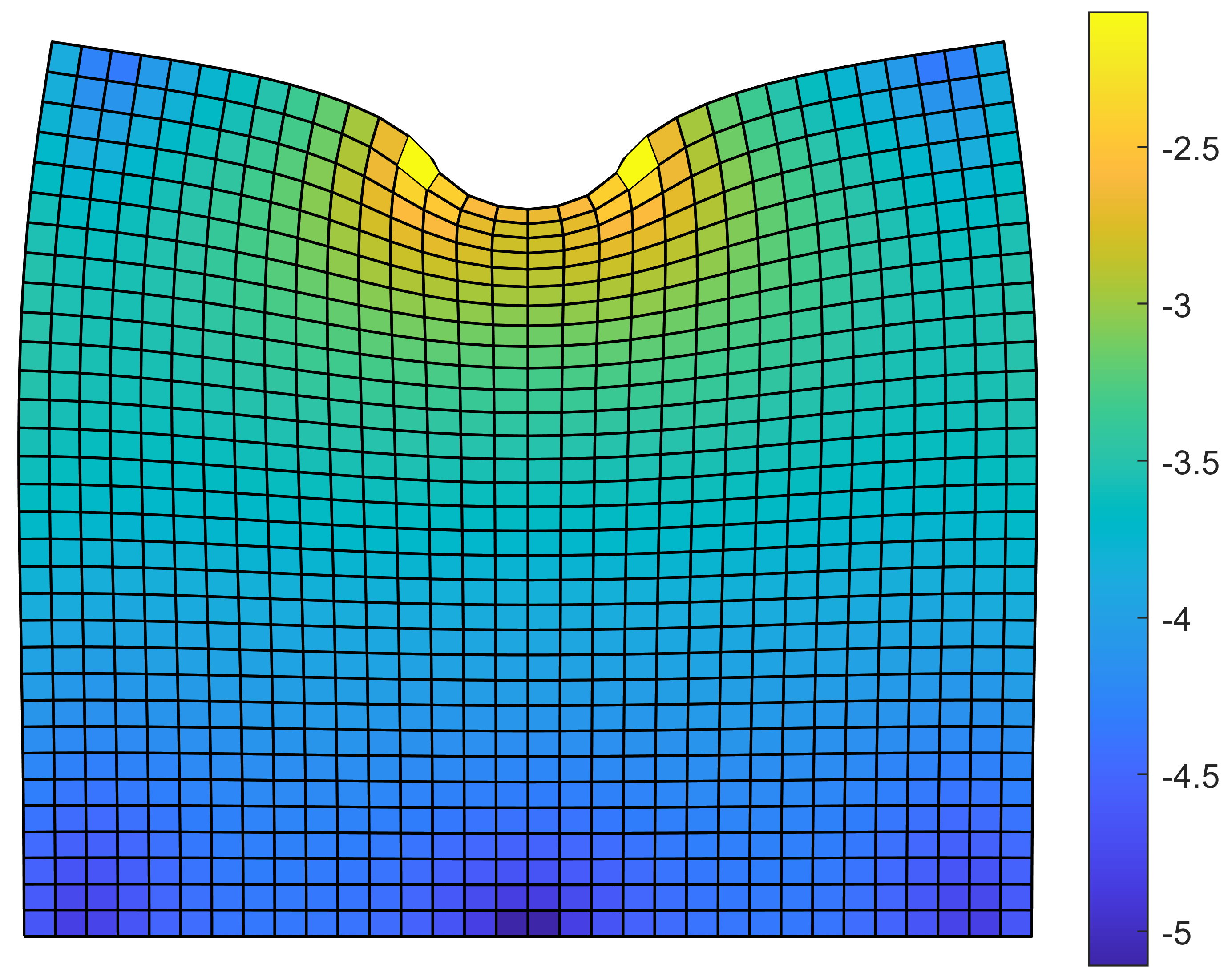}
			\caption{Structured - Step 1}
		\end{subfigure}%
		\begin{subfigure}[t]{0.33\textwidth}
			\centering
			\includegraphics[width=0.95\textwidth]{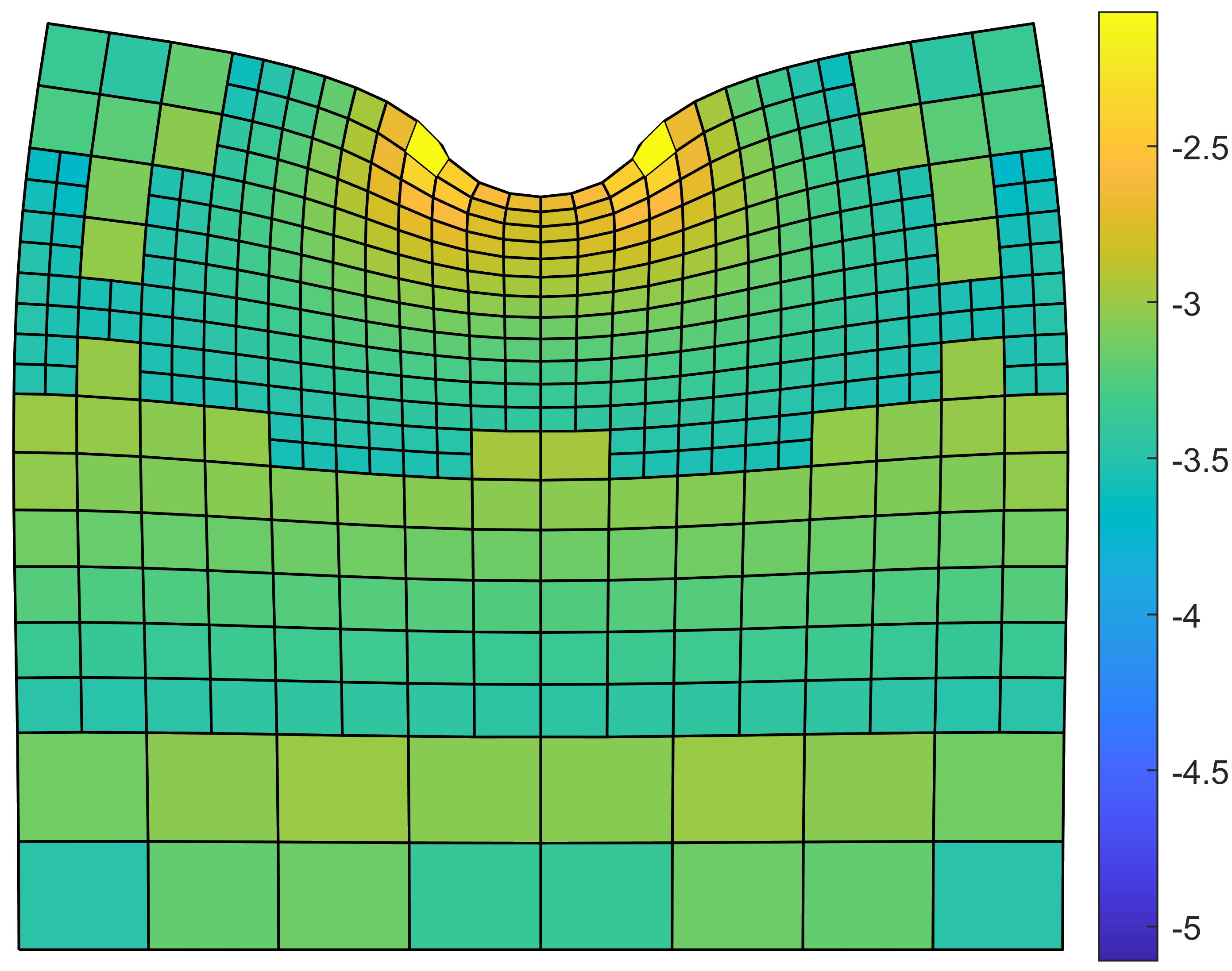}
			\caption{Structured - Step 6}
		\end{subfigure}%
		\begin{subfigure}[t]{0.33\textwidth}
			\centering
			\includegraphics[width=0.95\textwidth]{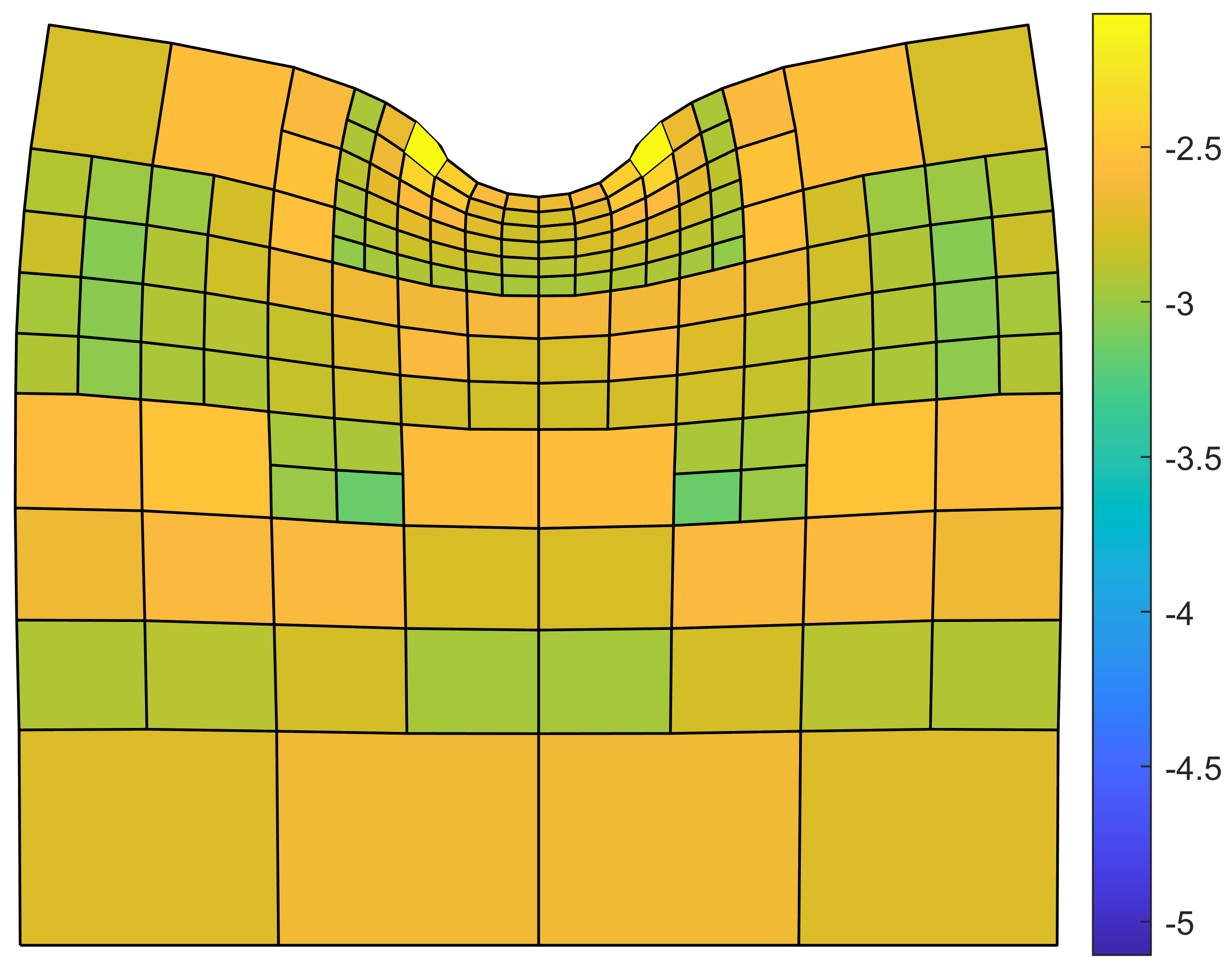}
			\caption{Structured - Step 12}
		\end{subfigure}
		\vskip \baselineskip 
		\begin{subfigure}[t]{0.33\textwidth}
			\centering
			\includegraphics[width=0.95\textwidth]{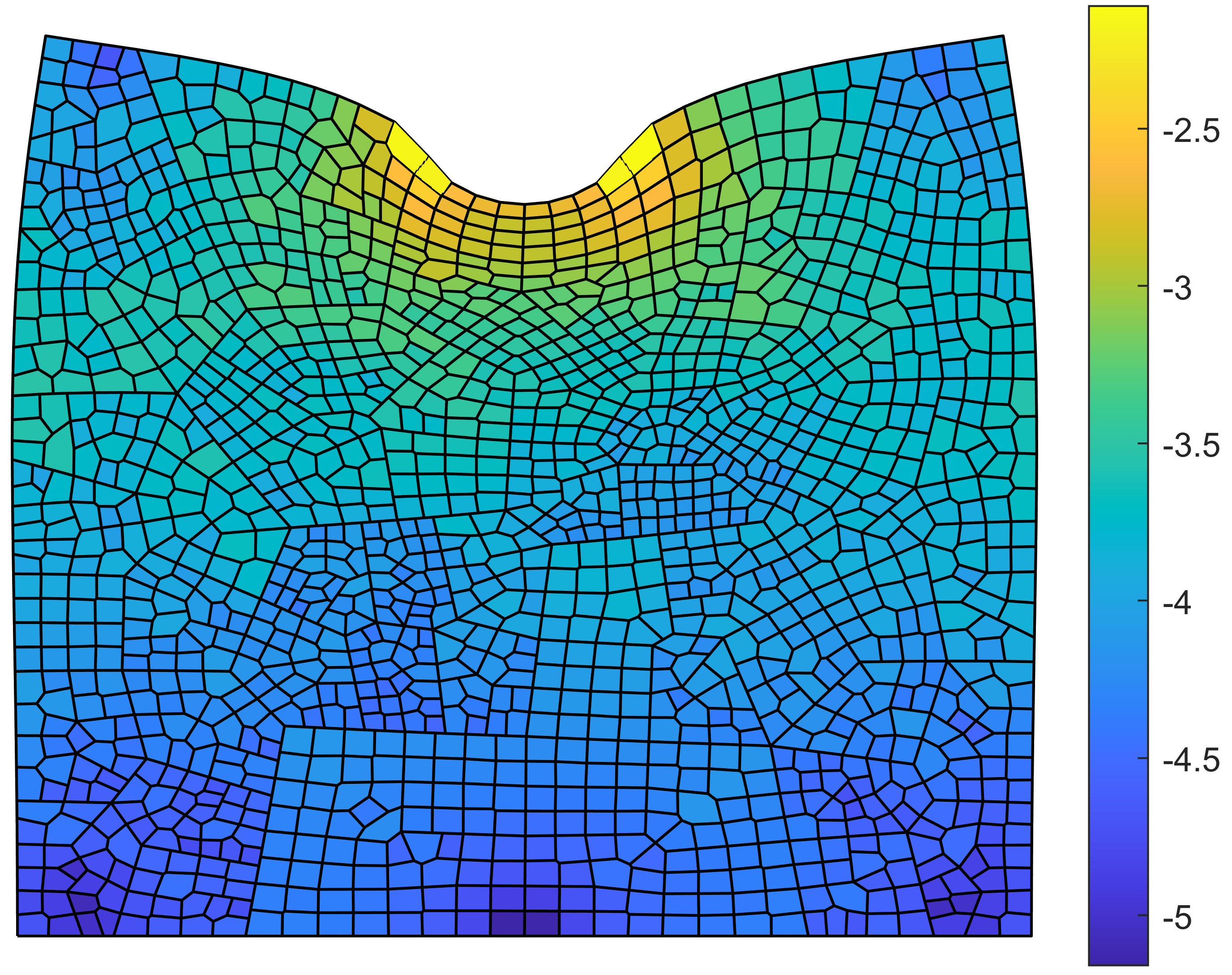}
			\caption{Voronoi - Step 1}
		\end{subfigure}%
		\begin{subfigure}[t]{0.33\textwidth}
			\centering
			\includegraphics[width=0.95\textwidth]{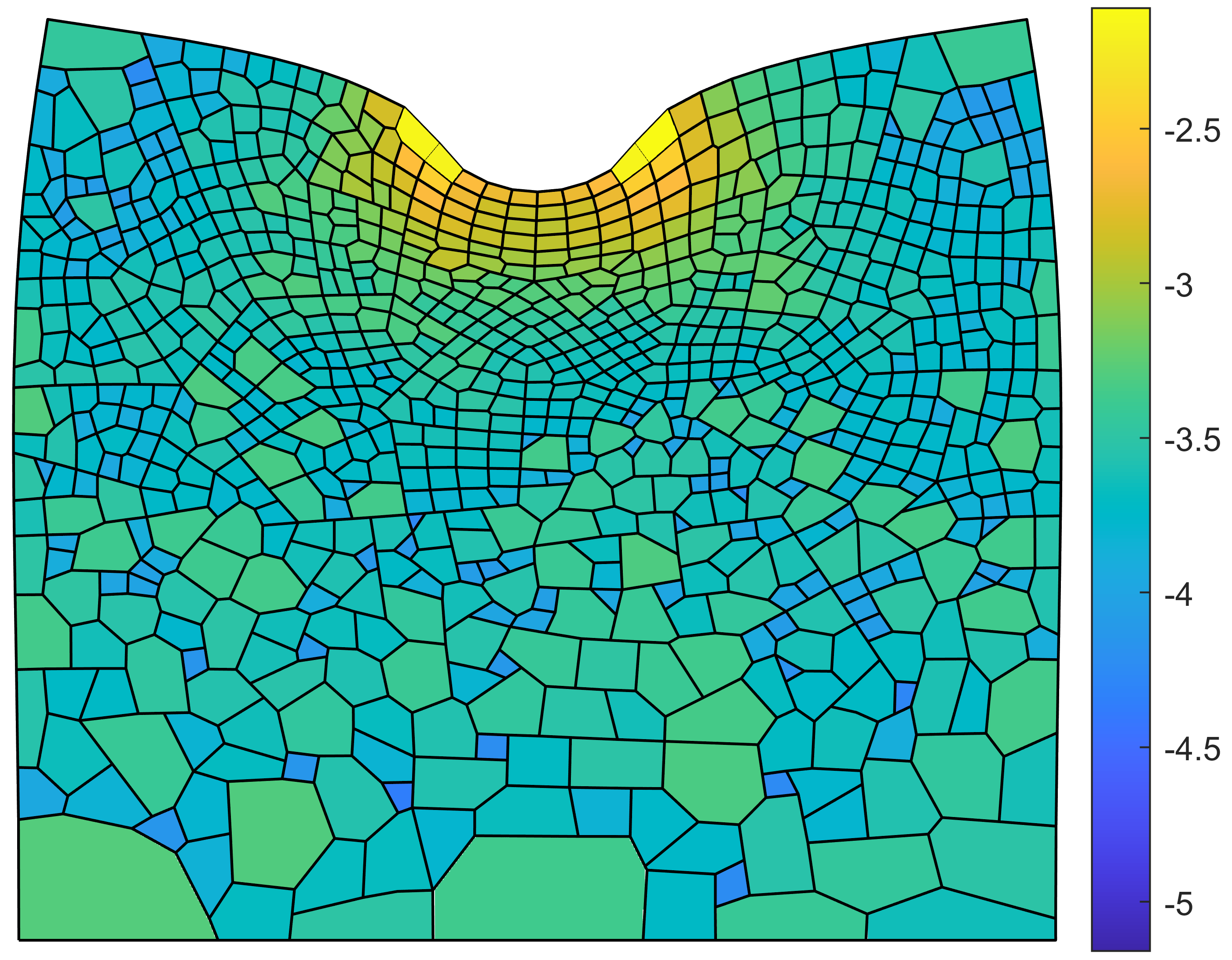}
			\caption{Voronoi - Step 6}
		\end{subfigure}%
		\begin{subfigure}[t]{0.33\textwidth}
			\centering
			\includegraphics[width=0.95\textwidth]{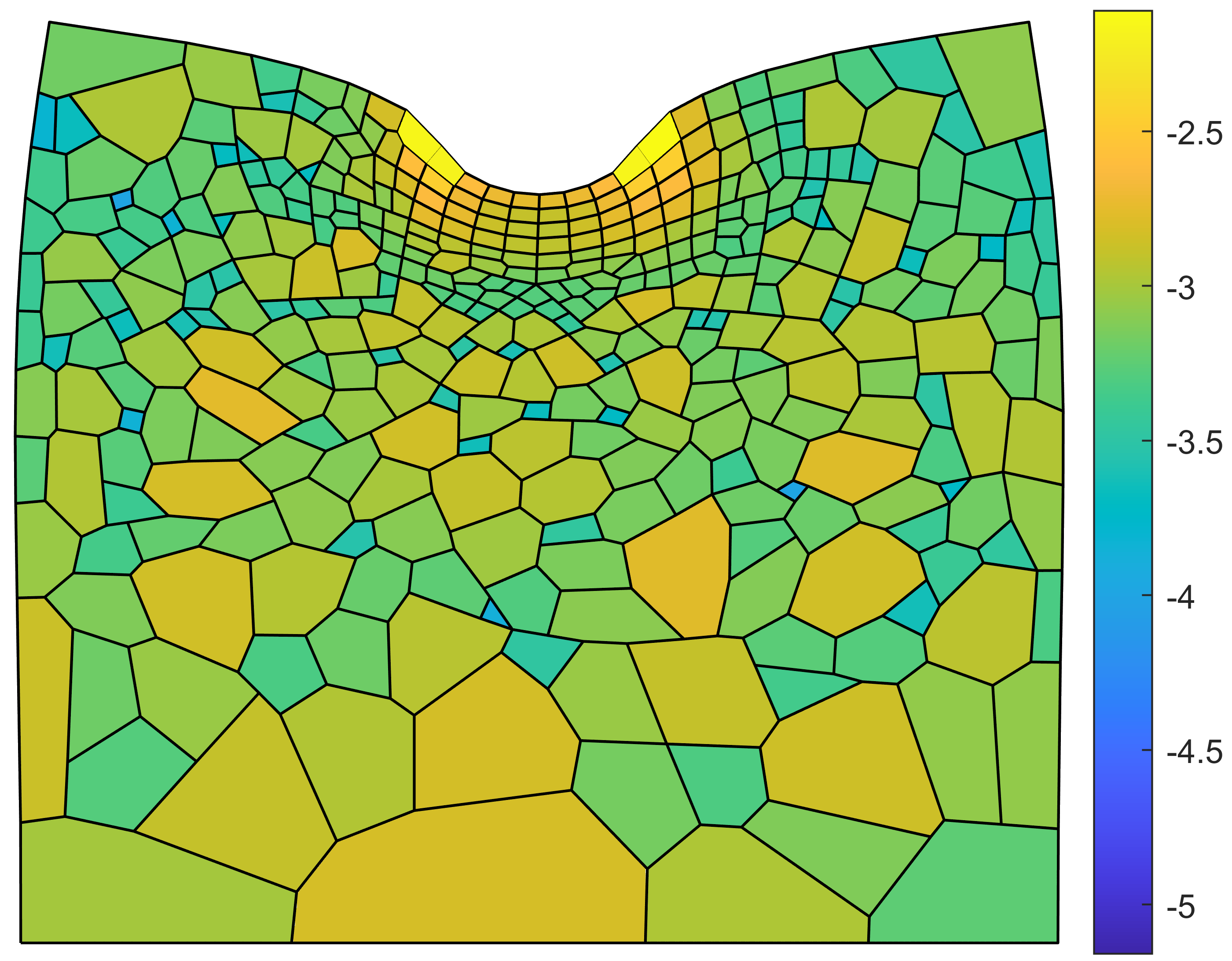}
			\caption{Voronoi - Step 12}
		\end{subfigure}
		\caption{${\mathcal{H}^{1}}$ error distribution during the coarsening process for the punch problem on structured and Voronoi meshes using the displacement-based coarsening procedure with ${T=20\%}$.
			\label{fig:PunchErrorMeshes}}
	\end{figure}

	The distribution of the ${\mathcal{H}^{1}}$ error over the domain during the mesh coarsening process for the punch problem is depicted in Figure~\ref{fig:PunchErrorMeshes}. The ${\mathcal{H}^{1}}$ error is depicted in a logarithmic scale on structured and Voronoi meshes for the case of the displacement-based coarsening procedure with ${T=20\%}$. The error distribution exhibited in step~1, i.e. Figure~\ref{fig:PunchErrorMeshes}(a), demonstrates that the mesh evolution illustrated in Figure~\ref{fig:PunchMeshes} is indeed sensible, as discussed. The mesh evolution closely reflects the error distribution and the regions with the lowest errors are coarsened most aggressively.
	Furthermore, the error distribution over the domain becomes increasingly even as the coarsening procedure progresses. Since an optimal mesh would have a perfectly even error distribution, this behaviour demonstrates that the mesh becomes closer to optimal during the coarsening procedure. It is noted that due to the discrete increases in element size during coarsening, a perfectly even error distribution is not possible as it would require precise resizing of individual elements. Nevertheless, the improved error distribution during coarsening demonstrates the high degree of efficacy of the proposed procedure.
	
	\FloatBarrier
	\FloatBarrier
	
\begin{figure}[ht!]
		\centering
		\begin{subfigure}[t]{0.5\textwidth}
			\centering
			\includegraphics[width=0.95\textwidth]{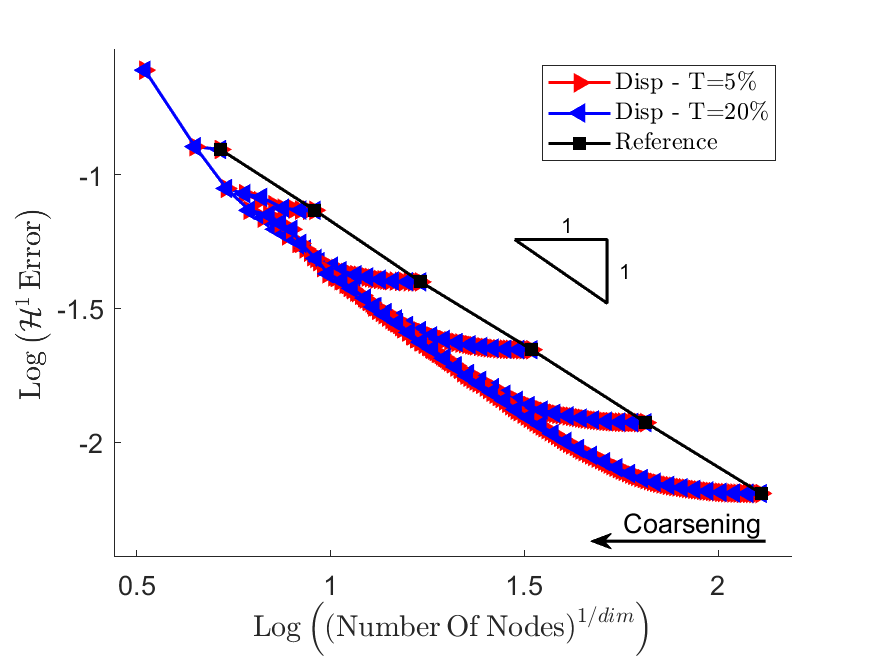}
			\caption{Displacement-based indicator - Structured meshes}
		\end{subfigure}%
		\begin{subfigure}[t]{0.5\textwidth}
			\centering
			\includegraphics[width=0.95\textwidth]{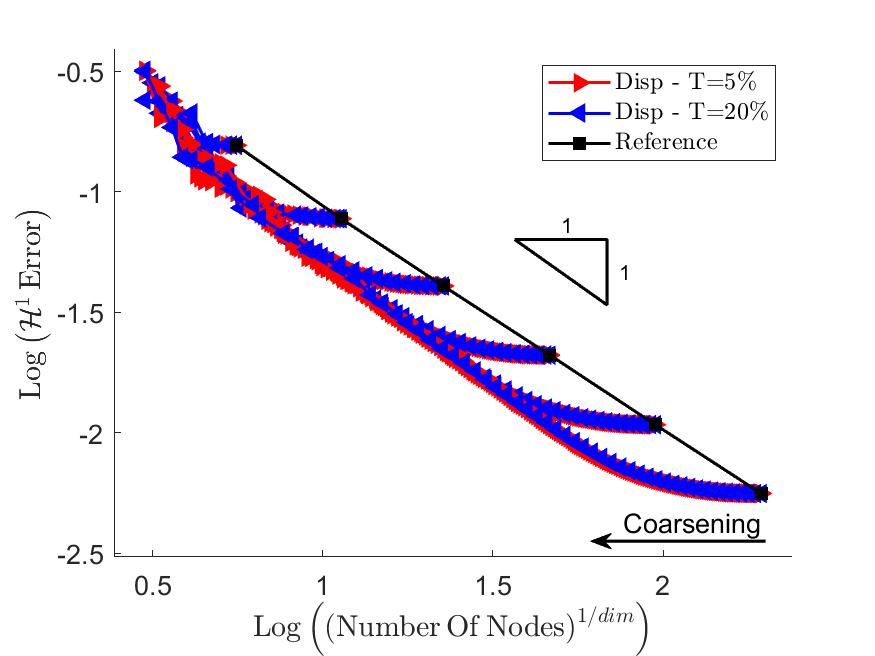}
			\caption{Displacement-based indicator - Voronoi meshes}
		\end{subfigure}
		\vskip \baselineskip 
		\begin{subfigure}[t]{0.5\textwidth}
			\centering
			\includegraphics[width=0.95\textwidth]{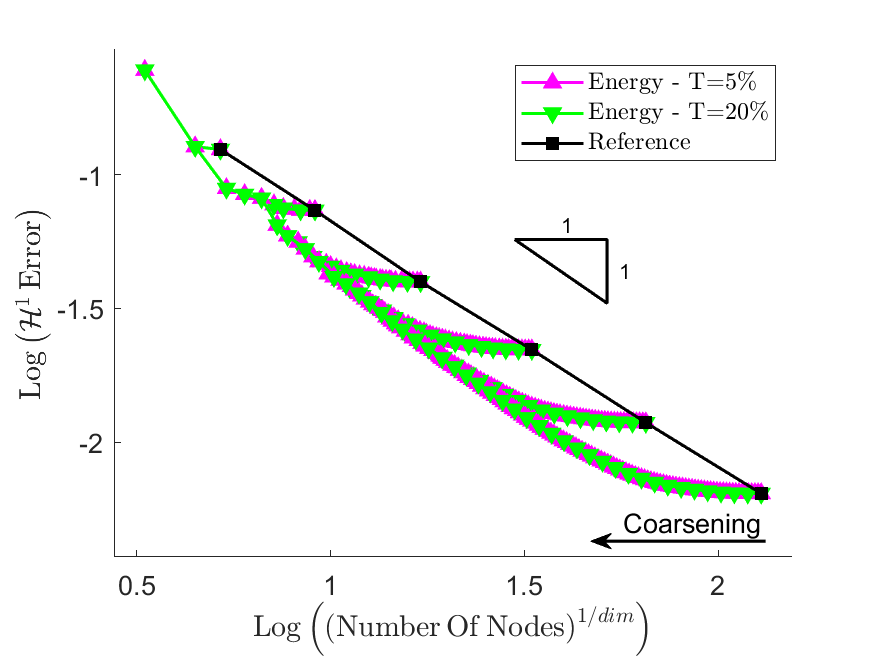}
			\caption{Energy error-based indicator - Structured meshes}
		\end{subfigure}%
		\begin{subfigure}[t]{0.5\textwidth}
			\centering
			\includegraphics[width=0.95\textwidth]{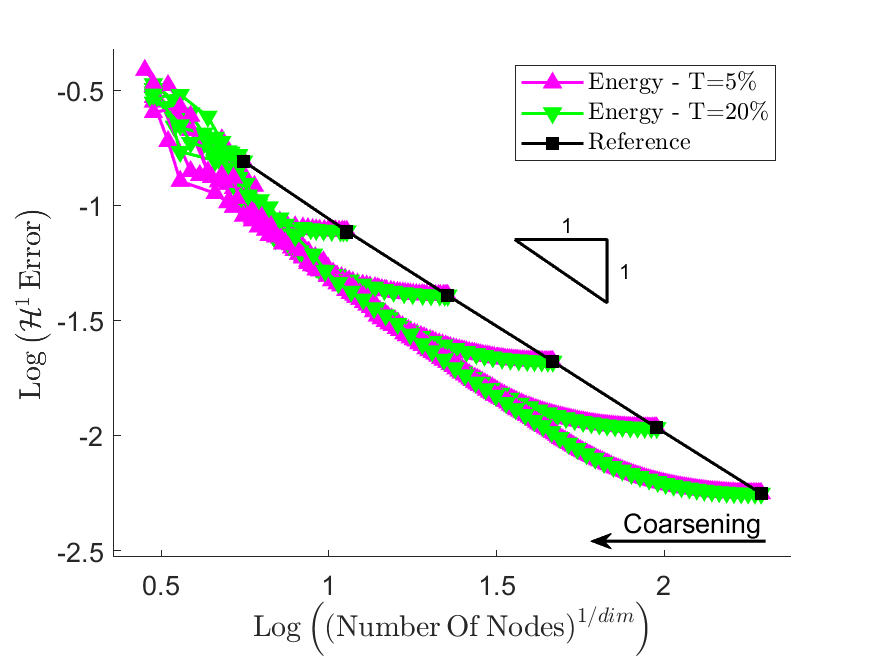}
			\caption{Energy error-based indicator - Voronoi meshes}
		\end{subfigure}
		\caption{$\mathcal{H}^{1}$ error vs $n_{\rm v}$ for the punch problem.
			\label{fig:PunchConvergenceNumberOfNodes}}
	\end{figure}

	The convergence behaviour in the ${\mathcal{H}^{1}}$ error norm of the VEM for the punch problem using the displacement-based and energy error-based coarsening procedures is depicted in Figure~\ref{fig:PunchConvergenceNumberOfNodes} on a logarithmic scale. Here, the convergence behaviour of the displacement-based and energy error-based procedures are presented on the top and bottom rows of figures respectively, with the results generated on structured and Voronoi meshes presented in the left and right columns of figures respectively. In each case results are presented for coarsening thresholds ${T=5\%}$ and ${T=20\%}$ to demonstrate the effect of the choice of $T$ on the convergence behaviour. Additionally, uniform initial meshes of various discretization densities are considered.
	Even though the punch problem is `less challenging´ and can be suitably analysed using uniform meshes the benefit of the coarsening procedures is clear. Both coarsening procedures successfully coarsen the least critical elements in the domain, which eliminates the least important degrees of freedom. 
	The benefit of the coarsening procedures is clear from the path of the ${\mathcal{H}^{1}}$ error curves. The curve is initially almost horizontal as the least significant nodes and elements are coarsened and almost no error is introduced. As the procedure continues the coarsening spreads over the domain which inevitably does begin to introduce some error and the gradients of the ${\mathcal{H}^{1}}$ error curves increase. However, the efficiency of the coarsened solutions is superior to that of uniform meshes (indicated by the black curve denoting the reference uniform mesh approach).
	That is, for a given number of degrees of freedom, the coarsened meshes yield a lower ${\mathcal{H}^{1}}$ error than the uniform meshes. Additionally, similar behaviour is exhibited in the cases of both structured and Voronoi meshes. 
	Interestingly, the behaviour does not appear to be significantly influenced by the choice of $T$, with nearly identical performance exhibited in the cases of ${T=5\%}$ and ${T=20\%}$ for both coarsening procedures and on both mesh types. Finally, the efficacy of the coarsening procedures increases with increasing density of the initial uniform mesh. This behaviour is sensible because the more dense the initial mesh is, the more elements and the more degrees of freedom there are in the least critical portions of the domain. These regions can then be efficiently coarsened while introducing the least amount of ${\mathcal{H}^{1}}$ error possible.
	
	\FloatBarrier
	\FloatBarrier
	
	\subsection{Plate with hole}
	\label{subsec:PlateWithHoleTraction}
	
\begin{figure}[ht!]
		\centering
		\begin{subfigure}[t]{0.45\textwidth}
			\centering
			\includegraphics[width=0.95\textwidth]{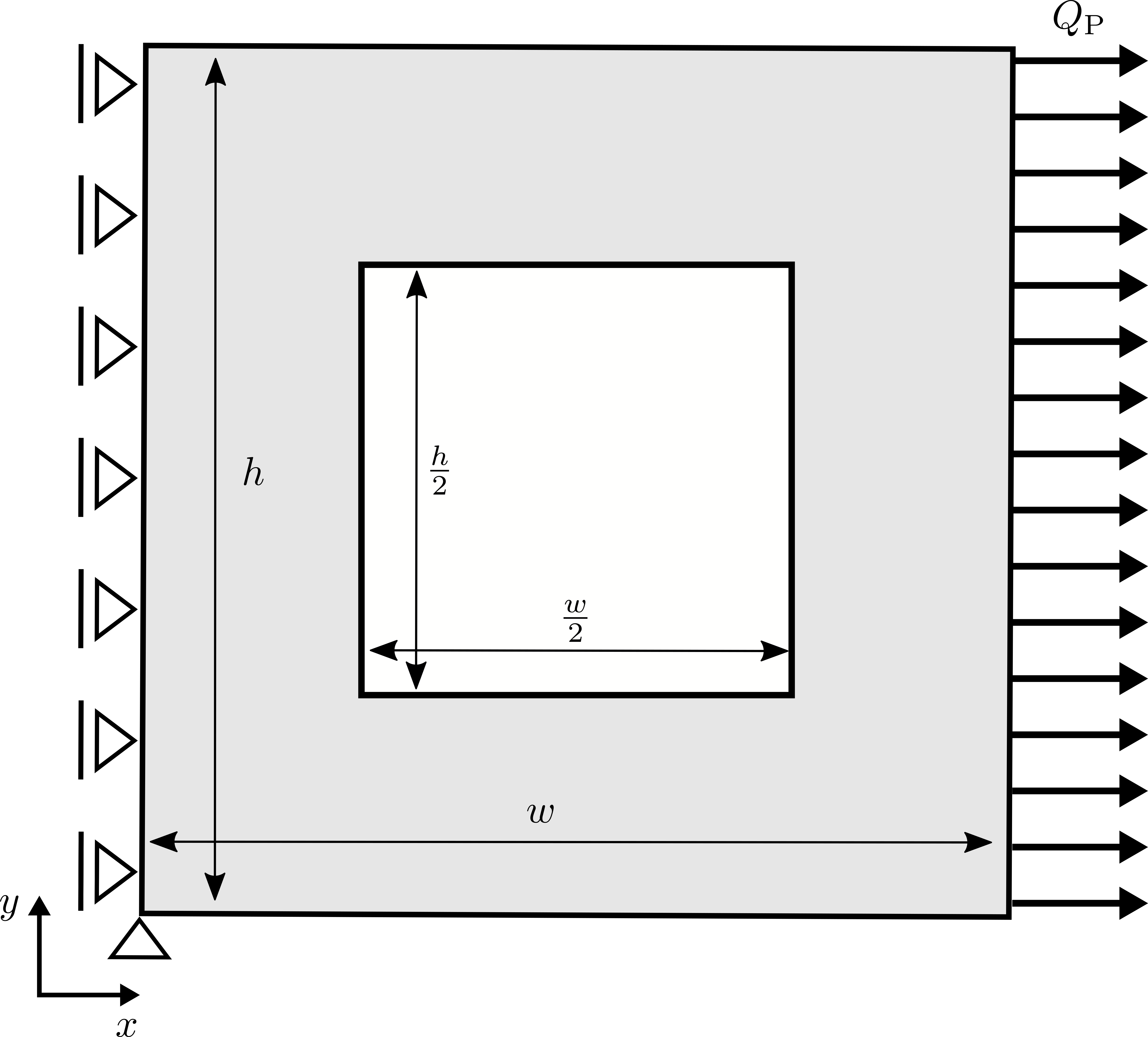}
			\caption{Problem geometry}
		\end{subfigure}%
		\begin{subfigure}[t]{0.55\textwidth}
			\centering
			\includegraphics[width=0.95\textwidth]{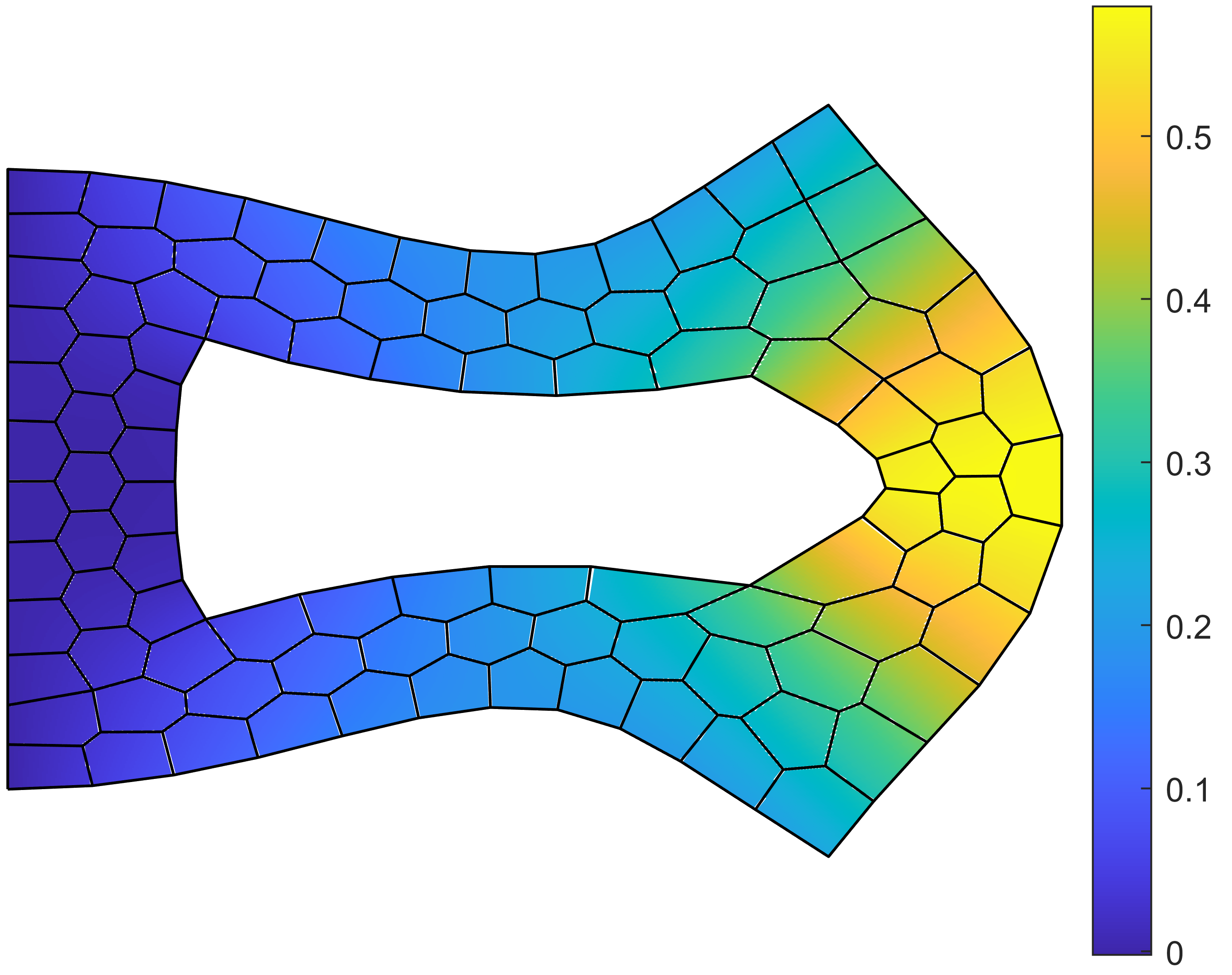}
			\caption{Deformed configuration}
		\end{subfigure}
		\caption{Plate with hole (a) geometry, and (b) deformed configuration of a Voronoi mesh. 
			\label{fig:PlateWithHoleTractionGeometry}}
	\end{figure}

	The plate with hole problem comprises a domain of width ${w=1~\rm{m}}$ and height ${h=1~\rm{m}}$ with a centrally located square hole (see Figure~\ref{fig:PlateWithHoleTractionGeometry}(a)). The left-hand edge of the plate is constrained horizontally and the bottom left-hand corner is fully constrained. The right-hand edge is subject to a prescribed traction of ${Q_{\rm P}=0.2~\frac{\rm N}{\rm m}}$. A sample deformed configuration of the plate with a Voronoi mesh is depicted in Figure~\ref{fig:PlateWithHoleTractionGeometry}(b) with the horizontal displacement ${u_{x}}$ plotted on the colour axis.
	The boundary conditions of the plate with hole problem are simple and themselves do not introduce any challenges in modelling this problem. However, the sharp corners of the hole act as stress concentrators which introduce high stresses and complex deformation. Thus, this problem is used to provide insight into the efficacy of the proposed coarsening procedures in cases of `more challenging´ problems.
	
	\FloatBarrier
	\FloatBarrier
	
\begin{figure}[ht!]
		\centering
		\begin{subfigure}[t]{0.33\textwidth}
			\centering
			\includegraphics[width=0.95\textwidth]{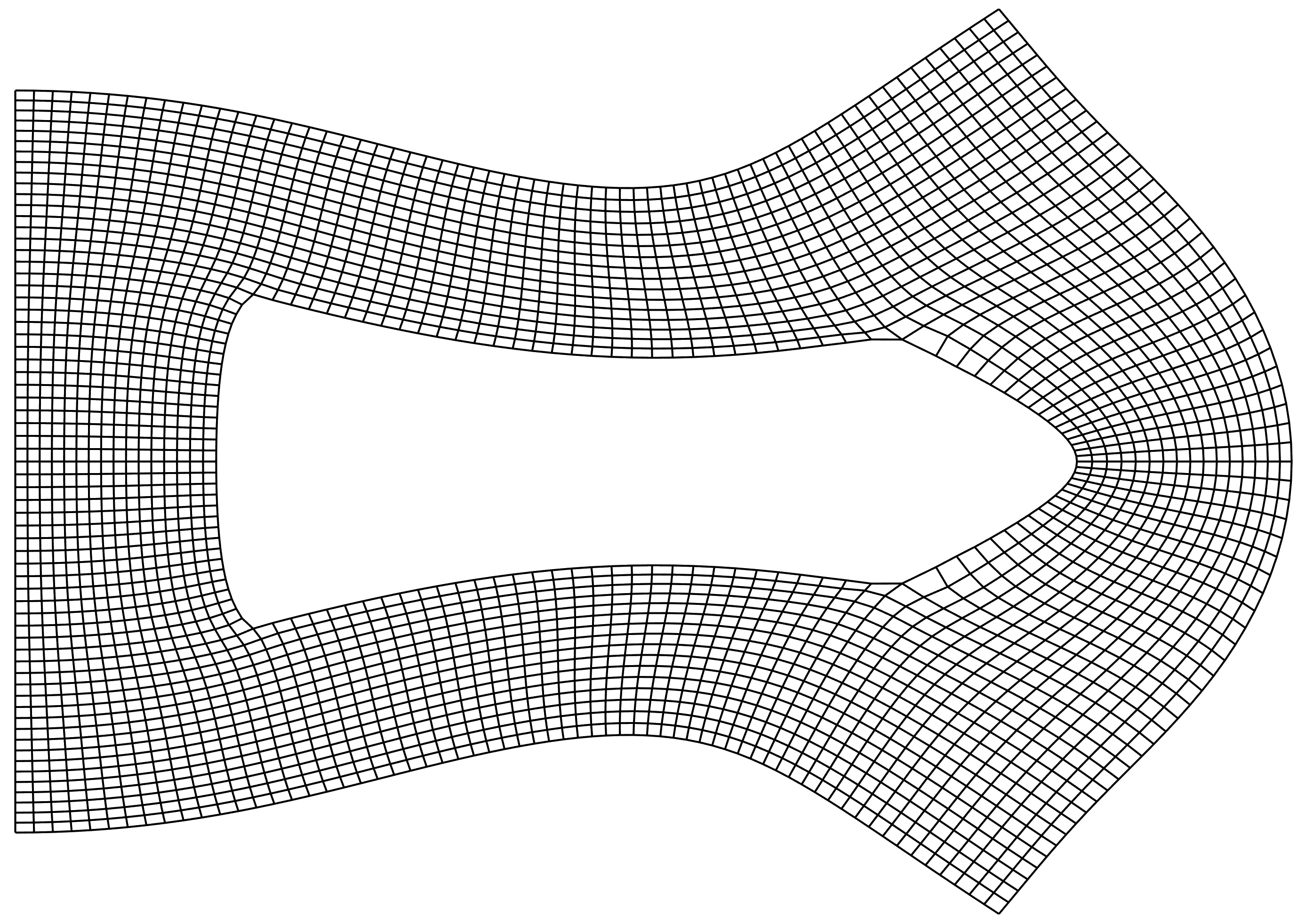}
			\caption{Structured - Step 1}
		\end{subfigure}%
		\begin{subfigure}[t]{0.33\textwidth}
			\centering
			\includegraphics[width=0.95\textwidth]{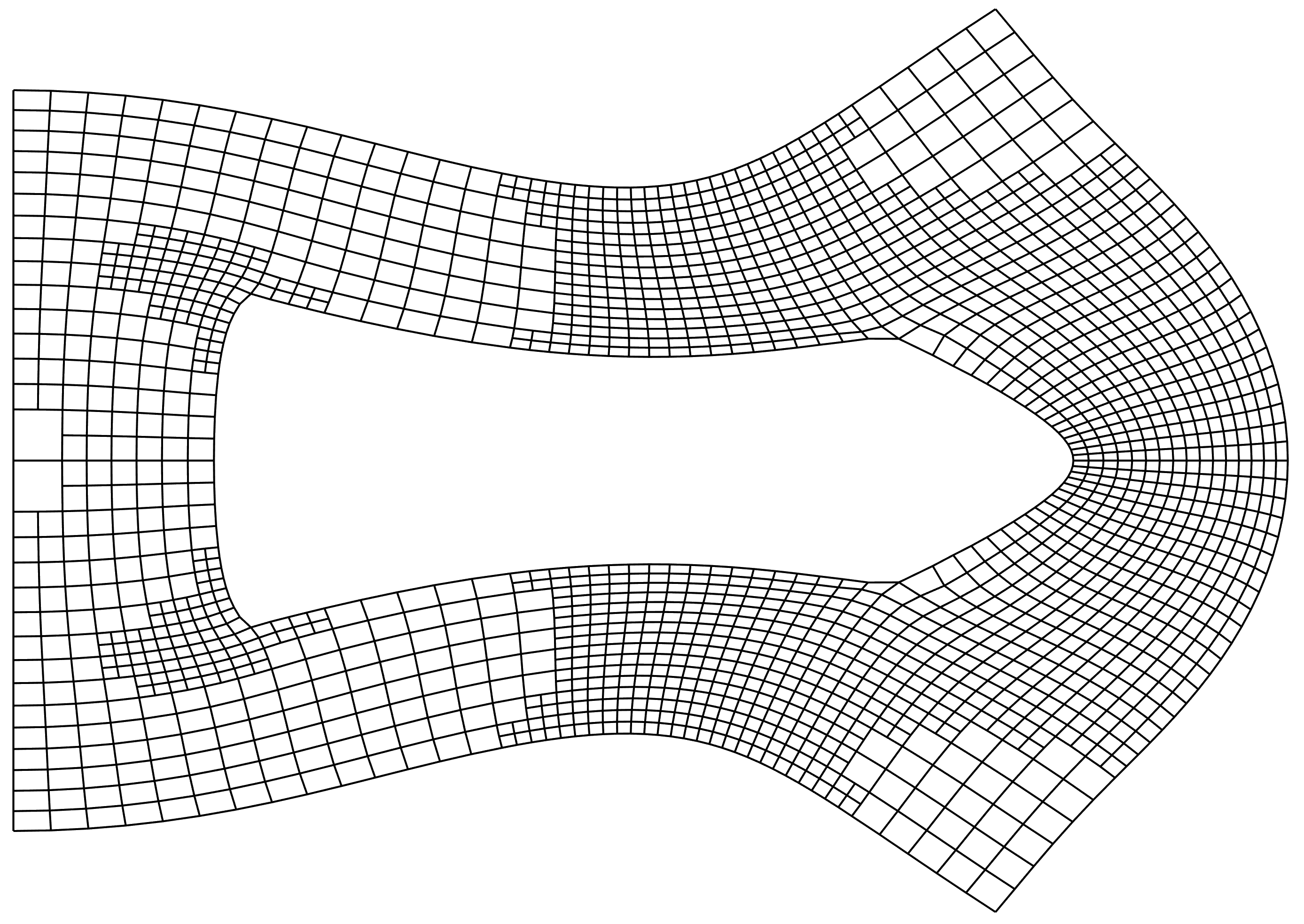}
			\caption{Structured - Step 4}
		\end{subfigure}%
		\begin{subfigure}[t]{0.33\textwidth}
			\centering
			\includegraphics[width=0.95\textwidth]{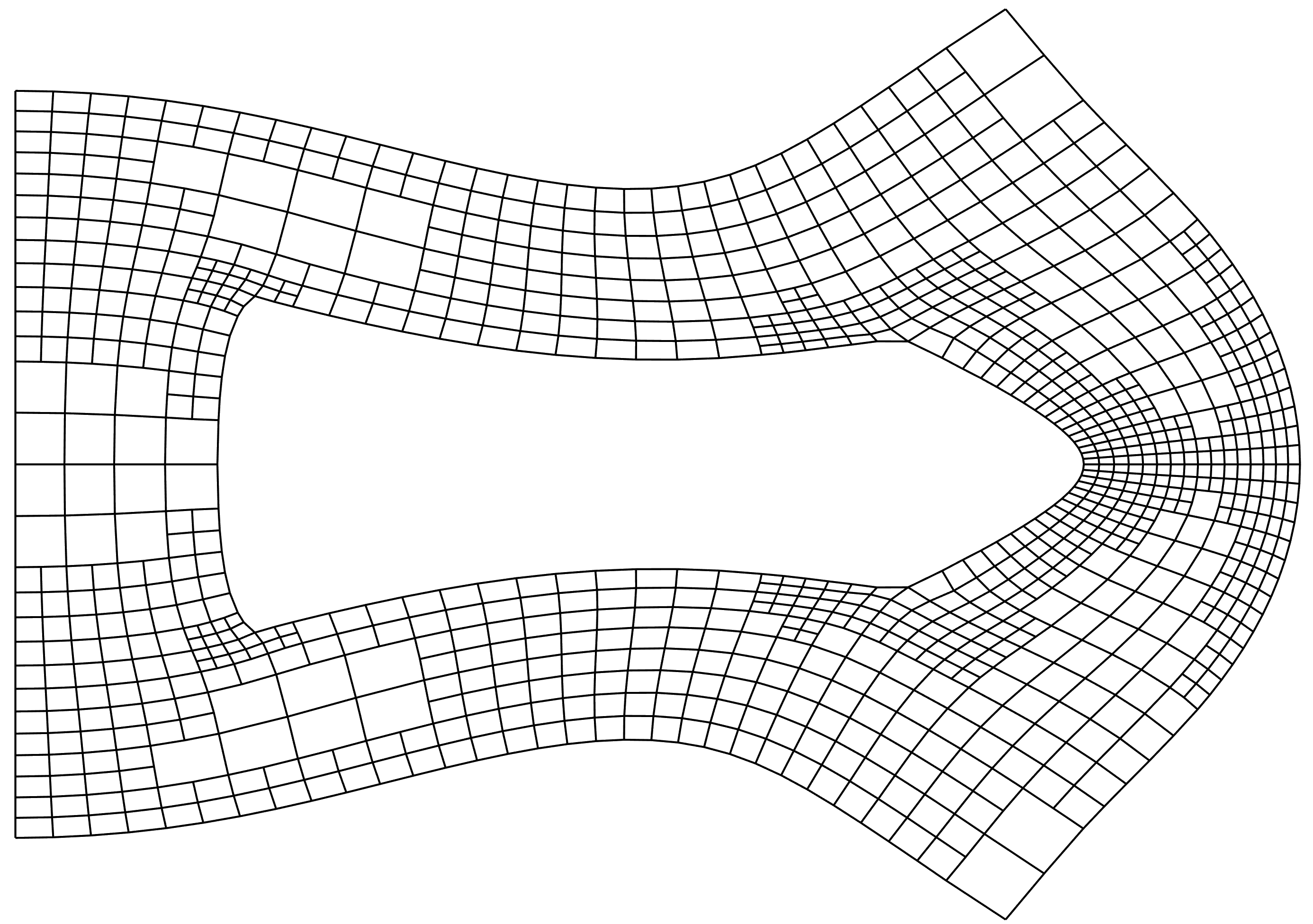}
			\caption{Structured - Step 8}
		\end{subfigure}
		\vskip \baselineskip 
		\begin{subfigure}[t]{0.33\textwidth}
			\centering
			\includegraphics[width=0.95\textwidth]{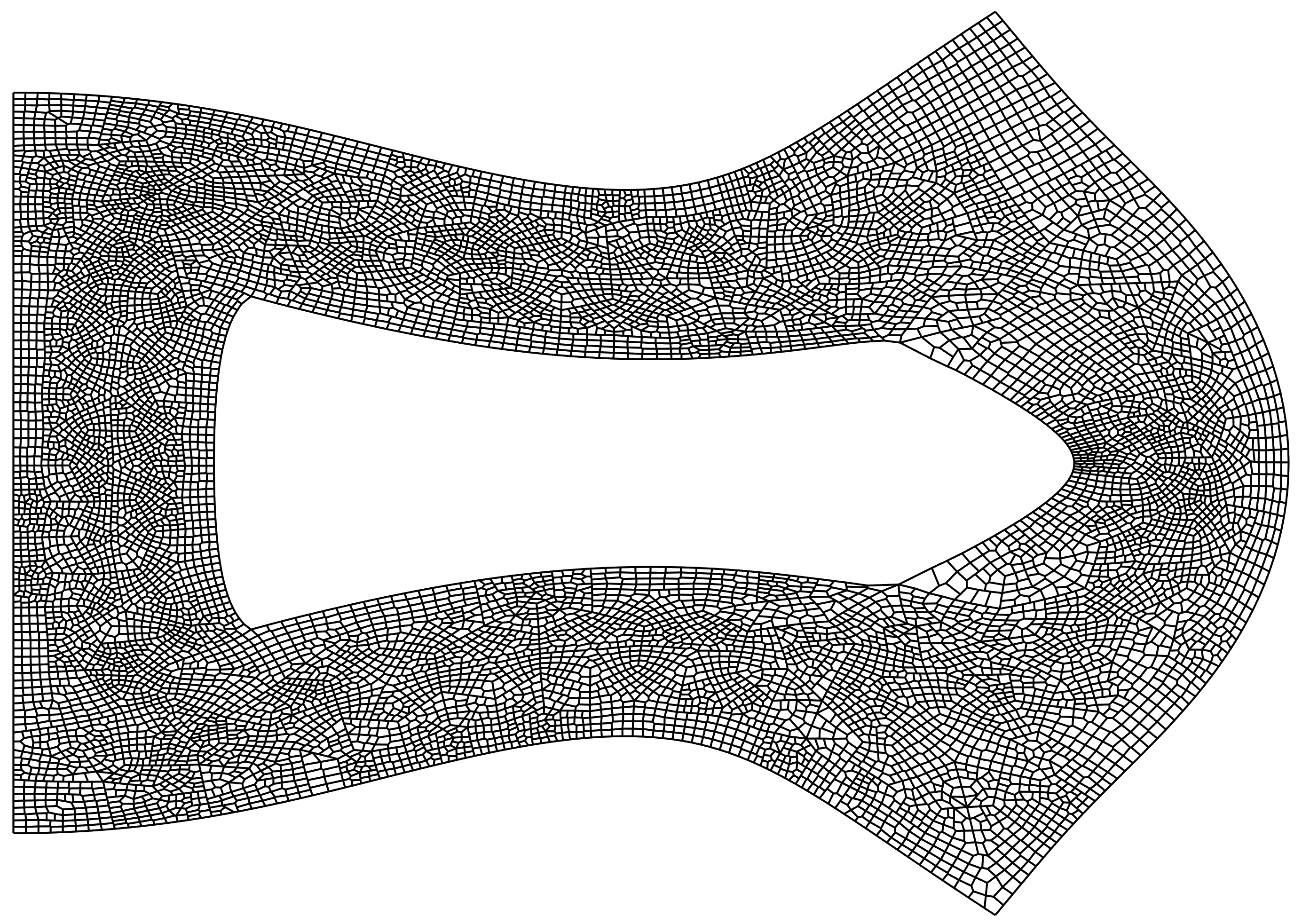}
			\caption{Voronoi - Step 1}
		\end{subfigure}%
		\begin{subfigure}[t]{0.33\textwidth}
			\centering
			\includegraphics[width=0.95\textwidth]{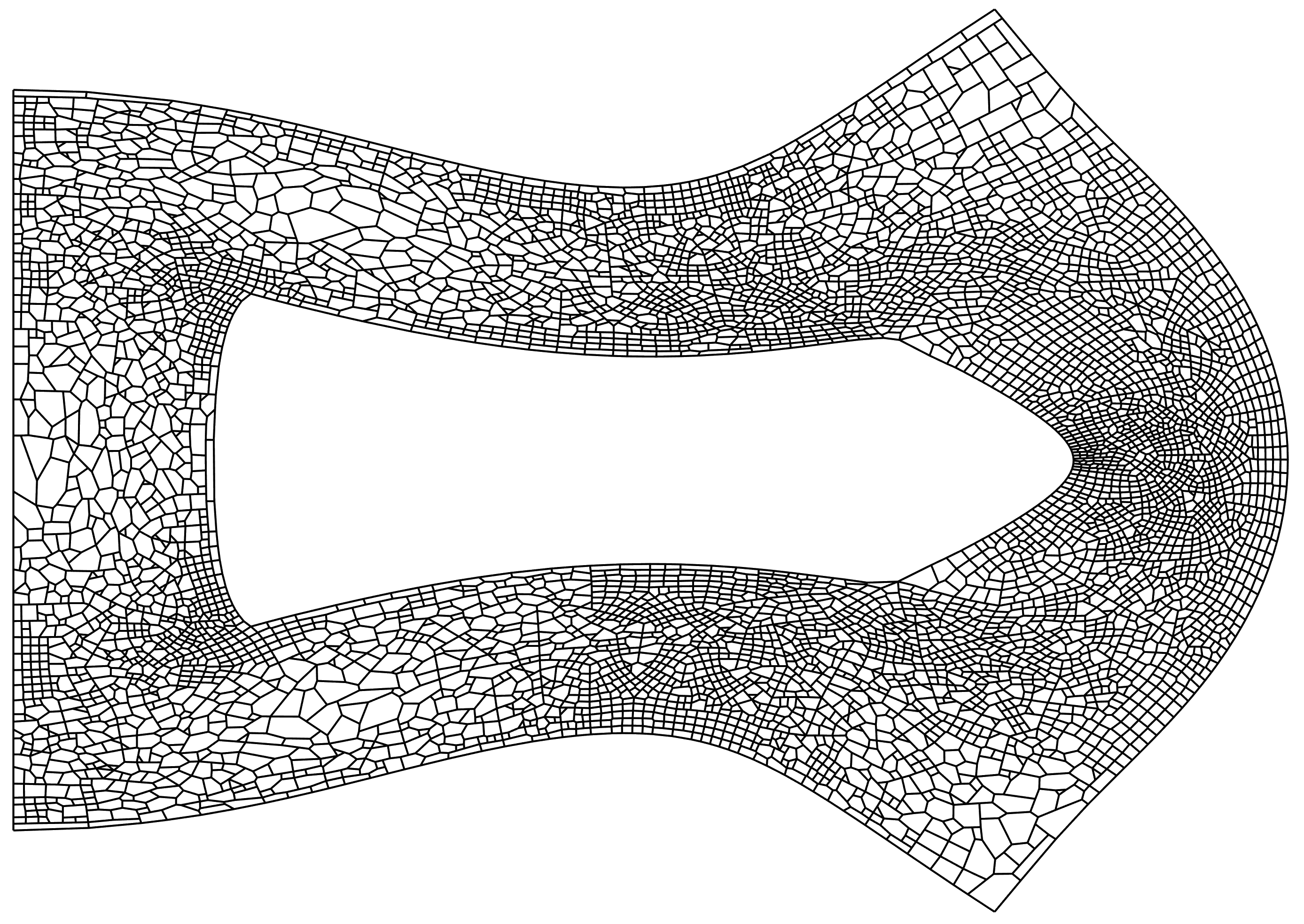}
			\caption{Voronoi - Step 6}
		\end{subfigure}%
		\begin{subfigure}[t]{0.33\textwidth}
			\centering
			\includegraphics[width=0.95\textwidth]{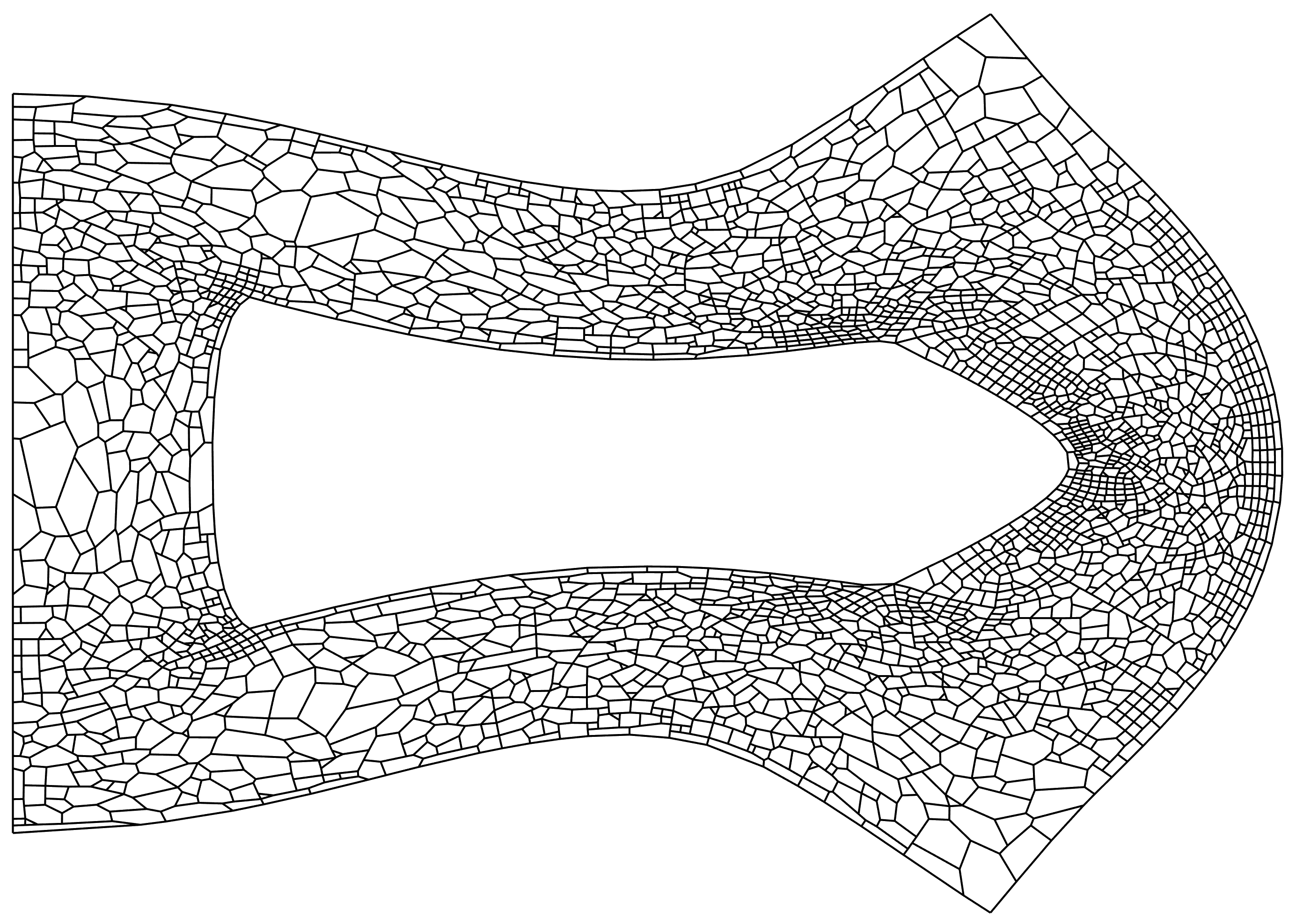}
			\caption{Voronoi - Step 12}
		\end{subfigure}
		\caption{Mesh coarsening process for the plate with hole problem on structured and Voronoi meshes using the energy error-based coarsening procedure with ${T=20\%}$.
			\label{fig:PlateWithHoleMeshes}}
	\end{figure}

	The mesh evolution during the coarsening process for the plate with hole problem is depicted in Figure~\ref{fig:PlateWithHoleMeshes} for the case of the energy error-based coarsening procedure with ${T=20\%}$ on structured and Voronoi meshes.
	Meshes are shown at various coarsening steps with step~1 corresponding to the initial mesh.
	The coarsening behaviour is similar on both structured and Voronoi meshes. The mesh remains fine in the regions around the corners of the hole where the deformation is relatively complex. Additionally, the meshes remain fine in the regions undergoing the most severe deformations, i.e the right-hand edge of the hole, which experiences significant compression, and the right-hand edge of the plate, which experiences significant tension.
	Furthermore, since the right-hand portion of the domain experiences significantly more deformation than the left-hand side, the regions around the right-hand corners of the hole experience higher stresses and more complex deformation than the regions around the left-hand corners of the hole. This behaviour is reflected in the coarsening process with finer meshes preserved in a larger region surrounding the right-hand corners of the hole than the left-hand corners.
	Additionally, meshes are coarsened the most in regions experiencing relatively little and/or simple deformation such as the middle of the left-hand portion of the plate, and the top and bottom right-hand corners of the plate.
	Thus, the mesh evolution is sensible for this problem. 
	
	\FloatBarrier
	\FloatBarrier
	
\begin{figure}[ht!]
		\centering
		\begin{subfigure}[t]{0.33\textwidth}
			\centering
			\includegraphics[width=0.95\textwidth]{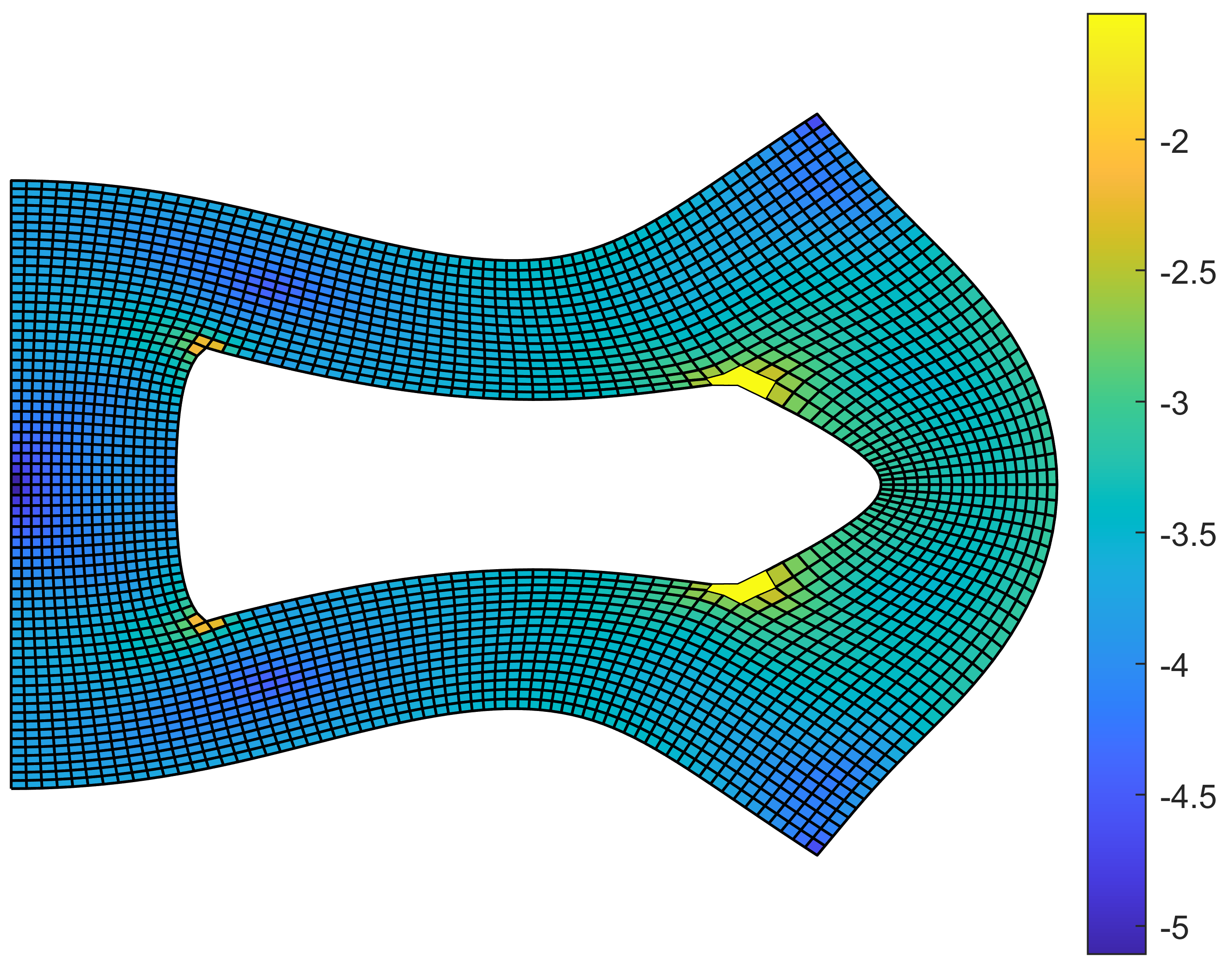}
			\caption{Structured - Step 1}
		\end{subfigure}%
		\begin{subfigure}[t]{0.33\textwidth}
			\centering
			\includegraphics[width=0.95\textwidth]{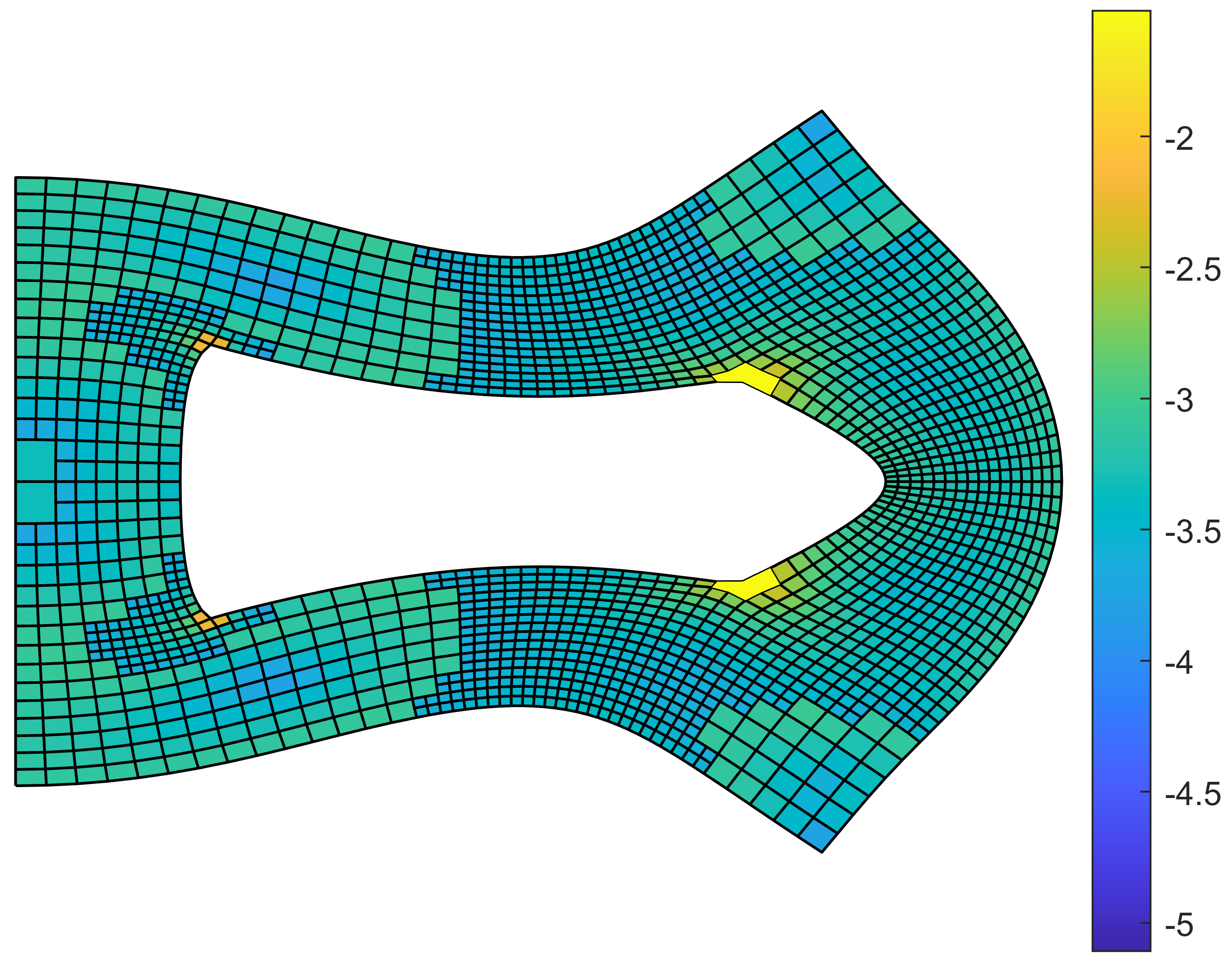}
			\caption{Structured - Step 4}
		\end{subfigure}%
		\begin{subfigure}[t]{0.33\textwidth}
			\centering
			\includegraphics[width=0.95\textwidth]{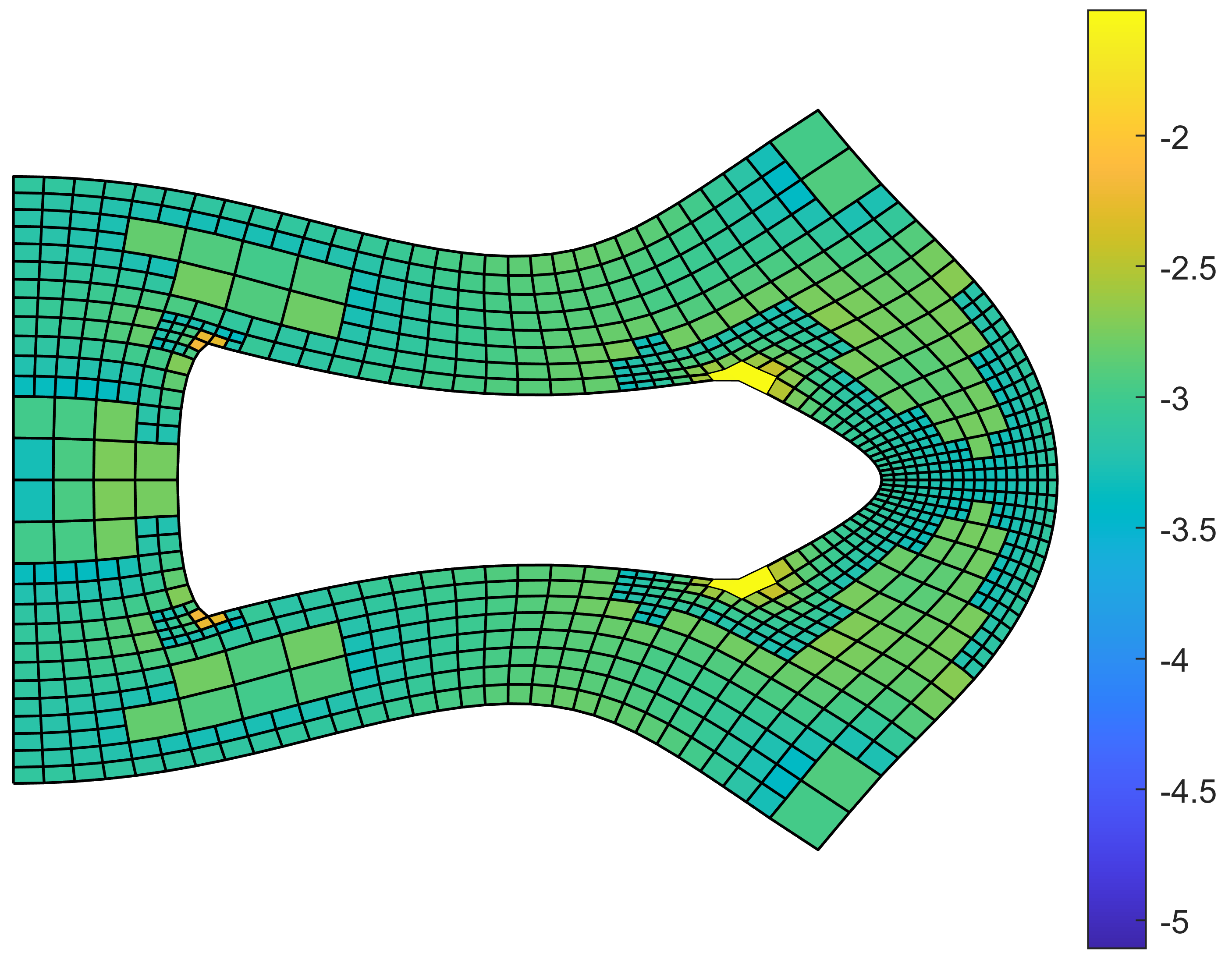}
			\caption{Structured - Step 8}
		\end{subfigure}
		\vskip \baselineskip 
		\begin{subfigure}[t]{0.33\textwidth}
			\centering
			\includegraphics[width=0.95\textwidth]{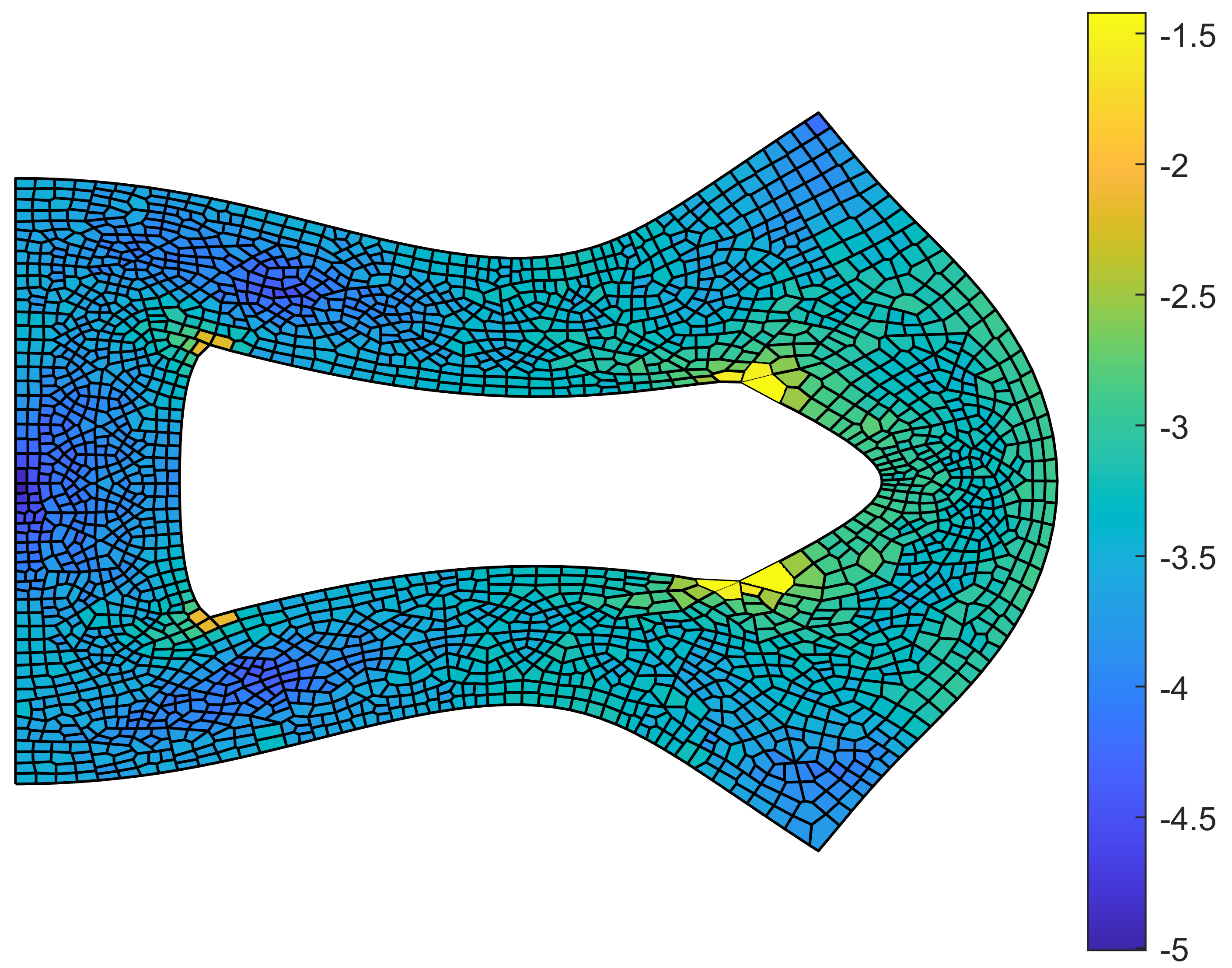}
			\caption{Voronoi - Step 1}
		\end{subfigure}%
		\begin{subfigure}[t]{0.33\textwidth}
			\centering
			\includegraphics[width=0.95\textwidth]{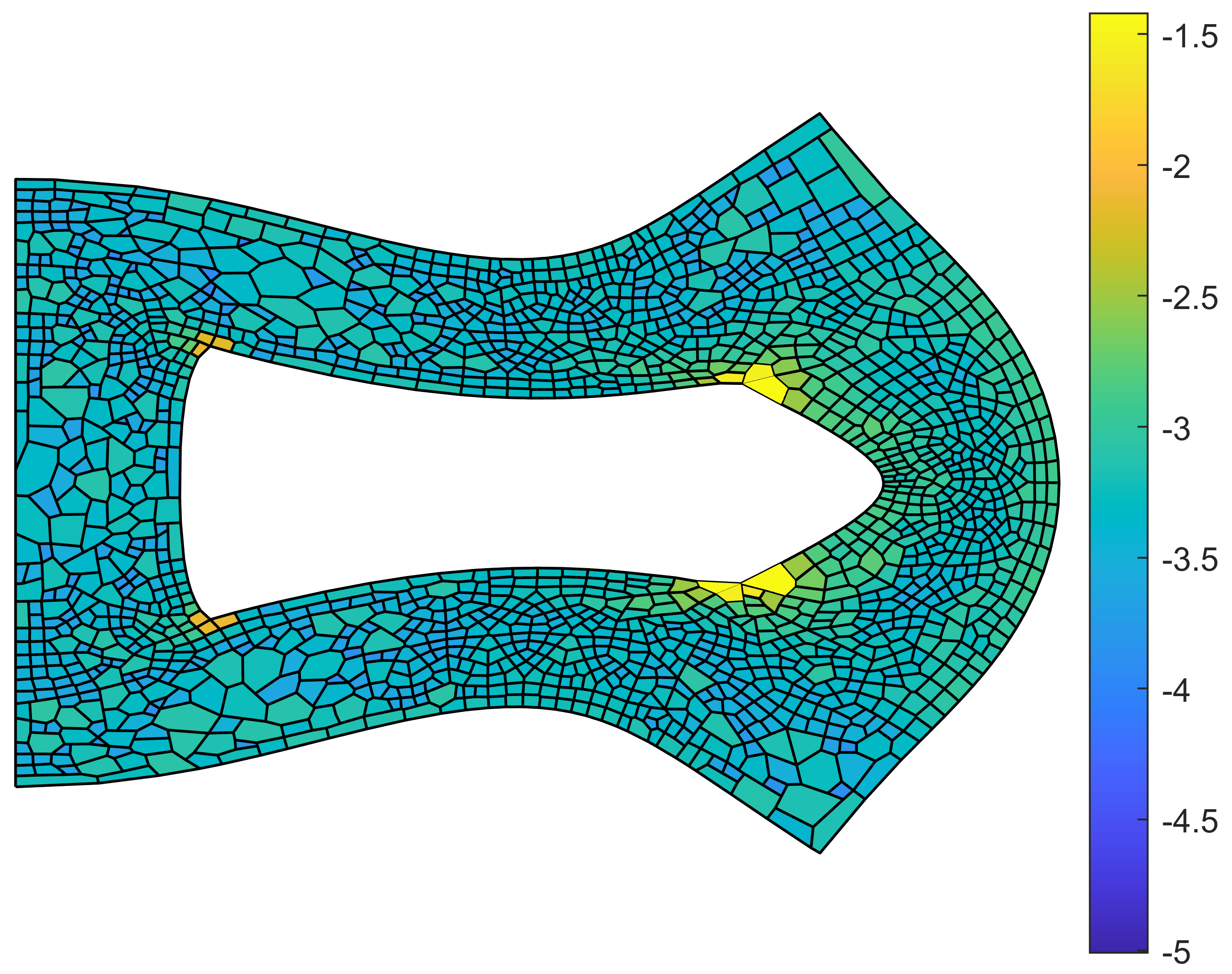}
			\caption{Voronoi - Step 4}
		\end{subfigure}%
		\begin{subfigure}[t]{0.33\textwidth}
			\centering
			\includegraphics[width=0.95\textwidth]{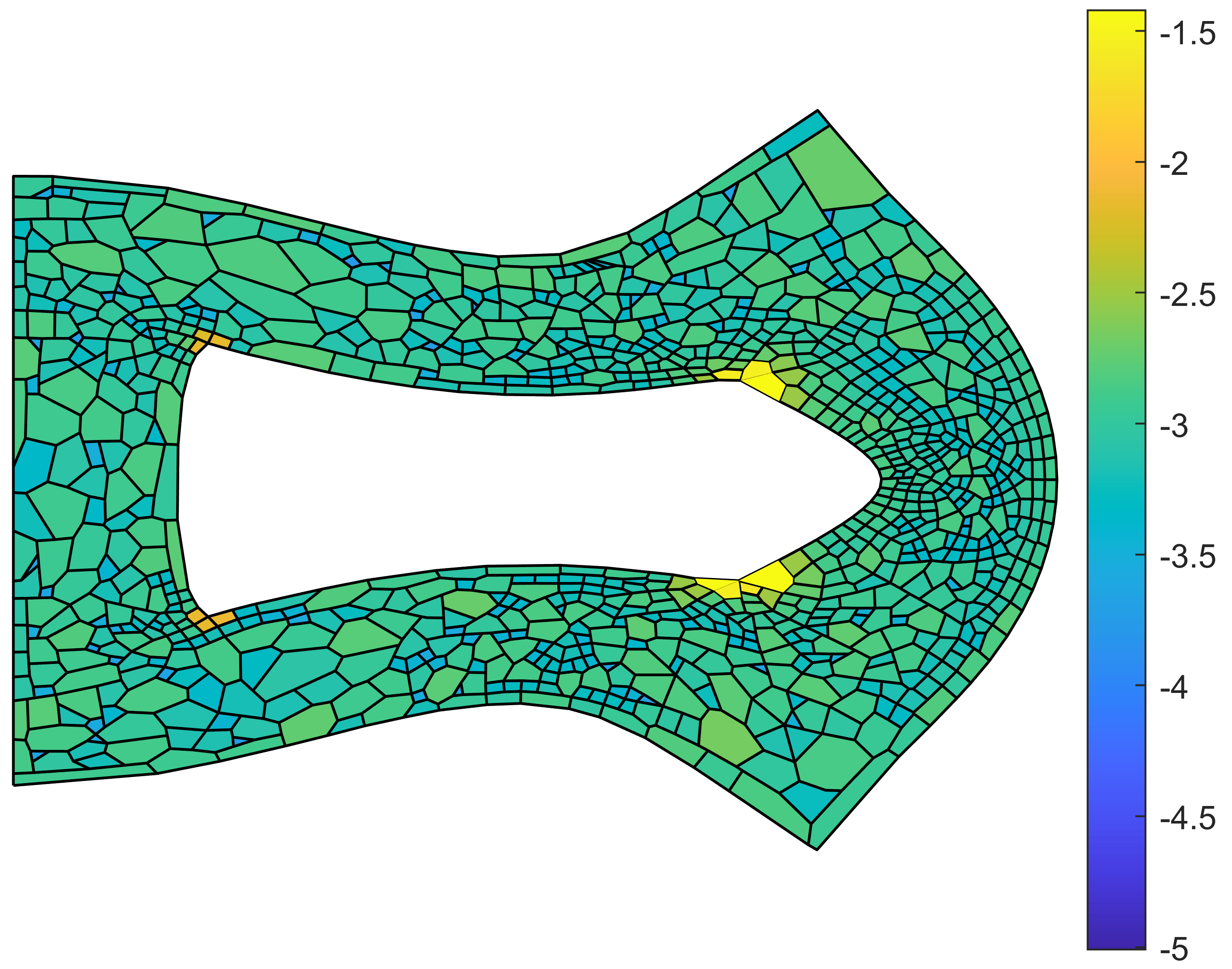}
			\caption{Voronoi - Step 8}
		\end{subfigure}
		\caption{${\mathcal{H}^{1}}$ error distribution during the coarsening process for the plate with hole problem on structured and Voronoi meshes using the energy error-based coarsening procedure with ${T=20\%}$.
			\label{fig:PlateWithHoleErrorMeshes}}
	\end{figure}

	The distribution of the ${\mathcal{H}^{1}}$ error over the domain during the mesh coarsening process for the plate with hole problem is depicted in Figure~\ref{fig:PlateWithHoleErrorMeshes}. The ${\mathcal{H}^{1}}$ error is depicted in a logarithmic scale on structured and Voronoi meshes for the case of the energy error-based coarsening procedure with ${T=20\%}$. 
	As observed in the case of the punch problem, the error distribution exhibited in step~1, i.e. Figure~\ref{fig:PlateWithHoleErrorMeshes}(a), demonstrates that the mesh evolution illustrated in Figure~\ref{fig:PlateWithHoleMeshes} is sensible. The mesh evolution reflects the error distribution and the regions with the lowest errors are the most coarsened. 
	During the mesh coarsening process the distribution of error, again, becomes more even over the domain. However, in the case of this `more challenging´ problem, the complex regions around the corners of the hole induce much higher and more localized stresses than those of the punch problem. Thus, even though the error is spread more evenly over much of the domain, the complex regions stand out as error hot spots.
	
	\FloatBarrier
	\FloatBarrier
	
\begin{figure}[ht!]
		\centering
		\begin{subfigure}[t]{0.5\textwidth}
			\centering
			\includegraphics[width=0.95\textwidth]{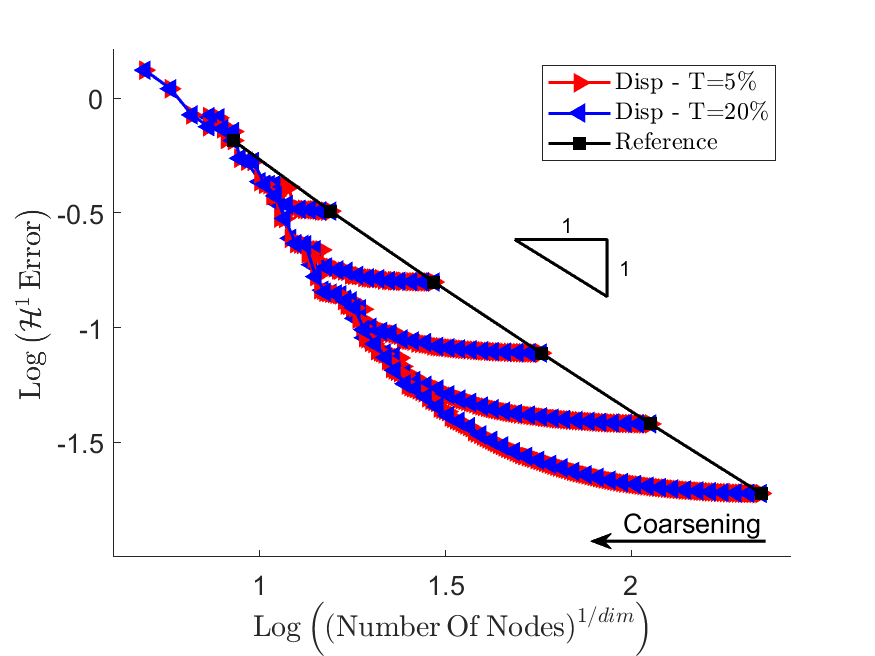}
			\caption{Displacement-based indicator - Structured meshes}
		\end{subfigure}%
		\begin{subfigure}[t]{0.5\textwidth}
			\centering
			\includegraphics[width=0.95\textwidth]{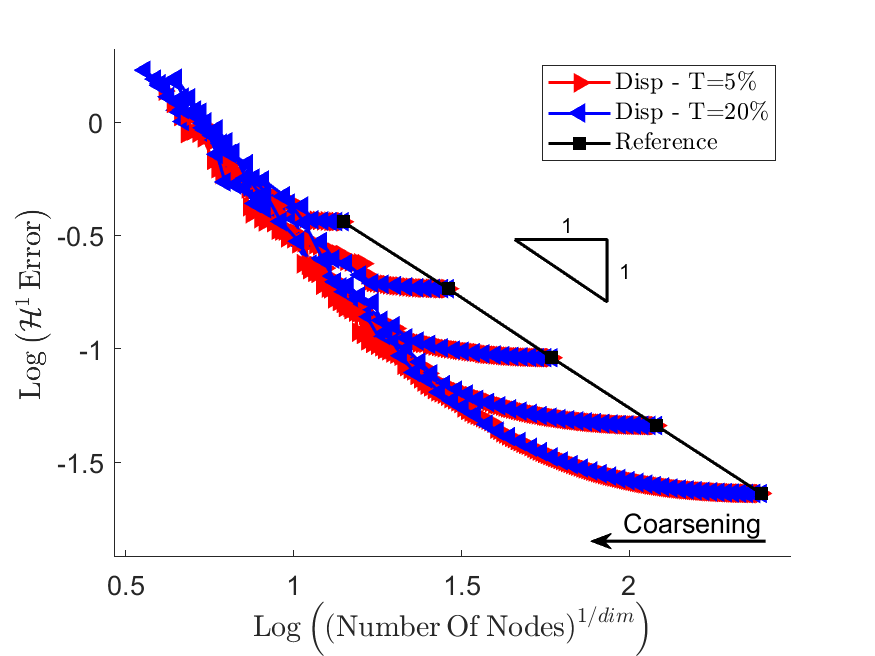}
			\caption{Displacement-based indicator - Voronoi meshes}
		\end{subfigure}
		\vskip \baselineskip 
		\begin{subfigure}[t]{0.5\textwidth}
			\centering
			\includegraphics[width=0.95\textwidth]{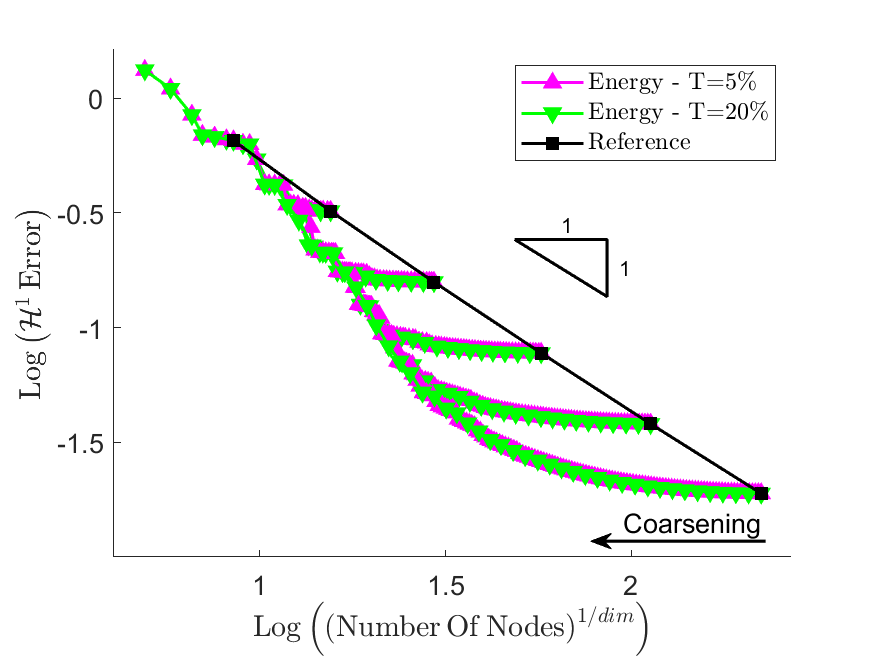}
			\caption{Energy error-based indicator - Structured meshes}
		\end{subfigure}%
		\begin{subfigure}[t]{0.5\textwidth}
			\centering
			\includegraphics[width=0.95\textwidth]{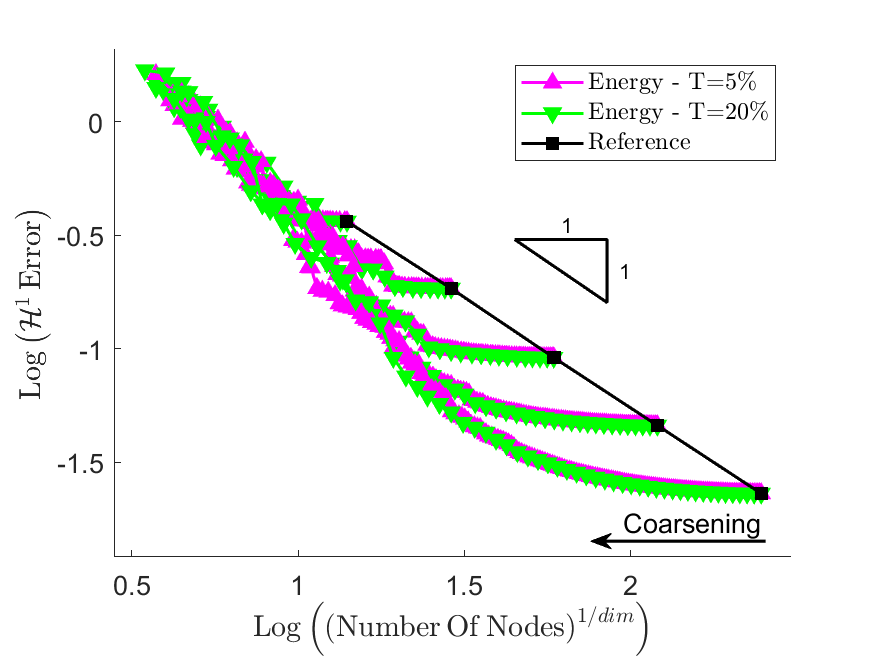}
			\caption{Energy error-based indicator - Voronoi meshes}
		\end{subfigure}
		\caption{$\mathcal{H}^{1}$ error vs $n_{\rm v}$ for the plate with hole problem.
			\label{fig:PlateWithHoleConvergenceNumberOfNodes}}
	\end{figure}

	The convergence behaviour in the ${\mathcal{H}^{1}}$ error norm of the VEM for the plate with hole problem using the displacement-based and energy error-based coarsening procedures is depicted in Figure~\ref{fig:PlateWithHoleConvergenceNumberOfNodes} on a logarithmic scale. Here, the convergence behaviour of the displacement-based and energy error-based procedures are presented on the top and bottom rows of figures respectively, with the results generated on structured and Voronoi meshes presented in the left and right columns of figures respectively. In each case results are presented for coarsening thresholds ${T=5\%}$ and ${T=20\%}$ to demonstrate the effect of the choice of $T$ on the convergence behaviour. Additionally, uniform initial meshes of various discretization densities are considered.
	The behaviour exhibited in Figure~\ref{fig:PlateWithHoleConvergenceNumberOfNodes} is qualitatively similar to that observed in Figure~\ref{fig:PunchConvergenceNumberOfNodes} for the case of the punch problem. The benefit of the coarsening procedures is, again, clear with the solutions obtained from the coarsened meshes exhibiting higher accuracy than those obtained from uniform meshes comprising the same number of degrees of freedom. Furthermore, similar behaviour is exhibited in the cases of both structured and Voronoi meshes, the behaviour does not appear to be significantly influenced by the choice of $T$, and the efficacy of the coarsening procedures increases with increasing density of the initial mesh.
	The most significant difference between the results presented in Figure~\ref{fig:PunchConvergenceNumberOfNodes} and Figure~\ref{fig:PlateWithHoleConvergenceNumberOfNodes} is the increased efficacy of the coarsening procedures exhibited in the case of the plate with hole problem compared to the punch problem. The reason for this difference is the relative complexities of the two problems. The plate with hole problem has significantly more localized complexity around the corners of the hole, while the punch problem is comparatively simpler to model. Thus, in the plate with hole problem most of the error is localized around the corners of the hole where the fine mesh is preserved. The rest of the domain can be efficiently coarsened while introducing very little error. Thus, this `more challenging´ problem and its more localized complexity is better suited to adaptive coarsening and the proposed coarsening procedures exhibit greater efficacy.
	
	\FloatBarrier
	\FloatBarrier

	\subsection{L-shaped domain}
	\label{subsec:LDomain}
\begin{figure}[ht!]
		\centering
		\begin{subfigure}[t]{0.45\textwidth}
			\centering
			\includegraphics[width=0.95\textwidth]{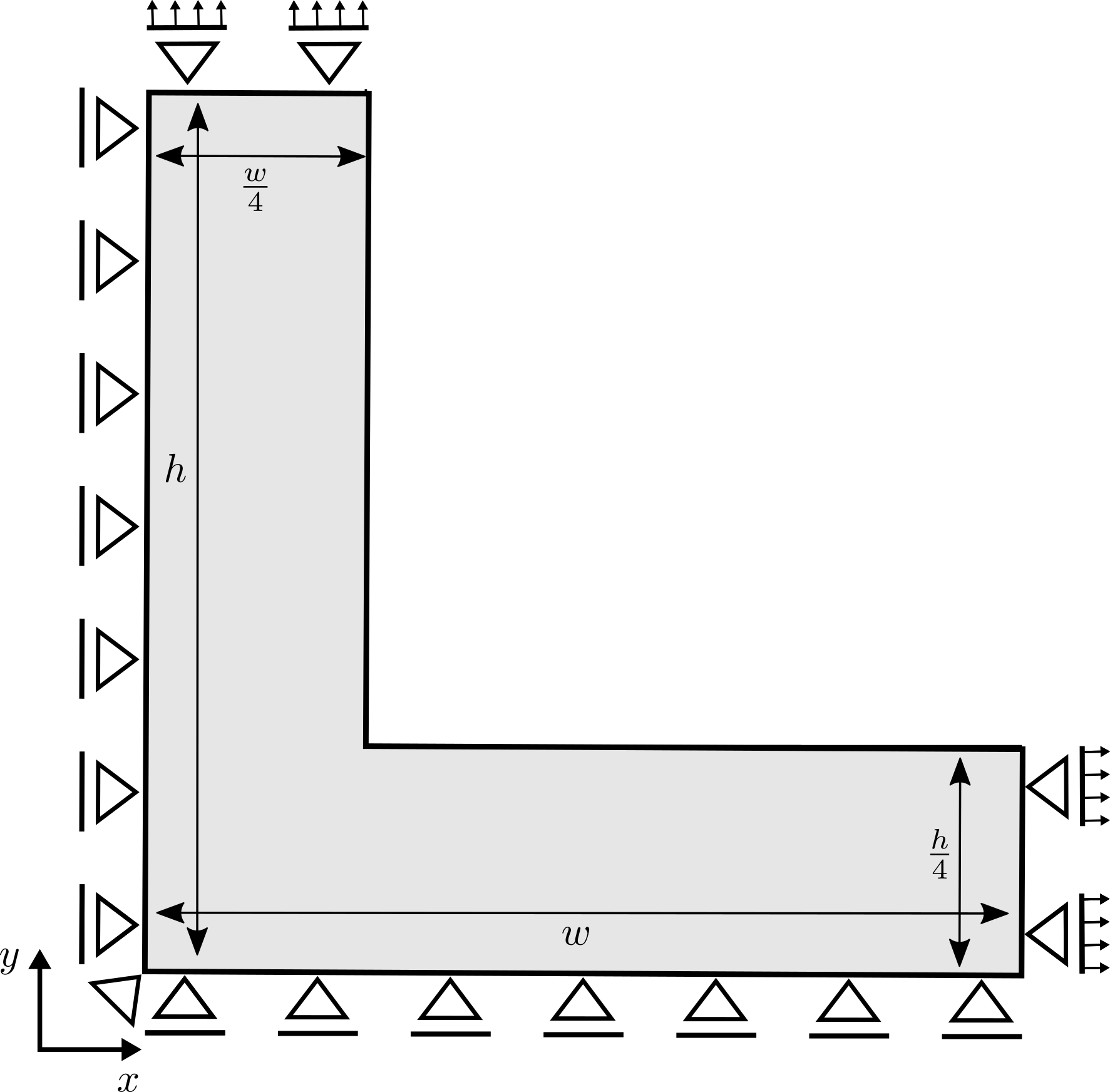}
			\caption{Problem geometry}
		\end{subfigure}%
		\begin{subfigure}[t]{0.55\textwidth}
			\centering
			\includegraphics[width=0.85\textwidth]{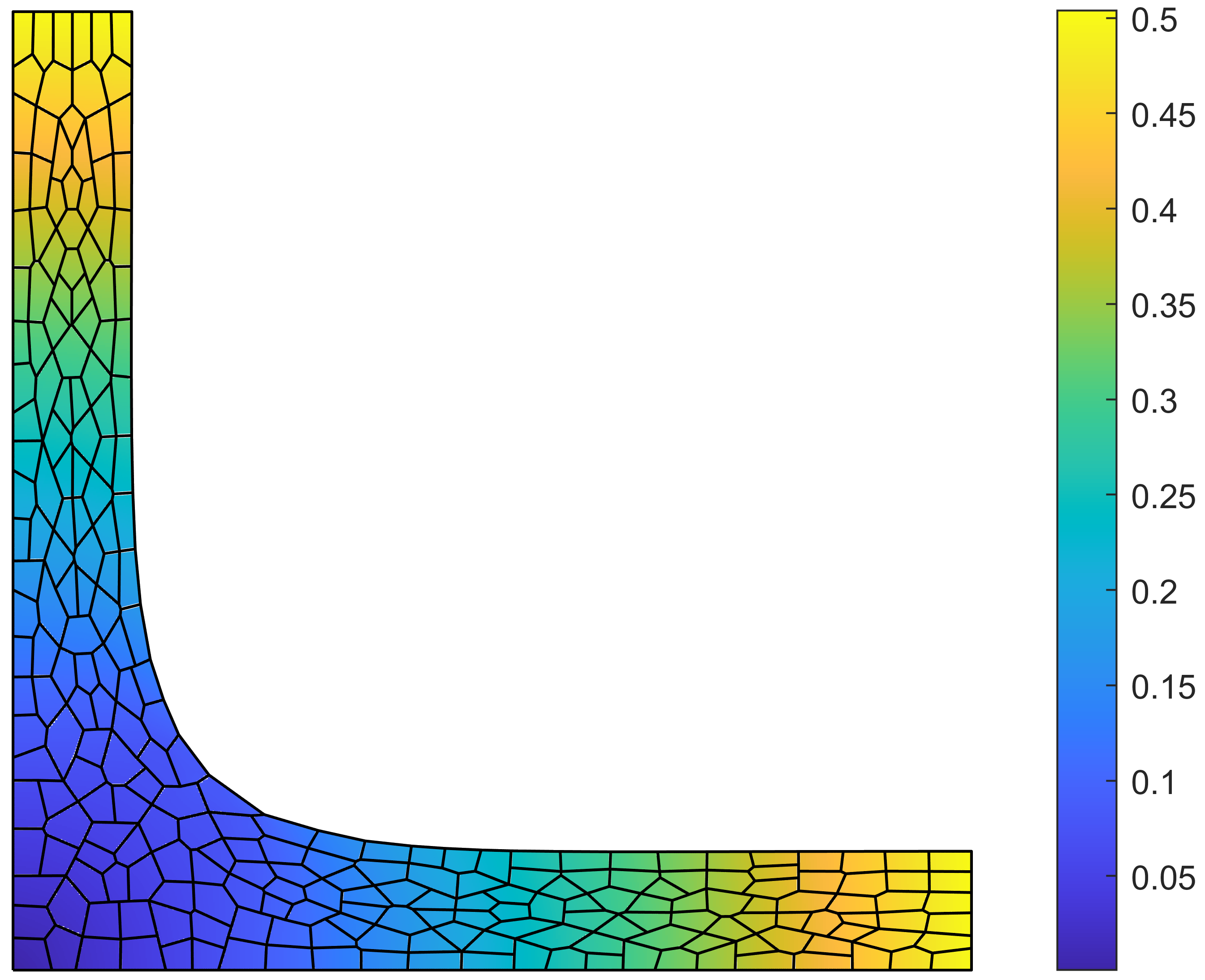}
			\caption{Deformed configuration}
		\end{subfigure}
		\caption{L-shaped domain (a) geometry and (b) sample deformed configuration.
			\label{fig:LShapedDomainGeometry}}
	\end{figure}

	The L-shaped domain problem comprises a domain of width ${w=1~\rm{m}}$ and height ${h=1~\rm{m}}$ where the horizontal and vertical thickness of the L are ${\frac{w}{4}}$ and ${\frac{h}{4}}$ respectively. The bottom and left-hand edges of the domain are constrained vertically and horizontally respectively, with the bottom left corner fully constrained. The upper and right-hand edges are subject to prescribed displacements of ${\bar{u}_{y}=0.5~\rm{m}}$ and ${\bar{u}_{x}=0.5~\rm{m}}$ respectively, with the displacements of the edges unconstrained in the $x$- and $y$-directions respectively (see Figure~\ref{fig:LShapedDomainGeometry}(a)). A sample deformed configuration of the L with a Voronoi mesh is depicted in Figure~\ref{fig:LShapedDomainGeometry}(b) with the displacement magnitude $|\bu|$ plotted on the colour axis. The non-convex corner of the L-shaped geometry introduces a strong singularity that provides a significant challenge for numerical analysis techniques. Thus, this problem is used to provide insight into the efficacy of the proposed coarsening procedures in cases of `very challenging´ problems. 
	
	\FloatBarrier
	\FloatBarrier
	
\begin{figure}[ht!]
		\centering
		\begin{subfigure}[t]{0.33\textwidth}
			\centering
			\includegraphics[width=0.95\textwidth,height=0.95\textwidth]{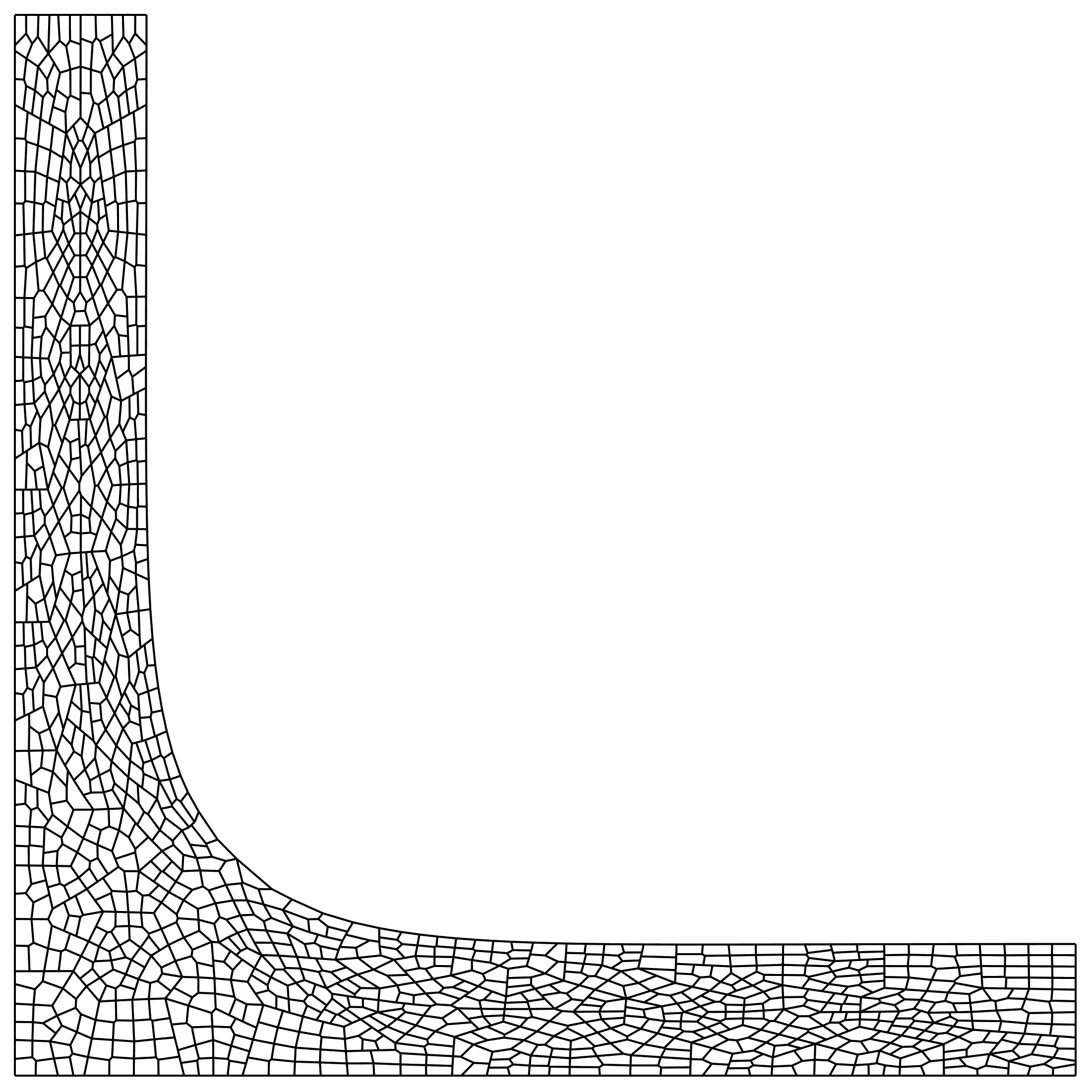}
			\caption{Displacement-based - Step 1}
		\end{subfigure}%
		\begin{subfigure}[t]{0.33\textwidth}
			\centering
			\includegraphics[width=0.95\textwidth,height=0.95\textwidth]{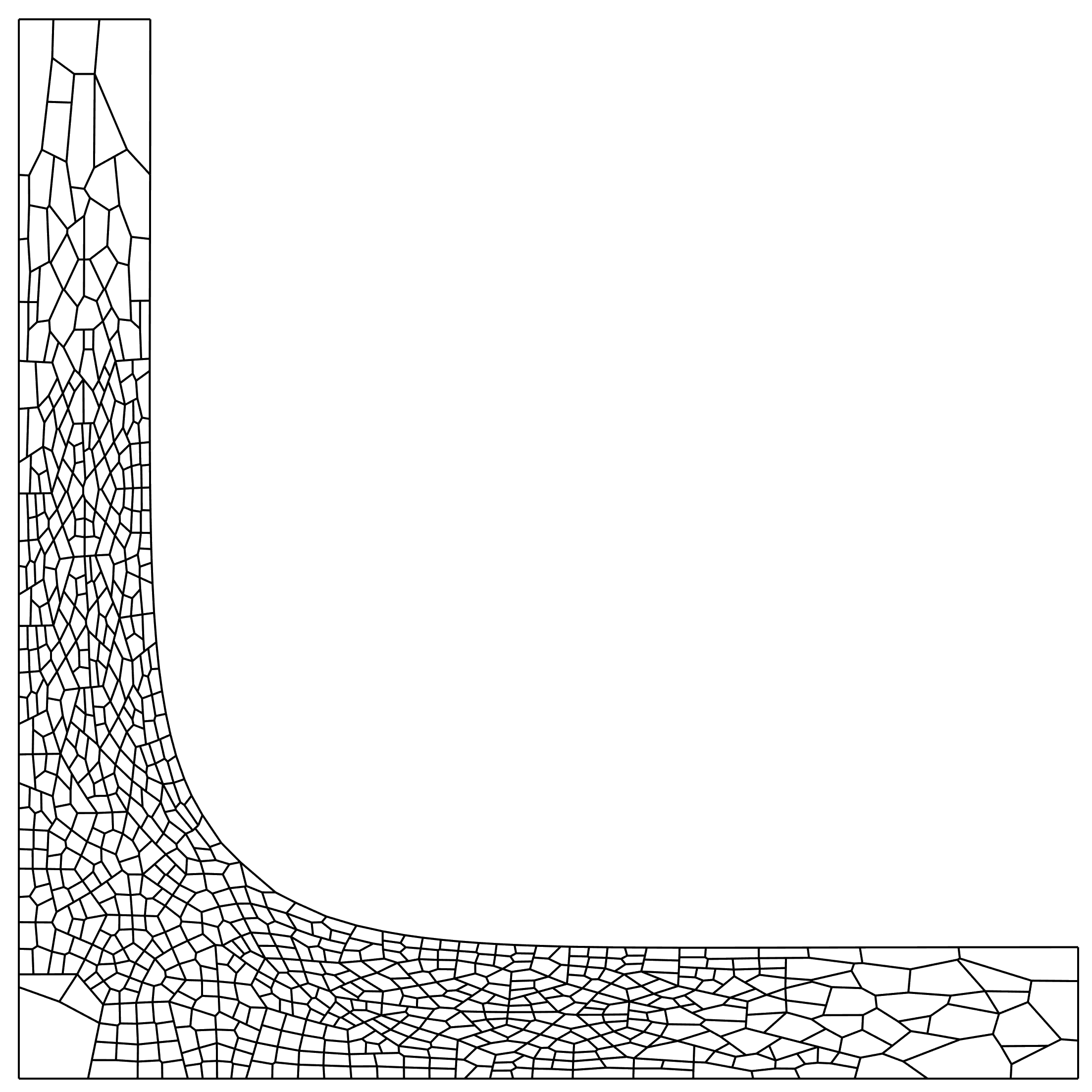}
			\caption{Displacement-based - Step 5}
		\end{subfigure}%
		\begin{subfigure}[t]{0.33\textwidth}
			\centering
			\includegraphics[width=0.95\textwidth,height=0.95\textwidth]{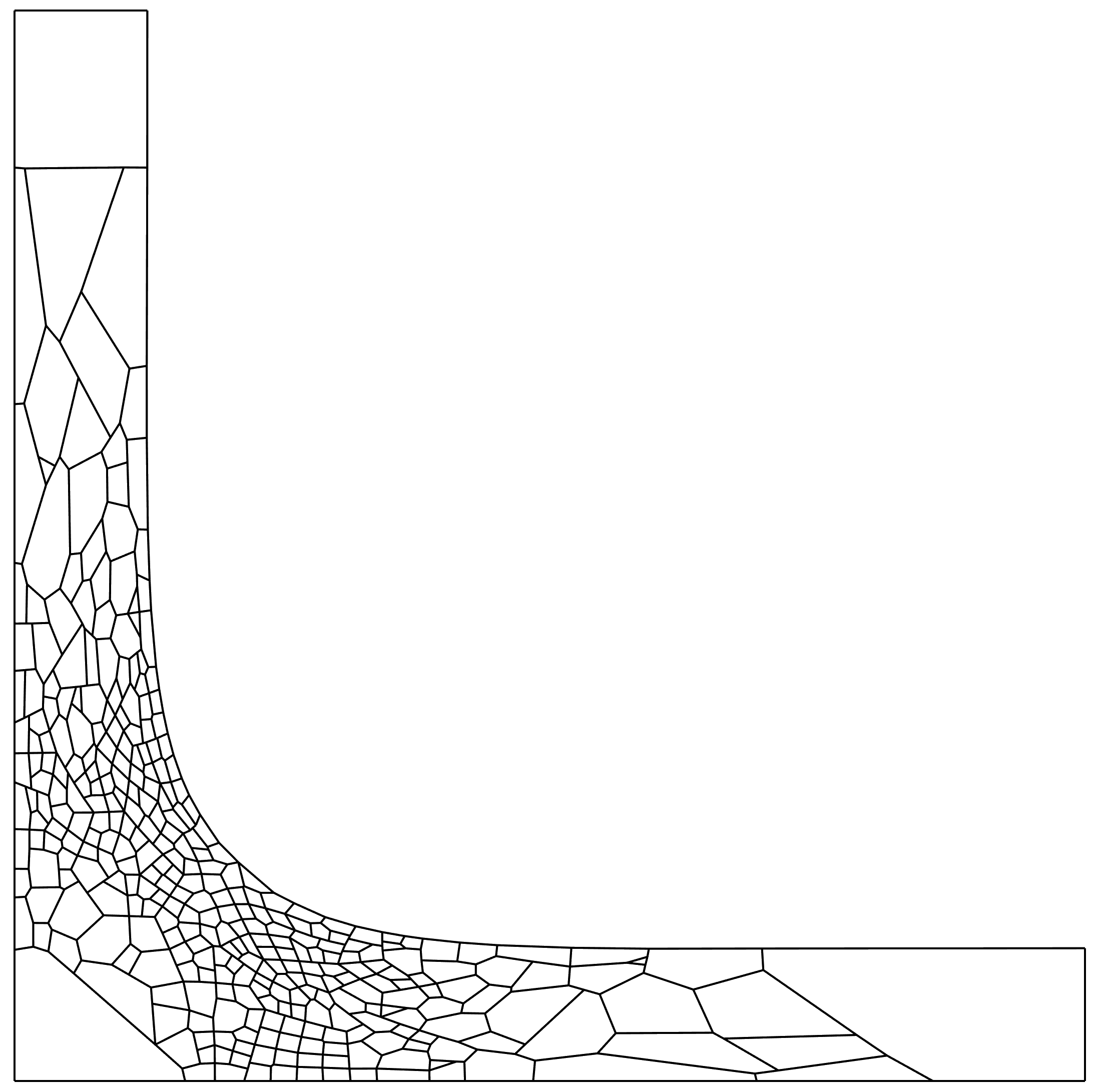}
			\caption{Displacement-based - Step 12}
		\end{subfigure}
		\vskip \baselineskip 
		\begin{subfigure}[t]{0.33\textwidth}
			\centering
			\includegraphics[width=0.95\textwidth,height=0.95\textwidth]{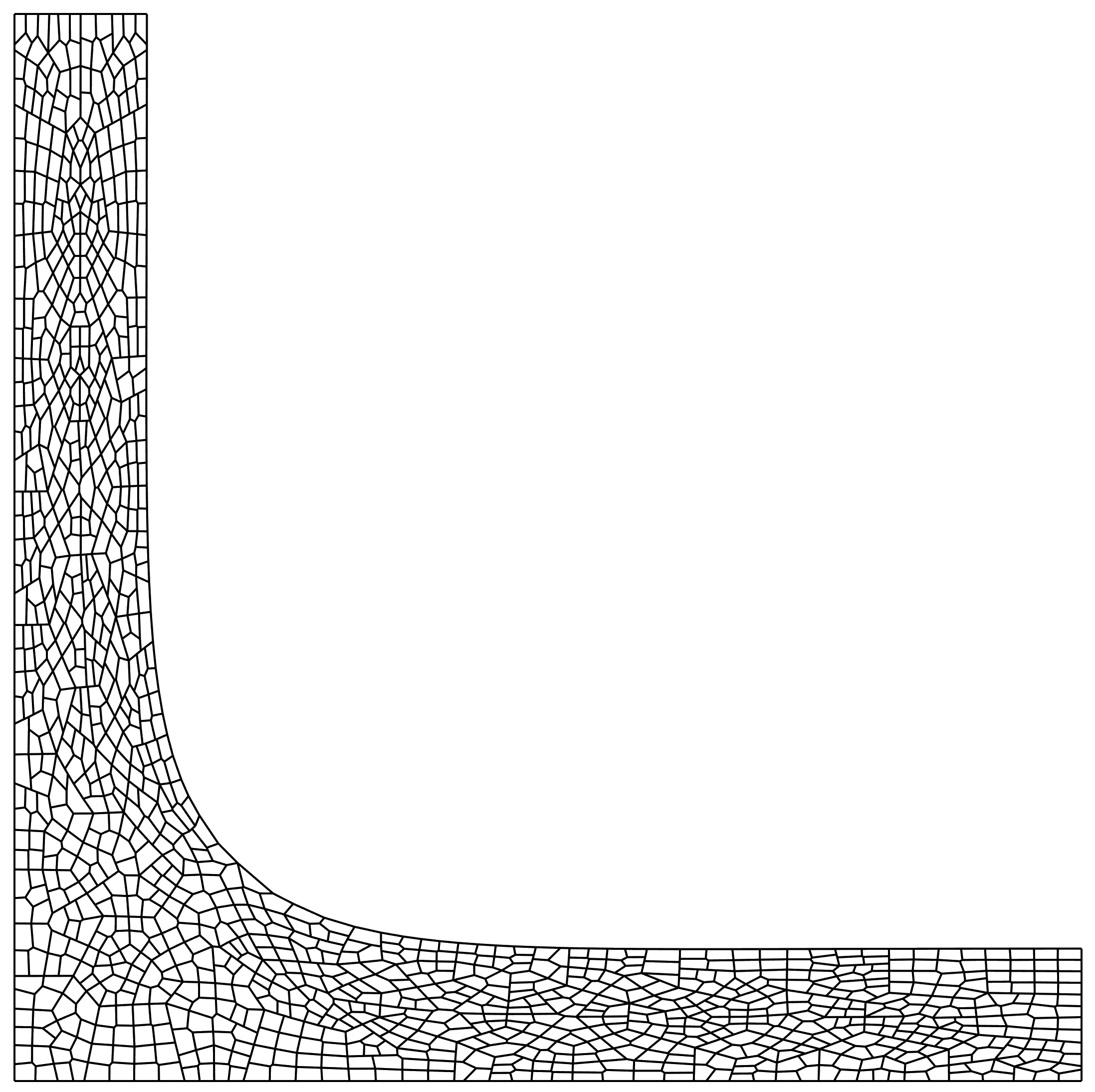}
			\caption{Energy error-based - Step 1}
		\end{subfigure}%
		\begin{subfigure}[t]{0.33\textwidth}
			\centering
			\includegraphics[width=0.95\textwidth,height=0.95\textwidth]{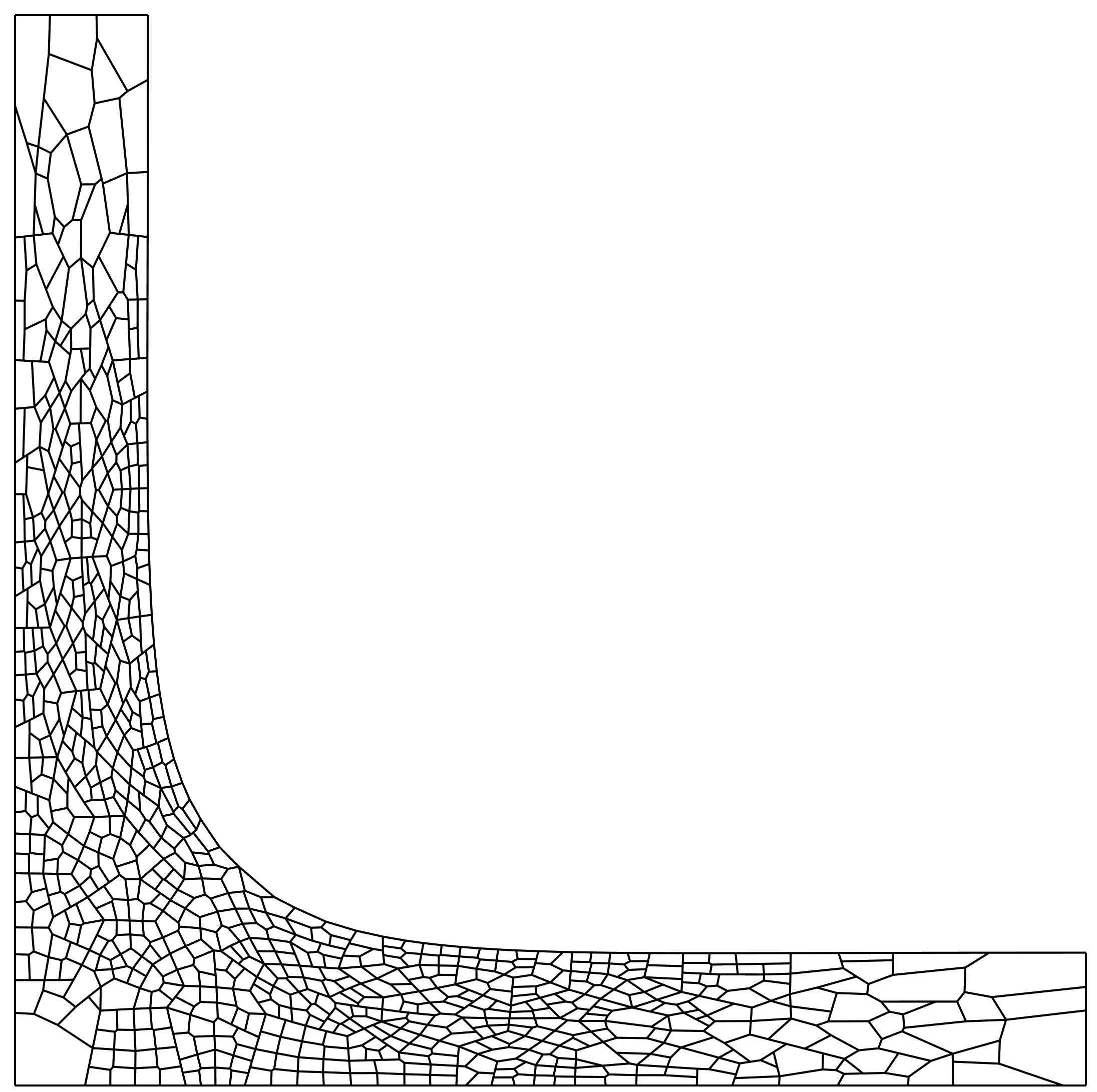}
			\caption{Energy error-based - Step 5}
		\end{subfigure}%
		\begin{subfigure}[t]{0.33\textwidth}
			\centering
			\includegraphics[width=0.95\textwidth,height=0.95\textwidth]{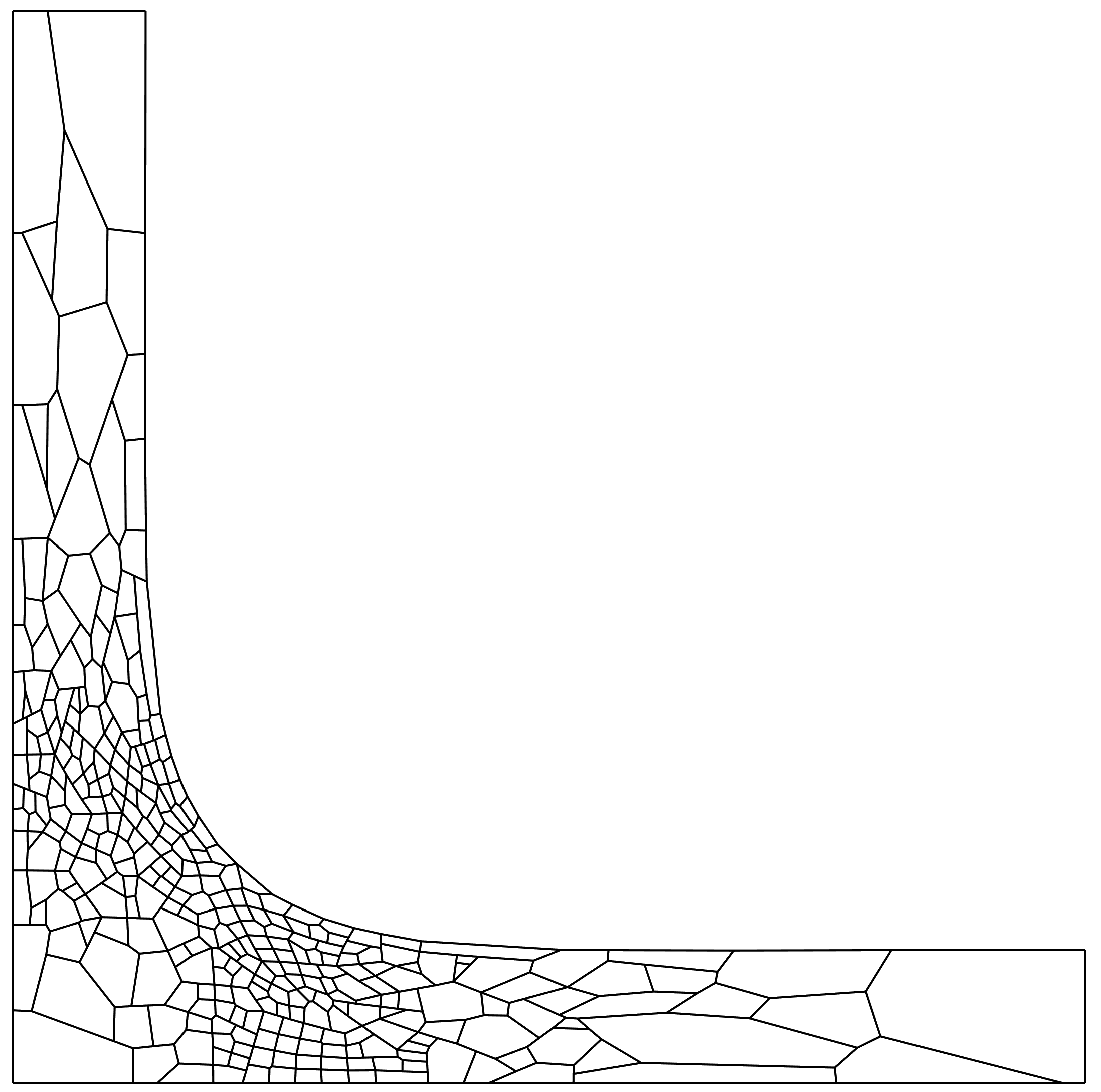}
			\caption{Energy error-based - Step 12}
		\end{subfigure}
		\caption{Mesh coarsening process for the L-shaped domain problem on Voronoi meshes for the displacement-based and energy error-based coarsening procedures with ${T=20\%}$.
			\label{fig:LDomainMeshes}}
	\end{figure}

	The mesh evolution during the coarsening process for the L-shaped domain problem is depicted in Figure~\ref{fig:LDomainMeshes} for the cases of the displacement-based and energy error-based coarsening procedures with ${T=20\%}$ on  Voronoi meshes.
	Meshes are shown at various coarsening steps with step~1 corresponding to the initial mesh.
	The coarsening behaviour is similar for both the displacement-based and energy-error based coarsening procedures. The mesh remains fine around the corner of the L. This is expected as the corner of the L induces a singularity in the solution and, as such, is the most complex region of the domain and requires the finest discretization possible. Conversely, the rest of the domain experiences very simple, almost linear, deformation, particularly in the regions furthest from the corner. Thus, a very low discretization density is required in these regions and they are increasingly coarsened as the number of coarsening steps/iterations performed increases.
	This dichotomy between fine and coarse mesh regions is expected for this problem, thus, indicating sensible mesh evolution. Furthermore, the meshes presented here generated by coarsening finer uniform initial meshes are qualitatively very similar to the meshes presented in \cite{Huyssteen2022} that were generated using adaptive refinement of initially uniform coarse meshes on the same example problem. This similarity exists for all of the example problems presented in this work and further emphasises the efficacy of the proposed coarsening procedures.
	
	\FloatBarrier
	\FloatBarrier
	
\begin{figure}[ht!]
		\centering
		\begin{subfigure}[t]{0.33\textwidth}
			\centering
			\includegraphics[width=0.95\textwidth]{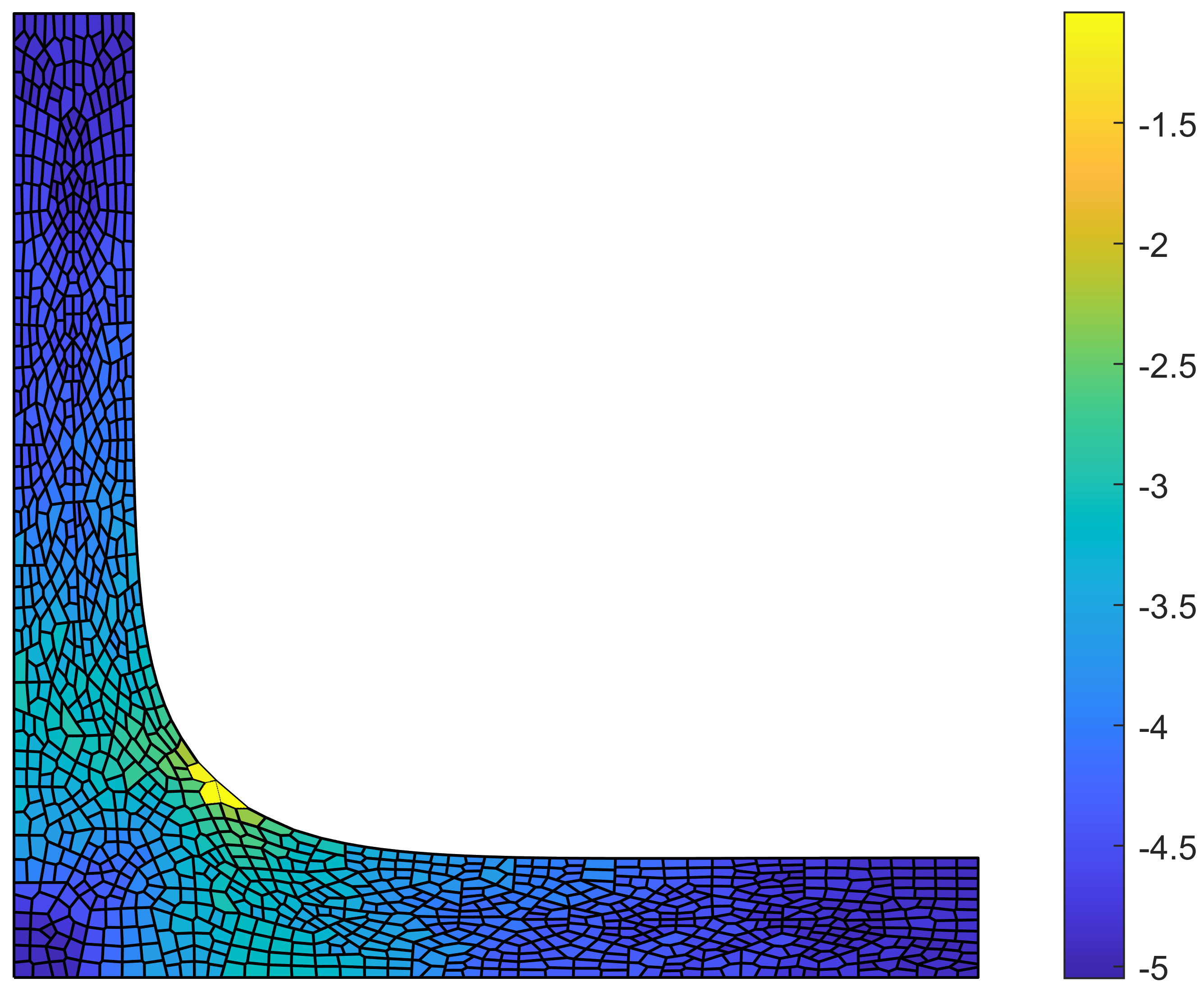}
			\caption{Displacement-based - Step 1}
		\end{subfigure}%
		\begin{subfigure}[t]{0.33\textwidth}
			\centering
			\includegraphics[width=0.95\textwidth]{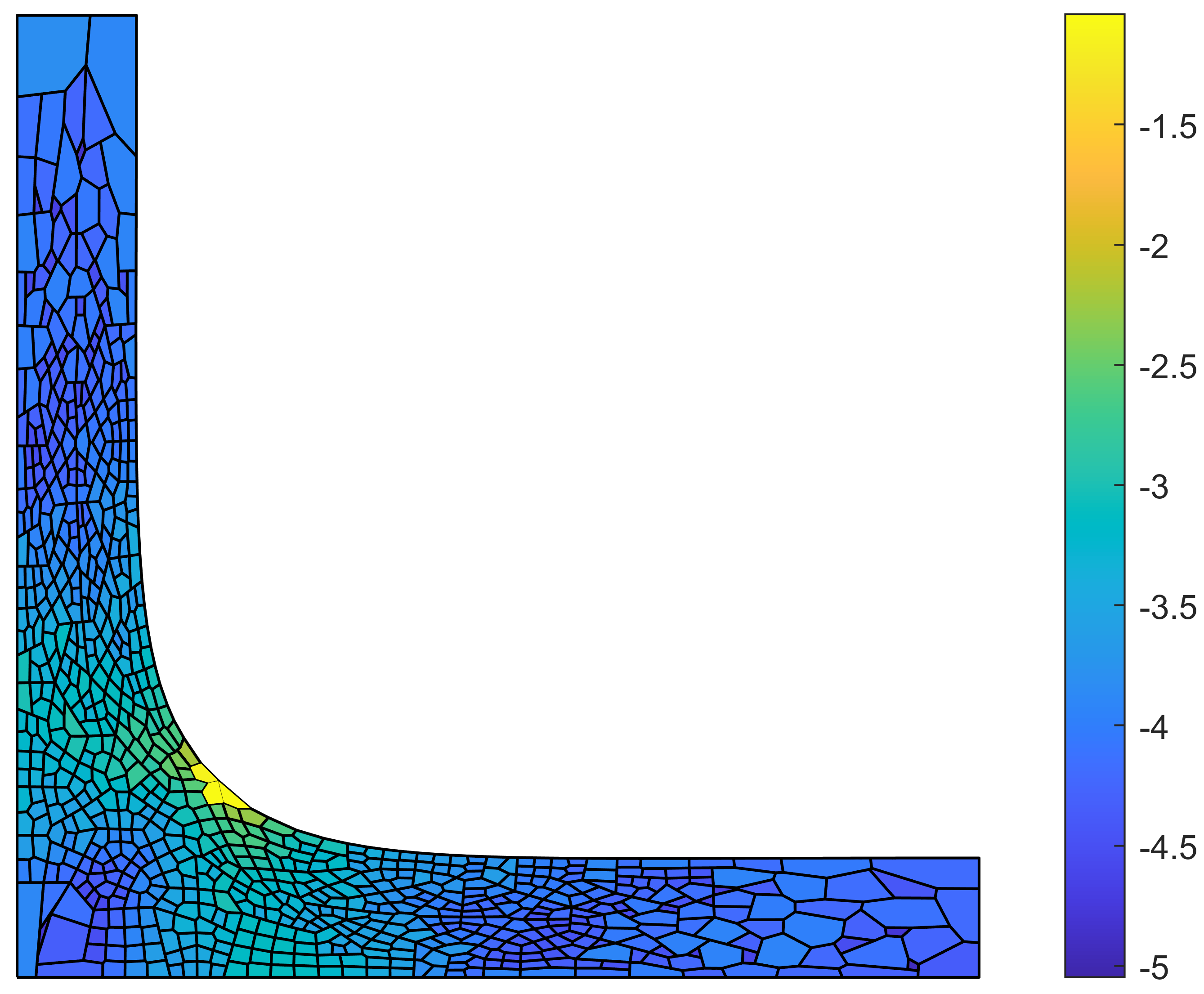}
			\caption{Displacement-based - Step 5}
		\end{subfigure}%
		\begin{subfigure}[t]{0.33\textwidth}
			\centering
			\includegraphics[width=0.95\textwidth]{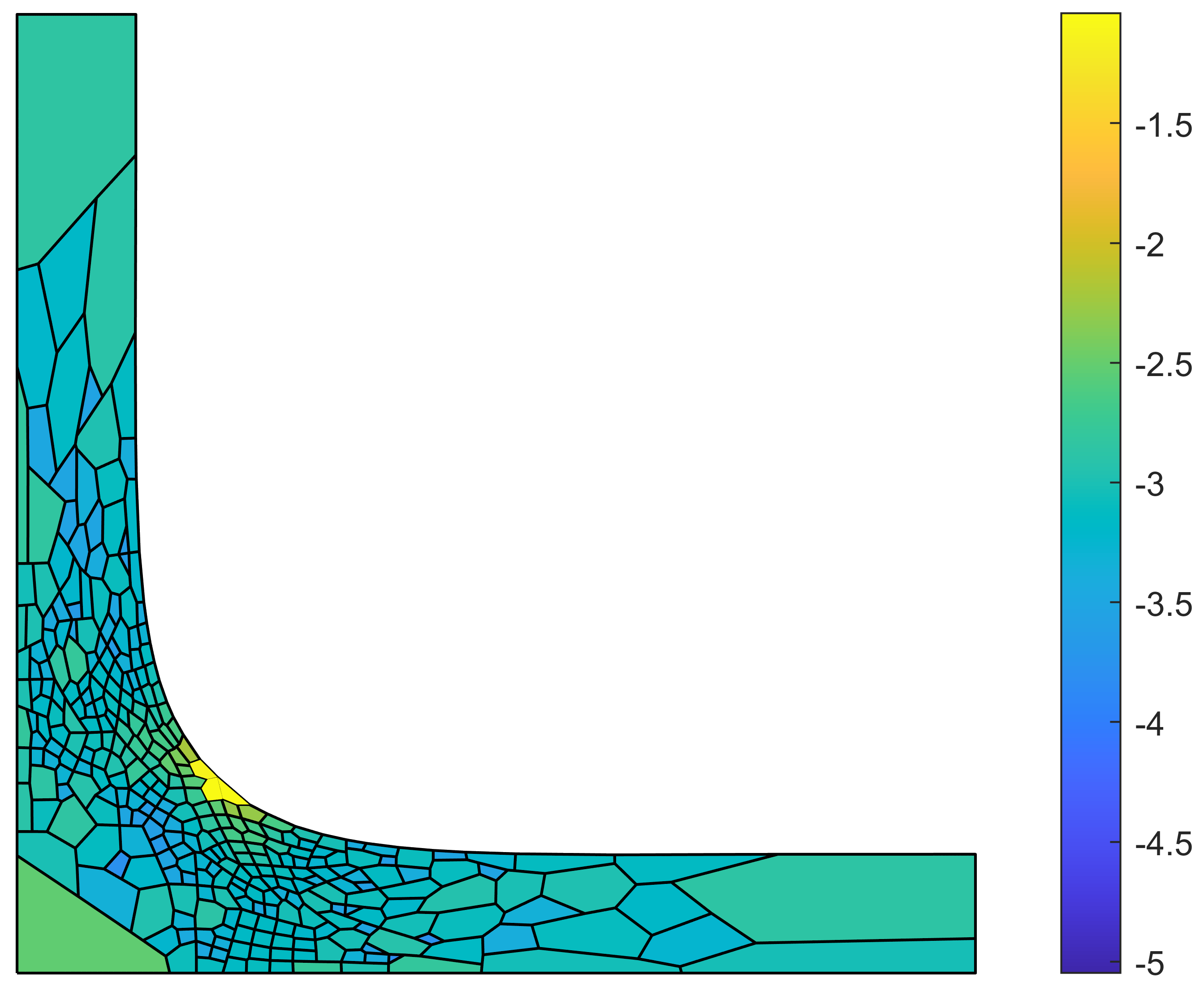}
			\caption{Displacement-based - Step 12}
		\end{subfigure}
		\vskip \baselineskip 
		\begin{subfigure}[t]{0.33\textwidth}
			\centering
			\includegraphics[width=0.95\textwidth]{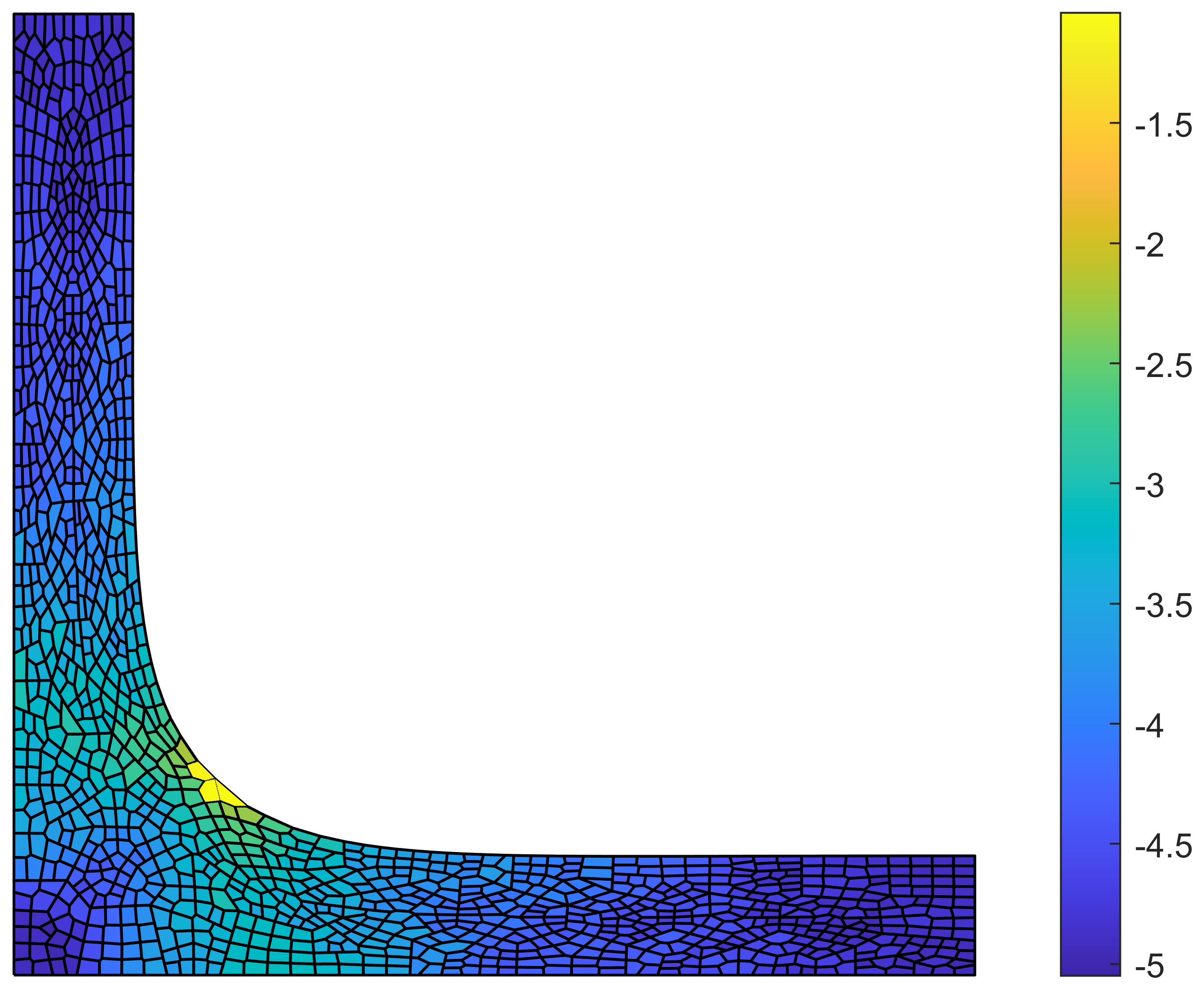}
			\caption{Energy error-based - Step 1}
		\end{subfigure}%
		\begin{subfigure}[t]{0.33\textwidth}
			\centering
			\includegraphics[width=0.95\textwidth]{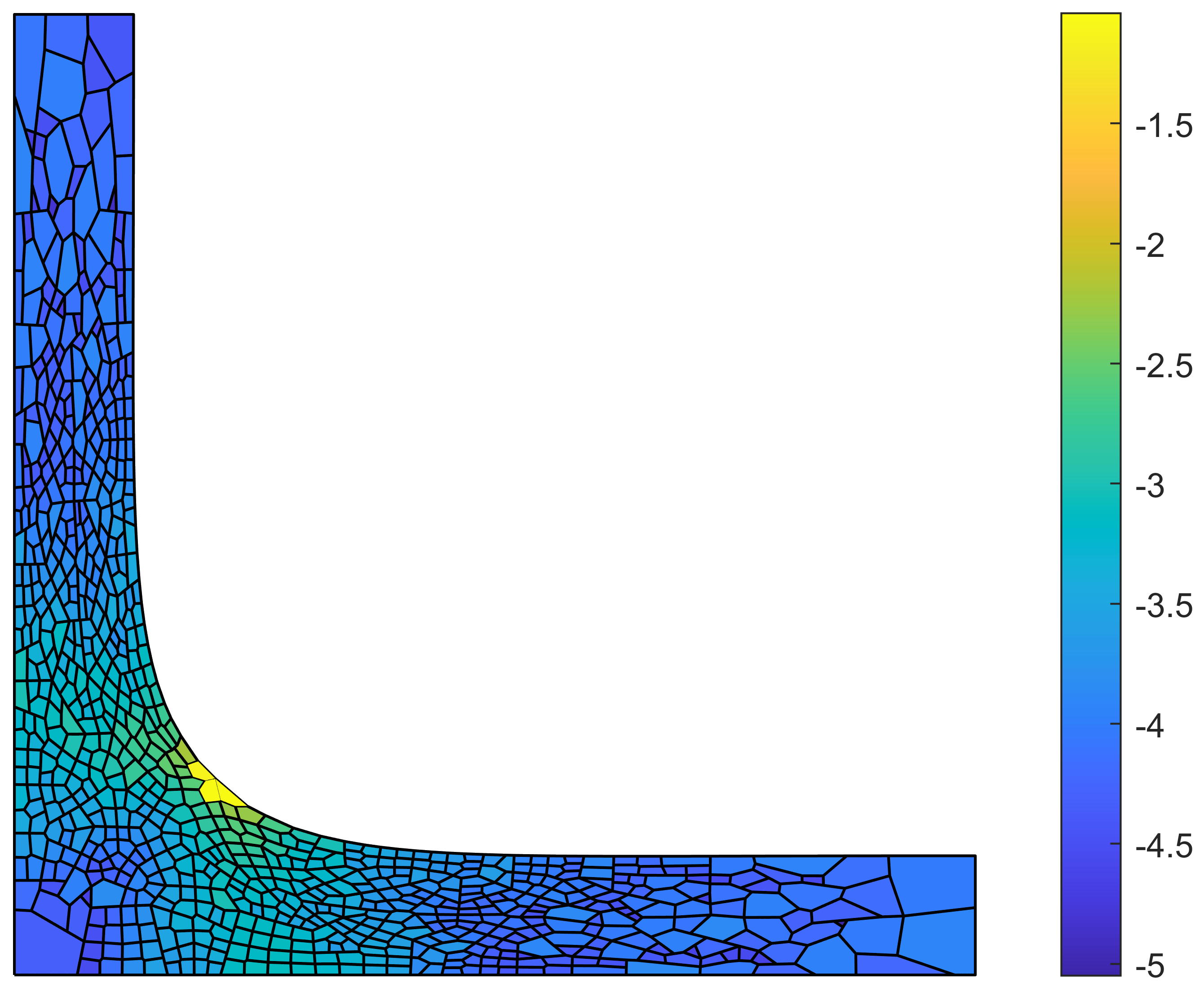}
			\caption{Energy error-based - Step 5}
		\end{subfigure}%
		\begin{subfigure}[t]{0.33\textwidth}
			\centering
			\includegraphics[width=0.95\textwidth]{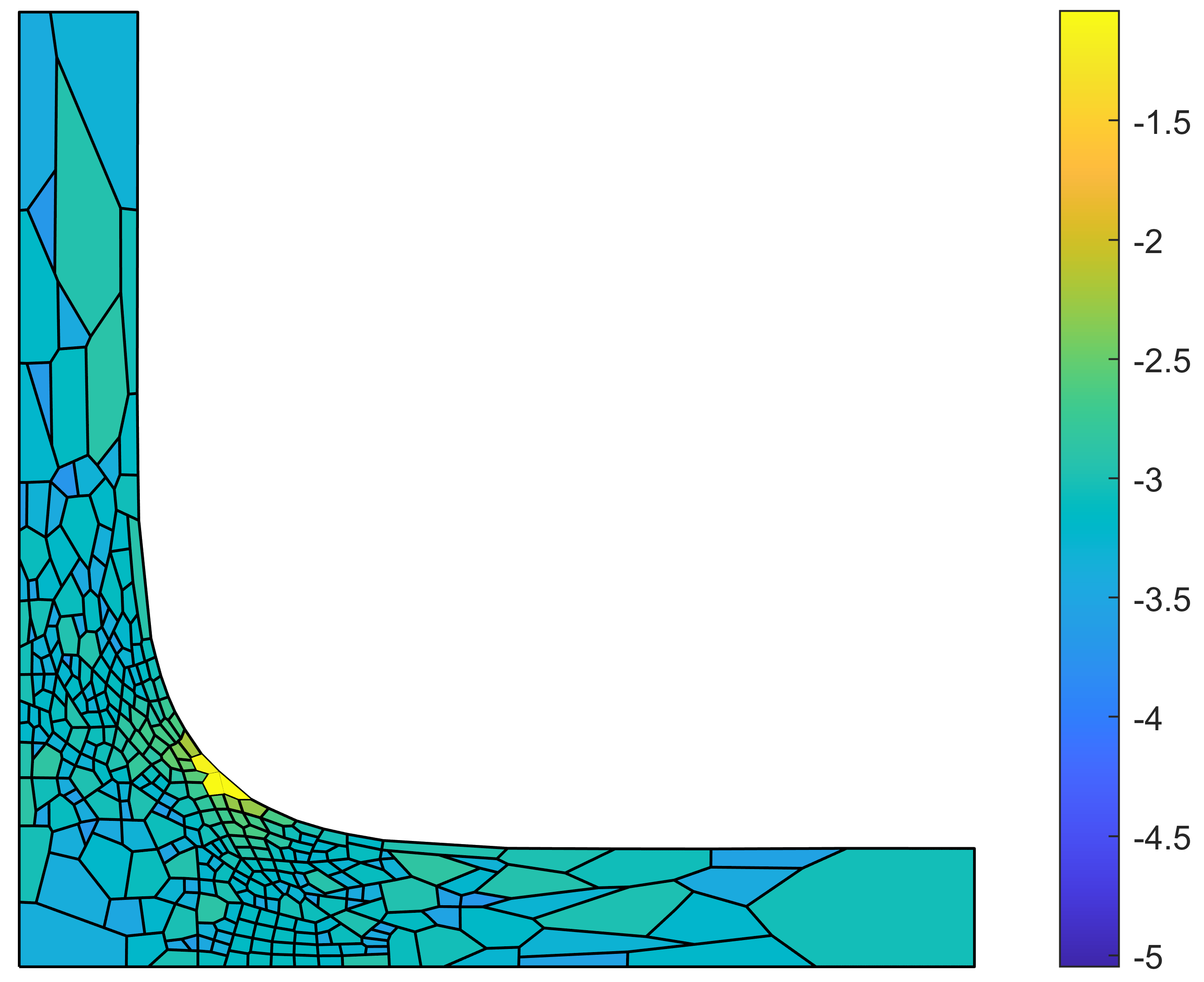}
			\caption{Energy error-based - Step 12}
		\end{subfigure}
		\caption{${\mathcal{H}^{1}}$ error distribution during the coarsening process for the L-shaped domain problem on Voronoi meshes using the displacement-based energy error-based coarsening procedure with ${T=20\%}$.
			\label{fig:LDomainErrorMeshes}}
	\end{figure}

	The distribution of the ${\mathcal{H}^{1}}$ error over the domain during the mesh coarsening process for the L-shaped domain problem is depicted in Figure~\ref{fig:LDomainErrorMeshes}. The ${\mathcal{H}^{1}}$ error is depicted in a logarithmic scale on Voronoi meshes for the cases of the displacement-based and energy error-based coarsening procedures with ${T=20\%}$. It is, again, clear from the error distribution in step~1, i.e. Figure~\ref{fig:LDomainErrorMeshes}(a), that the mesh evolution illustrated in Figure~\ref{fig:LDomainMeshes} is sensible.
	As seen in the previous examples, the error distribution over the domain becomes increasingly even during the mesh coarsening process. However, in this `very challenging´ problem the singularity at the corner of the L is a clear error hotspot. In problems containing singularities it is not possible to create a perfectly uniform error distribution. Nevertheless, the improved error distribution over the rest of the problem domain demonstrates the efficacy of the proposed coarsening procedures.
	
	\FloatBarrier
	\FloatBarrier
	
\begin{figure}[ht!]
		\centering
		\begin{subfigure}[t]{0.5\textwidth}
			\centering
			\includegraphics[width=0.95\textwidth]{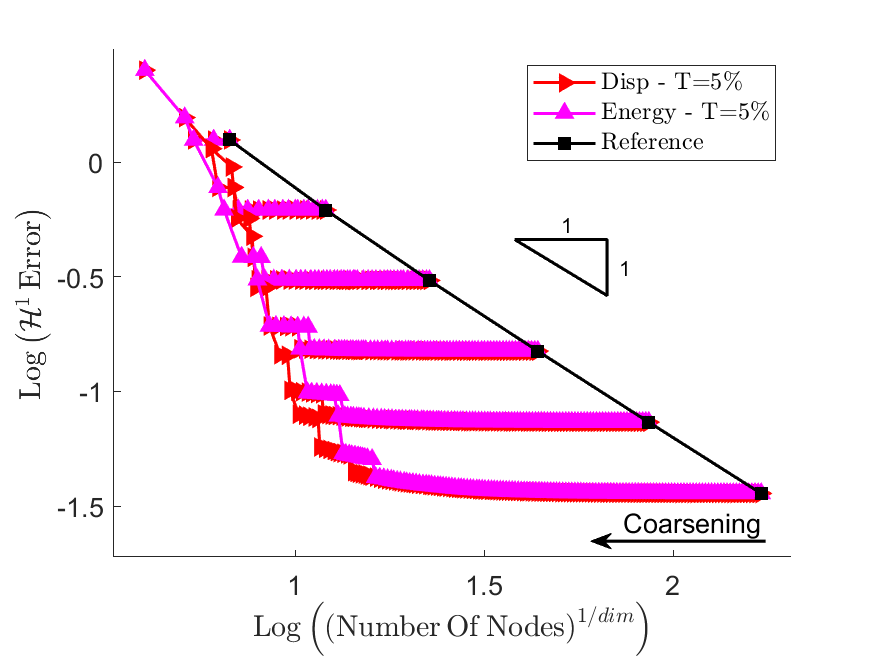}
			\caption{Structured meshes}
		\end{subfigure}%
		\begin{subfigure}[t]{0.5\textwidth}
			\centering
			\includegraphics[width=0.95\textwidth]{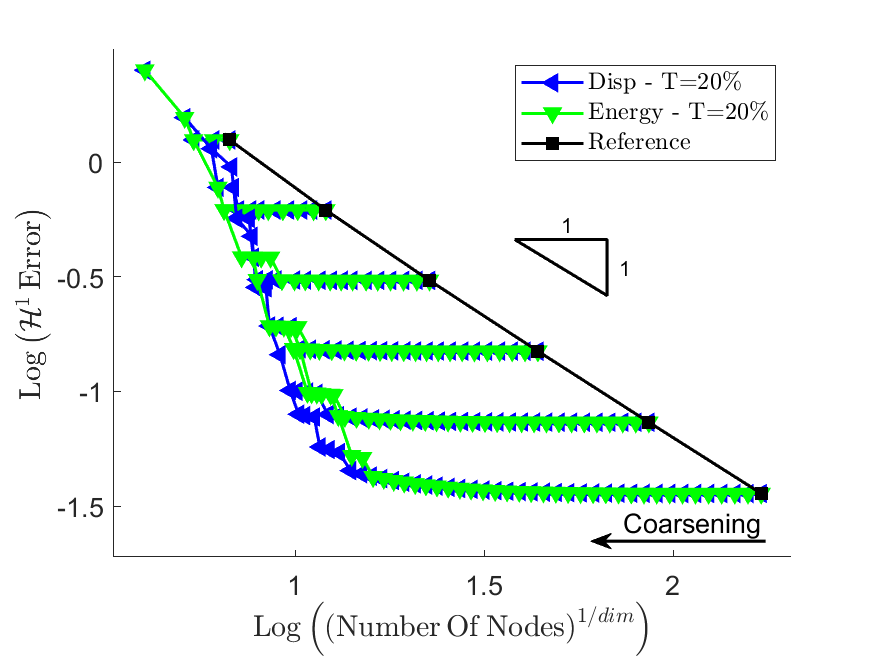}
			\caption{Structured meshes}
		\end{subfigure}
		\vskip \baselineskip 
		\begin{subfigure}[t]{0.5\textwidth}
			\centering
			\includegraphics[width=0.95\textwidth]{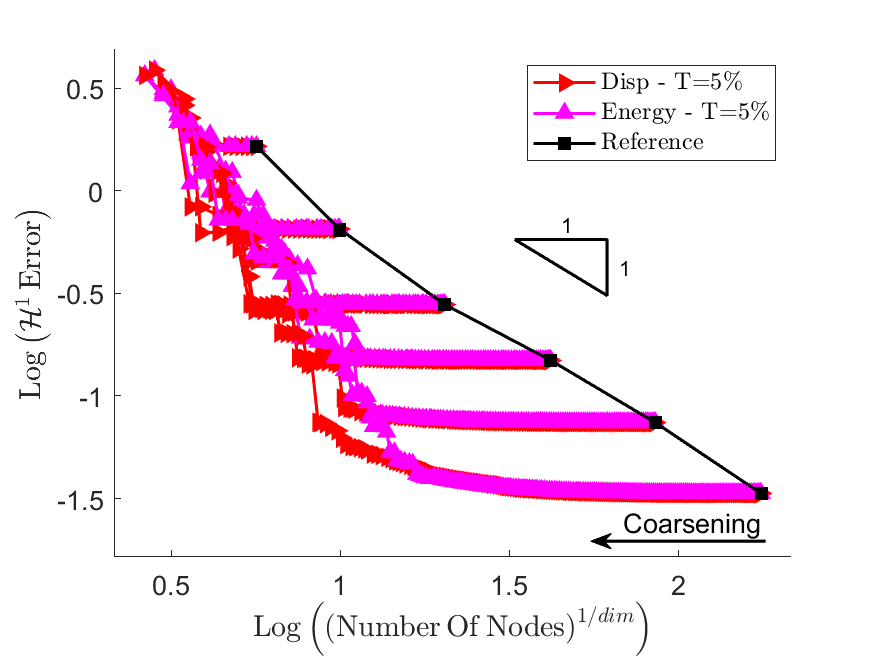}
			\caption{Voronoi meshes}
		\end{subfigure}%
		\begin{subfigure}[t]{0.5\textwidth}
			\centering
			\includegraphics[width=0.95\textwidth]{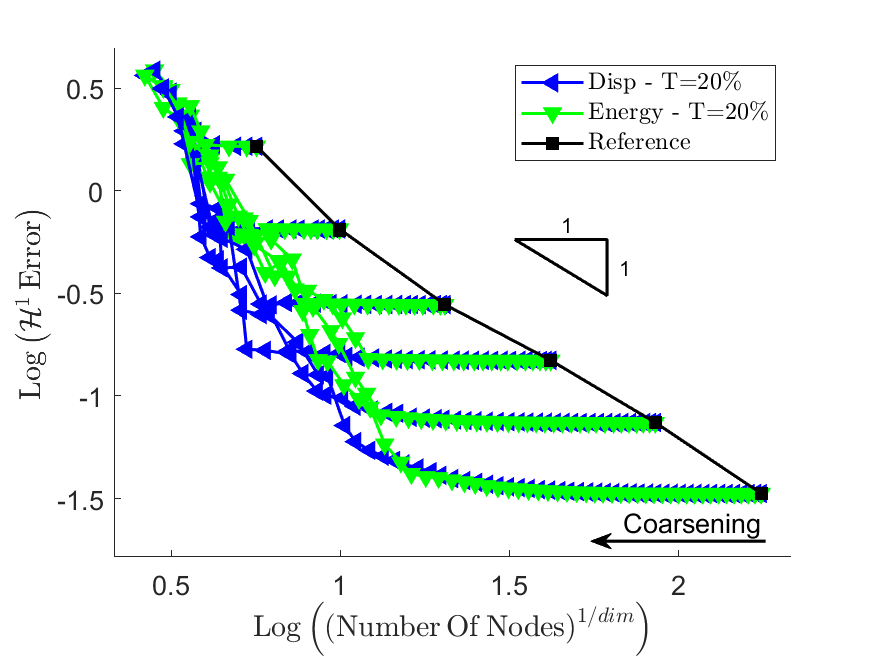}
			\caption{Voronoi meshes}
		\end{subfigure}
		\caption{$\mathcal{H}^{1}$ error vs $n_{\rm v}$ for the L-shaped domain problem.
			\label{fig:LDomainConvergenceNumberOfNodes}}
	\end{figure}

	The convergence behaviour in the ${\mathcal{H}^{1}}$ error norm of the VEM for the L-shaped domain problem using the displacement-based and energy error-based coarsening procedures is depicted in Figure~\ref{fig:LDomainConvergenceNumberOfNodes} on a logarithmic scale. Here, the convergence behaviour of the displacement-based and energy error-based procedures are plotted on the same axis to comparatively asses the performance of the procedures. These comparisons are made for the cases of structured (top row) and Voronoi meshes (bottom row), and for the coarsening thresholds ${T=5\%}$ (left column) and ${T=20\%}$ (right column).
	For this `very challenging´ problem the proposed coarsening procedures exhibit a very high degree of efficacy on both structured and Voronoi meshes, and for both choices of $T$. The coarsening procedures eliminate a significant portion of the number of degrees of freedom while introducing a negligible amount of error. The coarsening procedures exhibit good performance even on coarse initial uniform meshes. 
	Similar levels of efficacy are exhibited by the displacement-based and energy-error based coarsening procedures. Differences are only evident in the very coarse mesh range where the displacement-based procedure slightly outperforms the energy error-based procedure. This is most likely because the accuracy of the recovered/super-convergent stresses, used in the computation of the energy error-based indicator, depend on the accuracy of the global solution. Thus, as the global accuracy decays in the very coarse mesh range the accuracy and efficacy of the energy error-based indicator decay too.
	
	The coarsening procedure exhibits very good performance for a large portion of the range until, quite suddenly, the error increases very rapidly. This is because the coarsening procedure ensures that patches of elements are coarsened at every step. For much of the range coarsening is focused at the ends of the L, i.e. far from the corner of the L. Eventually, it is no longer possible to coarsen these regions while preserving the domain geometry, which means that the area around the corner of the L is then coarsened. Since this area coincides with a singularity, the coarsening causes significant error increases. In practice, the coarsening procedure would be terminated before the most critical regions of a domain are coarsened. However, for investigative purposes, and to test the coarsening procedures as thoroughly as possible, it is chosen to run the coarsening procedure until no more coarsening is possible.
	Somewhat similar behaviour is also observed in the other example problems, however, it is the most evident/exaggerated in this example problem due to the strong singularity at the corner of the L and the simple deformation in the rest of the domain.
	
	The results presented here, and in the other example problems, demonstrate the power of the proposed coarsening procedures in improving the efficiency of a simulation by removing the least critical elements and degrees of freedom while maintaining a high level of solution accuracy. Furthermore, it has been found that the efficacy of the coarsening procedures is not particularly sensitive to the choice of the coarsening threshold $T$. Thus, higher values of $T$ can be used to coarsen the mesh using fewer coarsening steps, and significantly faster, than lower values while yielding a similar degree of efficacy.
	
	\FloatBarrier
	\FloatBarrier
	
	\section{Discussion and conclusion} 
	\label{sec:Conclusion}
	In this work two novel mesh coarsening indicators have been proposed that are suitable for virtual element applications. Additionally, a simple procedure for selecting patches of elements qualifying for coarsening has been presented along with a novel mesh coarsening procedure that is suitable for both structured and unstructured/Voronoi meshes. 
	
	The proposed displacement-based and energy error-based coarsening indicators are computed over patches of elements and are motivated respectively by trying to; identify groups/patches of elements over which the displacement field is approximately linear, and predict the approximate error in the energy norm that would result from the coarsening of a particular patch of elements. 
	The mesh coarsening procedure involves constructing the geometry of a coarsened element by creating a convex hull around the patch of elements identified for coarsening. The geometry of the convex hull is created using a novel edge straightening procedure.
	
	The proposed mesh coarsening procedures were studied numerically over a range of benchmark problems of varying complexity. 
	For each problem the mesh evolution during coarsening was analysed along with the distribution of error over the problem domain.
	In terms of the mesh evolution, the efficacy of the proposed coarsening procedures was evident in several ways. Firstly, it was observed that the mesh density was significantly reduced in what were identified as the least critical portions of the domains. Secondly, similar performance was demonstrated in the cases of both structured and unstructured/Voronoi meshes. Thirdly, it was demonstrated that the coarsening procedures produced a more even/uniform error distribution over a problem domain compared to an initial uniformly discretized mesh. Finally, it was noted that the meshes generated through the coarsening of initially fine uniform meshes were very similar to meshes generated via adaptive refinement of coarse uniform initial meshes on the same example problems. 
	Furthermore, the efficacy of the proposed mesh coarsening procedures was studied in the ${\mathcal{H}^{1}}$ error norm. Here, the performance of the procedures was investigated on structured and unstructured/Voronoi meshes and was compared to a reference approach comprising meshes of uniform discretization density. Additionally, the influence of the initial mesh density on the performance of the coarsening procedures was investigated.
	The numerical results demonstrated the high degree of efficacy of both proposed coarsening procedures. The procedures efficiently reduced the number of elements and degrees of freedom in the least critical portions of a problem domain, while preserving the fine discretization in the most critical regions of a domain. Through this process the coarsening procedures yielded significantly more efficient solutions than those generated on uniform meshes. That is, for a given number of degrees of freedom the solutions obtained from adaptively coarsened meshes are significantly more accurate than those obtained from uniform meshes.
	Furthermore, it was found that the coarsening procedure was able to improve the efficiency of a solution even when starting from very coarse initial uniform meshes.
	
	The good performance exhibited by the proposed adaptive coarsening procedures over a range of challenging example problems on both structured and unstructured/Voronoi meshes demonstrates its versatility, efficacy and suitability for application to elastic analyses using the virtual element method. 
	
	Future work of interest would be the extension to combined adaptive refinement and coarsening procedures. Additionally, the extension to non-linear problems, higher-order formulations, and problems in three-dimensions would be of great interest.
	
	\section*{Conflict of interest}
	
	The authors declare that they have no known competing financial interests or personal relationships that could have appeared to influence the work reported in this paper.
	
	\section*{Acknowledgements}
	
	This work was carried out with support from the German Science Foundation (DFG) and the National Scientific and Technical Research Council of Argentina (CONICET) through project number DFG 544/68-1 (431843479).
	The authors acknowledge with thanks this support.

	\bibliographystyle{elsarticle-num}
	\bibliography{VEM_References}

\end{document}